\documentclass[11pt]{article}

\usepackage{epic}
\usepackage{eepic}
\usepackage{amssymb}
\usepackage{color}

\title{{\bf Subfactor realisation\\ of modular invariants}}

\author{
{\sc David E.\ Evans}\\
{\footnotesize School of Mathematics, University of Wales, Cardiff,}\\
{\footnotesize PO Box 926, Senghennydd Road, Cardiff CF24 4YH, Wales, U.K.}\\
{\footnotesize e-mail: {\tt EvansDE@cf.ac.uk}}\\
{\sc Paulo R. Pinto}\\
{\footnotesize Department of Mathematics, Inst.\ Super.\ T\'ecnico,}\\
{\footnotesize Av. Rovisco Pais, 1049-001 Lisboa, Portugal}\\
{\footnotesize e-mail: {\tt ppinto@math.ist.utl.pt}} }

\date{13 February, 2003 and revised 2 April, 2003}

\begin{document}
 \maketitle

\newcommand\bproof {\noindent {\it Proof. }}
\newcommand\eproof {\hspace*{\fill}\nolinebreak\hspace*{\fill}
{$\Box$}\par\vspace{3mm}}
\newcommand\labl[1]{\label{#1}\ee\\[-9mm]\hspace*{\fill}\nolinebreak
                    \hspace*{\fill}{\small\tt {#1}}\par\vspace{1mm}\noindent}
\newcommand\lablth[1]{\label{#1}}
\newcommand\lablsec[1]{\label{#1}\hspace*{\fill}
                   {\tiny\sf Sec:}{\small\tt {#1}}}

\font\gotico=eufm10


\newcommand\erf[1] {Eq.\ \nolinebreak (\ref{#1})}
\newcommand\eerf[1] {Eq.\ \nolinebreak (#1)}
\newcommand\sect[1] {Sect.\ \nolinebreak \ref{#1}}
\newcommand\fig[1] {Fig.\ \nolinebreak \ref{#1}}

\def\bbC           {\Bbb{C}}
\def\bbM           {\Bbb{M}}
\def\bbN           {\Bbb{N}}
\def\bbNo          {\Bbb{N}_0}
\def\bbR           {\Bbb{R}}
\def\bbT           {\Bbb{T}}
\def\bbZ           {\Bbb{Z}}
\def\SUz           {{\mathit{SU}}(2)}
\def\ran           {\rangle}
\def\lan           {\langle}
\def\ch            {{\chi}}
\def\MXN           {{}_M {\mathcal{X}}_N}
\def\MXM           {{}_M {\mathcal{X}}_M}
\def\MXMa          {{}_M^{} {\mathcal{X}}_M^\a}
\def\MXMo          {{}_M^{} {\mathcal{X}}_M^0}
\def\MXMop         {{}_{M_+}^{}\!\! {\mathcal{X}}_{M_+}^0}
\def\MXMom         {{}_{M_-}^{}\!\! {\mathcal{X}}_{M_-}^0}
\def\MXMopm        {{}_{M_\pm}^{}\!\! {\mathcal{X}}_{M_\pm}^0}
\def\MXMp          {{}_M^{} {\mathcal{X}}_M^+}
\def\MXMm          {{}_M^{} {\mathcal{X}}_M^-}
\def\MXMpm         {{}_M^{} {\mathcal{X}}_M^\pm}
\def\MXMmp         {{}_M^{} {\mathcal{X}}_M^\mp}
\def\MXMpp         {{}_{M_+}^{} {\mathcal{X}}_{M_+}^+}
\def\MXMmm         {{}_{M_-}^{} {\mathcal{X}}_{M_-}^-}
\def\MXMppm        {{}_{M_+}^{} {\mathcal{X}}_{M_+}^\pm}
\def\MXMpmpm       {{}_{M_\pm}^{} {\mathcal{X}}_{M_\pm}^\pm}
\def\NXN           {{}_N {\mathcal{X}}_N}
\def\MaXMa          {{}_{M_a} {\mathcal{X}}_{M_a}}
\def\MaXMb           {{}_{M_a} {\mathcal{X}}_{M_b}}
\def\MbXMa           {{}_{M_b} {\mathcal{X}}_{M_a}}
\def\NXMc           {{}_{N} {\mathcal{X}}_{M_c}}
\def\McXN           {{}_{M_c} {\mathcal{X}}_{N}}
\def\MbXMb           {{}_{M_b} {\mathcal{X}}_{M_b}}
\def\NXM           {{}_N {\mathcal{X}}_M}
\def\Tr            {\hbox{Tr}}
\def\Det           {\hbox{Det}}
\def\a             {\alpha}
\def\la            {\lambda}
\def\tr            {{\mathrm{tr}}}

\def\per           {positive energy representation}
\newcommand\ddA[1] {$\mathrm{A}_{#1}$}
\newcommand\ddD[1] {$\mathrm{D}_{#1}$}
\newcommand\ddE[1] {$\mathrm{E}_{#1}$}


\newtheorem{definition}{Definition}[section]
\newtheorem{lemma}[definition]{Lemma}
\newtheorem{corollary}[definition]{Corollary}
\newtheorem{theorem}[definition]{Theorem}
\newtheorem{proposition}[definition]{Proposition}
\newtheorem{remark}[definition]{Remark}


\begin{center}
{\footnotesize \it Dedicated to Rudolf Haag on the occasion
of his eightieth birthday}
\end{center}

\begin{abstract}

We study the problem of realising modular invariants by braided
subfactors and the related problem of classifying nimreps. We
develop the fusion rule structure of these modular invariants.
This structure is useful tool in the analysis of modular data from
quantum double subfactors, particularly those of the double of
cyclic groups, the symmetric group on 3 letters and the double of
the subfactors with principal graph the extended Dynkin diagram
$D_5^{(1)}$. In particular for the double of $S_3$, 14 of the 48
modular modular invariants are nimless, and only 28 of the
remaining 34 nimble invariants can be realised by subfactors.
\end{abstract}


{\footnotesize
\tableofcontents
}


\section{Introduction}

A prominent problem in rational conformal field theory (RCFT) is
the classification of modular invariant partition functions
$Z(\tau)$=$\sum Z_{\la,\mu}\ch_\la(\tau)\ch_\mu^\ast(\tau)$ where
$\ch_\la(\tau)$=$\hbox{Tr}_{H_\la}(e^{2\pi i\tau (L_0-c/24)})$ is the
trace
in the irreducible representation $\la$ of the chiral algebra,
with conformal Hamiltonian $L_0$, $\mathrm{Im}(\tau)>0$ and $c$
is the central charge. This problem has been solved only for a
few models, although its mathematical formulation is very simple
in terms of the following modular data. For a given
finite-dimensional representation of the modular group
SL$(2,\bbZ)$, let $S=[S_{\la,\mu}]$ and $T=[T_{\la,\mu}]$ denote
the matrices representing the (images) of the generators
${0\,-1\choose 1\,\,\,0}$ and ${1\,1\choose 0\,1}$, respectively.
We further suppose that $T$ is diagonal, $S$ is symmetric, $S^2$
a permutation matrix and $S_{\la,0}\geq S_{0,0}> 0$ where ``0''
is a distinguished primary field. Then a coupling matrix $Z$ that
commutes with $S$ and $T$ subject to the constraints
$$Z_{\la,\mu}=0,1,2,3,\dots\ \hbox{ and } Z_{00}=1$$ is called a modular
invariant. These constraints reflect the physical background of
the problem. The condition $Z_{00}=1$ reflects the uniqueness of
the vacuum and will be relaxed in this work. Related to this
problem is the classification of non-negative integer matrix
representations (nimrep \cite[Page 735]{BEK2}) of
the original fusion rules. It arose in two a priori unrelated
contexts: by Cardy \cite{cardy} on boundary RCFT and di
Francesco-Zuber \cite{FZ} in an attempt to generalize the ADE
Coxeter graphs pattern. Lately, it has been further revived by
Ocneanu \cite{O}, Feng Xu \cite{xu}, B\"ockenhauer and Evans
\cite{BE1,BE2,BE3,BE4,BE5,BE6}, B\"ockenhauer, Evans and
Kawahigashi \cite{BEK1,BEK2,BEK3} and Evans \cite{E1,E2}. The
set of modular invariants for a given modular data is finite
\cite{Gan3, BE5} and has also been of intense study by Cappelli
et al. \cite{CIZ} and Gannon \cite{Gan2} for WZW models.
 The modular invariants are the
torus partition functions (one loop partition function of a
closed string) and the cylinder partition functions (one loop
partition function of an open string) are the nimreps of the
underlying RCFT. These two partition functions should be
compatible \cite{cardy}.
 We can find
modular data in a wide variety of contexts. In the sequel, we
will mainly be concerned with the ones arising from quantum
doubles of finite groups $G$ or those from the Wess-Zumino-Witten
(WZW) models for compact Lie groups $G$ at arbitrary levels
\cite{W}. For a fixed finite group $G$, the primary fields
are given by pairs $(a,\chi)$ where $a$ are representatives of
conjugacy classes of $G$ and $\chi$ are the characters of the
irreducible representations of the centraliser $C_G(a)$ of $a\in
G$. Then the matrices \cite{Gan}
\begin{eqnarray}
S_{(a,\chi),(a^\prime,\chi^\prime)}&=&{1\over |C_G(a)| |C_G(a^\prime)|}
\sum_{g\in G(a,a^\prime)}
\overline{\chi^\prime(g^{-1}ag)}\ \overline{\chi(ga^\prime g^{-1})},\nonumber\\
T_{(a,\chi),(a^\prime,\chi^\prime)}&=&\delta_{a,a^\prime}
\delta_{\chi,\chi^\prime}\chi(a)/\chi(e)
\label{st}
\end{eqnarray} where
$G(a,a^\prime)=\{g\in G| aga^\prime g^{-1}=ga^\prime g^{-1}a\}$
and $e$ denotes the identity of $G$ form the modular data of the
quantum double (subfactor). This modular data can be twisted for
every $[k]\in H^3(G,S^1)$ \cite{Gan}, and the subfactor
interpretation of this data is in \cite{iz1}. The parameter $k$ is
regarded as the level in this setting \cite{Gan}. Then the quantum
double of $G$ at level 0 is precisely the model in \erf{st}. In
the WZW setting (or loop group setting), the matrices $S$ and $T$
are the Kac-Peterson matrices constructed from the representation
theory of unitary integrable highest weights modules over affine
Lie algebras or in an exponentiated form from the positive energy
representations of loop groups \cite{W}. It has been observed
\cite{Gan} that the number of modular invariants of both the
twisted group and the WZW data are small, but usually dwarfed by
the swarm of modular invariants of the (untwisted) group models
(cf.\ the 48 modular invariants for the $G=S_3$ with the 9 for the
twisted ones).

In the framework of braided subfactors we consider a finite
system of endomorphisms $\NXN$ of a type III factor $N$ with a
(non-degenerate) braiding on it, which leads to a Verlinde fusion
rule algebra and produce the modular data $S$ and $T$
\cite{BEK1,R1}. Sources of modular data from braided subfactors
arise from the WZW models and Ocneanu asymptotic subfactor which
is regarded as the subfactor analogue of Drinfel'd quantum double
construction. This {\it quantum double subfactor} is basically
the same as the Longo-Rehren inclusion \cite{LR} and is a way of
yielding braided systems from not necessarily commutative
systems.  The $S$- and $T$- matrices for the WZW subfactor models
have been constructed by A. Wassermann \cite{W} proving that they
indeed coincide with the Kac-Peterson $S$ and $T$ matrices. The
modular data from a quantum double subfactor was first
established by Ocneanu \cite[Section 12.6]{EK} using topological
insight, later by Izumi \cite{iz1} with an algebraic flavour (see
also \cite{sato1,sato2}).

Fixing a braided system of endomorphisms on a type III factor $N$,
we look for inclusions $\iota:N\hookrightarrow M$ such that its
dual canonical endomorphism $\theta=\bar{\iota}\iota$ decomposes
as a sum of endomorphisms from $\NXN$. To produce a modular
invariant from such an inclusion, we first employ Longo-Rehren
$\a^\pm$-induction method \cite{LR} of extending endomorphisms of
$N$ to  those in $M$ and then compute the dimensions of the
intertwiner spaces $Z_{\la,\mu}:=\lan\a_\la^+,\a_\mu^-\ran$. The
matrix $Z_{N\subset M}=[Z_{\la,\mu}]$ thus constructed from a
braided inclusion $N\subset M$ is a modular invariant
\cite{BEK1,E2}, see also Sect.\ \ref{frobsub} or Theorem
\ref{nice}. A modular invariant $Z_{N\subset M}$ encodes the rich
structure of the inclusion $N\subset M$. The induced systems,
which may well be non-commutative and hence not braided, encode
the quantum symmetries of the physical situation we started with.
Given the list of all modular invariants of a $S$- and $T$-model
it is an interesting problem to determine which can be produced
by braided inclusions through $\a$-induction. This is in general
a difficult task and in turn related to the problem of
classification of subfactors. The present work features the
several known  constraints a modular invariant has to fulfill in order
to be produced by an inclusion as well as new ones
arising from the fusion rule structure of those realised by subfactors.

In the course of previous work \cite{BE5, E1,EP2}, it was observed that modular
invariants satisfy remarkable fusion rules. If the sufferable modular
invariants $\{Z_{a}; a\}$  which can be realised by subfactors span the
commutant $\{S,T\}$ of the modular representation, then are certainly
decompositions
 $Z_{a}Z_{b} = \Sigma_{c}  {m_{ab}}^{c} Z_{c}$
 where the ${m_{ab}}^{c}$ are some complex coefficents.
 Certainly the product  $Z_{a}Z_{b}$ commutes with S and T and
 so is a modular invariant with non-negative integer entries
 but may not be normalised. It is therefore intriguing to ask
 whether any non-normalised modular invariant can be decomposed
 in this way with non-negative coefficients, in particular
 such products $Z_{a}Z_{b}$. It is in fact more natural to consider
 \begin{eqnarray}
 Z_{a}Z_{b}^{t}= \Sigma_{c}  {n_{ab}}^{c} Z_{c}
 \end{eqnarray}
 where $t$ is the transpose. Such decompositions were  noted by
 one of us in the context of $\SUz$ modular invariants \cite{BE5}
 particularly at level $16$. When $N \subset M$ is a braided subfactor
 realising a modular invariant $Z$, then the induced $M$-$M$ system
 has $\sum_{\lambda, \mu} Z_{\lambda,\mu}^2$ irreducible sectors, 
whilst the $N$-$M$ system has
 $\sum_{\lambda} Z_{\lambda,\lambda}$ 
irreducible sectors yielding the Cappelli-Itzykson-Zuber (CIZ)
 graph and corresponding nimrep in the case of $\SUz$ \cite{BEK2}.
 In the case of $Z= Z_{\rm E_{7}^Ë}$, the modular invariant associated to
 the Dynkin diagram $E_{7}$, the $M$-$M$ system thus has $17$ sectors which
 naturally decomposes into two induced $N$-$N$ orbits, on
 $\rm D_{10}$ and $\rm E_{7}$ graphs, $17 = 10 +  7$. This was the
 original motivation in \cite{BE5, BE6}
to write the numerical count $\sum_{\la, \mu} Z_{\la,\mu}^2$
as $\Tr{ZZ^t}$, the trace of a modular invariant, since
the full system $\MXM$ for  $Z= Z_{\rm E_{7}^Ë}$ now breaks up as
$\Tr Z_{\rm E_7}^2$= $\sum_{\la, \mu} Z_{\la,\mu}^2$ = $17 = 10 + 7$
= $\Tr Z_{\rm D_{10}}+ \Tr Z_{\rm E_7}$. Remarkably
not only is this numerical trace computation  consistent
but the matrices themselves satisfy
$Z_{\rm E_7}^2=Z_{\rm D_{10}}+Z_{\rm E_7}$.
This lead us to consider  the role of these
 modular invariant fusion rules and the role of the matrix
 $ZZ^{t}$ which is nonnormalised modular invariant and to understand
 the full graph of $Z$ as the CIZ graph for $ZZ^{t}$. These ideas
 were used as a powerful tool
 in understanding the full $M$-$M$ graph in terms of the $N$-$N$ orbits
 in many examples -
 see e.g. \cite[Section 6]{BE5}.
 In \cite{E1}, this analysis was taken further and in particular the
 programme to analyse and understand these fusion rules of modular
 invariants and the $M_{a}$-$M_{b}$ system in terms of the
 decomposition of the modular invariant $Z_{a}Z_{b}^{t}$
 into normalised modular invariants was pushed further.
 Here we complete this programme. One of the  tools we use is
 adopted from the work of \cite{ost1, frs} in the setting of
 Frobenius algebras which in turn had adapted
 the work on $\alpha$-induction
 \cite{BE1,BE2,BE3,BE4,BE5,BEK1,BEK2} and used it as a tool
 in a braided tensor category framework.
  The braided Frobenius algebra product of \cite{ost1}
  in the factor inclusion  context allows us to construct
  inclusions whose canonical endomorhpisms are the products of
  individual canonical endomorphisms. Then the central decomposition
  of the associated von Neumann algebra inclusion into subfactors yields the
  fusion rule algebra structure Theorem \ref{nice}.

We note that it is known for
sometime that some reasonable looking partition functions are
nevertheless unphysical because they cannot appear in any
consistent RCFT. With the braided inclusion approach we have an
efficient machine to see if a given partition function is in fact
physical.
  There is considerable evidence that the sufferable modular invariants
are precisely those of physical interest.  In the setting of algebraic
quantum field theory \cite{BE4}, if the given factor $N$
extends to a quantum field theoretical net of factors $\{N(I)\}$
indexed by proper intervals $I\subset \bbR$ and the system
$\NXN$ is obtained as restricting DHR-morphisms, then a braided
subfactor  $N \subset M$ will provide two local nets
of subfactors $\{N(I)\subset M_\pm(I) \}$ which provide
a rigorous formulation of the left and right maximal extensions
of the chiral algebra.
Indeed
 Rehren has shown \cite{R7}
 with chiral observables as light-cone nets
built in an observable net over 2D Minkowski space
that any braided extension $N\subset M$
determines an entire local 2D conformal field theory over
Minkowski space. The vacuum Hilbert space
of the 2D net decomposes upon restriction to the tensor
product of the left and right chiral observables precisely according
to the matrix $Z$
arising from $N\subset M$ through $\a$-induction.
In the algebraic quantum field theory
framework, $\alpha$-induction plays a critical role in taking
localised DHR endorphisms to solitonic ones but the
the neutral system $\MXMo$ however corresponds to proper DHR
endomorphisms.
 Braided subfactors provide consistent unitary 6$j$-symbols
or Frobenius *-algebras in braided tensor categories which
provide the data for computation of correlation functions
\cite[Section 5.3]{pz}, \cite{frs}. The 6$j$-symbol or connection
approach and bimodules is more amenable to the type II setting and statistical
mechanical framework whilst the Frobenius *-algebras or Q-systems
and sectors
are more amenable to the type III setting and conformal field theory
framework as explained in \cite{EP2}. The $N$-$M$ sectors (which is
the CIZ graph in the $\SUz$ setting) describes the boundary conditions, whilst
the full system of $M$-$M$ sectors describe the defect lines.

In this work we start the analysis of modular invariants from the
modular data arising from quantum double subfactors. For that we
study further the structure of fusion rules of modular invariants
\cite{BE6,E1,frs,EP2}. The first obvious case is the quantum
double of the finite group subfactor $M_{0}\subset M_{{0}} \rtimes
G$, where the finite group G acts  outerly on a type III factor
$M_{0}$ identified by Ocneanu and later by Izumi (see \cite{EK})
to be the group-subgroup subfactor $N=M_0\rtimes\Delta(G)\subset
M_0\rtimes{(G\times G)}=M$ where $\Delta(G)=\{(g,g):g\in G\}$
denotes the diagonal subgroup of $G\times G$. The case $G=\bbZ_2$
was studied by B\"ockenhauer-Evans \cite{BE4} (in the context of
the ${{\mathit {SO}}}(16n)$ WZW level 1 models). Here we consider
the next models $G=\bbZ_3$ and $S_3$ the cyclic group of order 3
and the symmetric group on 3 elements. Some obvious analogues in
\cite{E1} from the previously studied models (most notably the
$\SUz$ WZW level $k$ model by B\"ockenhauer, Evans and Kawahigashi
\cite{BEK2} and Ocneanu \cite{O} where all the normalised modular
invariants are indeed produced by subfactors) are no longer true
(Remark \ref{canformZ}). It has been announced by Ocneanu
\cite{ocn4} that all ${\mathit{SU}}(3)$ WZW modular invariants are
produced by subfactors. An outline of the history of WZW modular
invariants is in \cite{E2}. The quantum $S_3$ double model has 48
modular invariants. The list of 32 modular invariants announced in
\cite{Gan} is not complete, we found the complete list with 16
other invariants. This extended list has been confirmed to us by
T. Gannon \cite{tg}. The quantum $S_3$ double model shows a rich
and complex structure. In order to have a thorough study of this
model we were led to consider products of modular invariants
\cite{BE6,E1,frs} and in fact showing that the product of dual
canonical endomorphisms are still canonical endomorphisms in the
setting of braided systems whose modular invariant is the
matricial product (\sect{frobsub}). Using module category theory,
Ostrik \cite{ost} has computed the possible $\Tr(Z)$ and
$\Tr(ZZ^t)$ for the list of modular invariants $Z$ from the
quantum $S_3$ double model which arise from module categories of
the modular category, all labelled $(K,\psi)$ for $K$ subgroups of
$S_3\times S_3$ and $\psi$ are elements of the 2-cohomology group
$H^2(K, S^{1})$. With our approach we found module categories (of
the $\MXN$ system) giving rise to $\Tr(Z)=4$ which were missing in
\cite{ost} since the modular invariants realised by subfactors
must form a fusion rule algebra (see Theorem \ref{nice}). Indeed
there are two subgroups (conjugacy classes) of $S_3\times S_3$
producing trace 4 matrices.

This paper is organized as follows. In \sect{pre} we review some
mathematical tools, see e.g. \cite{BEK1}, necessary in the
sequel. It contains some elements of the nimrep theory following
\cite{Gan2} and on the subfactor side the paper \cite{BEK2} and
its generalizations. In this section we also propose some
terminology for modular invariants as in \sect{mixing}. In
\sect{frobsub}, we provide the fusion rule decomposition
of modular invariants using the central decomposition of
braided inclusions. In \sect{kosakietal}, the
relation between the subgroups of a finite group $G$ and the
intermediate subfactors $N=M_0\rtimes\Delta(G)\subset M_0\rtimes
(G\times G)=M$ is clarified. In \sect{cyclicprime} we prove that
all the 8 quantum $\bbZ_3$ double modular invariants are realised
by subfactors. \sect{s3} is devoted to the quantum $S_3$ double
model. With a numerical search and employing the estimate of
\cite[\eerf{1.6}]{BE6}, we found all the 48 quantum $S_3$ double
modular invariants and then the Verlinde fusion matrices. We
classify the modular invariants that have matching nimreps and
decide which modular invariants are indeed produced by braided
subfactors (we list the nimble and the sufferable ones in
Corollary \ref{list}). In \sect{pictures}, the full systems
arising from the sufferable modular invariants are
displayed.

\section{Preliminaries}
\label{pre}

In this section we recall
the general framework of \cite{BEK1,BEK2,Gan}.

\subsection{Braided inclusion theory}

We shall consider type III von Neumann algebras with finite dimensional
centres. A morphism between such algebras $A$ and $B$
shall be a faithful unital $\ast$-homomorphism $\rho:A\rightarrow B$,
called a $B$-$A$ morphism, and we write $\rho\in\hbox{Mor}(A,B)$.
We will only consider those of finite statistical dimension
or inclusions of finite index.
Then $d_\rho=[B:\rho(A)]^{1/2}$ is called the statistical (or
quantum) dimension of $\rho$; here $[B:\rho(A)]$ is the Jones
index of the inclusion $\rho(A)\subset B$. If $\rho$ and $\sigma$
are $B$-$A$ morphisms with finite statistical
dimensions, then the vector space of intertwiners
$\hbox{Hom}(\rho,\sigma)$=$\{t\in B: t\rho(a)=\sigma(a)t\,,\,\,
a\in A \}$ is finite-dimensional, and we denote its dimension by
$\lan\rho,\sigma\ran$. A morphism conjugate to $\rho$
will be denoted $\bar{\rho}: B\hookrightarrow A$.

Consider a type III inclusion $N\subset M$ with $N$
a factor, and denote by $\iota$ the inclusion map.
Then $\gamma=\iota\bar{\iota}$
and $\theta=\bar{\iota}\iota$ are the canonical and dual canonical
endomorphism of $N\subset M$ respectively. In the group case,
$M^G\subset M$ where we have a outer action $\a$ of the finite
group $G$ on say a factor $M$, then $\gamma$ decomposes as a
sector into the group elements, $[\gamma]=\oplus_g[\a_g]$, the
decomposition of $\theta$ as a sector is according to irreducible
representations of $G$,
$[\theta]=\oplus_{\pi\in\hat{G}}d_\pi[\pi]$ with $d_\pi$ being
the dimension of the irreducible representation of $\pi$.

Let $\NXN$ denote a finite system of irreducible endomorphisms of
a  factor $N$ in the sense that different elements of $\NXN$ are
inequivalent, for any $\la\in\NXN$ there is a representative
$\bar{\la}\in\NXN$ of the conjugate sector $[\bar{\la}]$ and
$\NXN$ is closed under composition and subsequent irreducible
decomposition. We denote by $\Sigma(\NXN)$ a set of
representative endomorphisms of integral sums of sectors from
$\NXN$ \cite{BEK1}. The quantity $\omega=\sum d_\la^2$ is
the global index of the system.

If $\NXN=\{\rho_\xi\}$ is a finite system and $j:  N\to
N^{\mathrm{opp}}$ is the natural anti-linear isomorphism,
$\rho_\xi^{\mathrm{opp}}=j\cdot \rho_\xi\cdot j$, we set
$B=N\otimes N^{\mathrm{opp}}$. Longo and Rehren \cite{LR} have
shown that there exists a subfactor $A\subset B$ such that
$\gamma=\bigoplus_\xi\ \rho_\xi\otimes\rho_\xi^{\mathrm{opp}}$ is
its canonical endomorphism (now referred to as the Longo-Rehren
inclusion). This notion was further translated into the type II$_1$
setting by Masuda \cite{masuda} who proved that Longo-Rehren
inclusion and Ocneanu asymptotic inclusion are essentially the
same object. Ocneanu \cite[Chapter 12]{EK} has constructed a
non-degenerate braiding on the $A$-$A$ system from the above
quantum double inclusion $A\subset B$ and later in a more algebraic
way by Izumi \cite{iz1, iz2}. For a subsystem $\Pi$ of $\NXN$
Izumi's Galois correspondence \cite[Proposition 2.4]{iz1} asserts
then that there exists an intermediate subfactor $A\subset
B_\Pi\subset B$ such that $\gamma_\Pi=
\bigoplus_{\xi\in\Pi}\rho_\xi\otimes\rho_\xi^{\mathrm{opp}}$ is
a canonical endomorphism of the inclusion $B_\Pi\subset B$.


Whenever we have a non-degenerate braiding on a system, from the
Hopf link and twist we can define $S$- and $T$-matrices of type
$\NXN\times\NXN$ \cite{R1,tur}, which satisfy the Verlinde formula
\cite{Ve},
\begin{eqnarray}
\sum_{\rho\in\NXN}{S_{\la,\rho}S_{\mu,\rho}S_{\nu,\rho}^\ast\over
S_{\rho,0}}=
\lan\la\mu,\nu\ran.
\end{eqnarray}
The fusion matrices $N_\la=[N_{\la,\mu}^\nu]_{\mu,\nu}$, where $
N_{\la,\mu}^\nu=\lan\la\mu,\nu\ran$, recover the original fusion
rules, i.e.
$$N_\la\cdot N_\mu=\sum_{\nu\in\NXN} N_{\la,\mu}^\nu N_\nu.$$


In our setting of a braided inclusion $\iota:N\hookrightarrow M$
of type III we are interested to extend a morphism
$\la\in\hbox{End}(N)$ to another morphism in $\hbox{End}(M)$. We
assume now that we have a (type III) inclusion
$\iota:N\hookrightarrow M$ together with a finite system
$\NXN\subset\hbox{End}(N)$ which is non-degenerately braided and
such that the dual canonical endomorphism
$\theta=\bar{\iota}\iota\in\Sigma(\NXN)$ for the inclusion $M$-$N$
morphism $\iota:N\hookrightarrow M$. In the above setting, we say
that $N\subset M$ is a braided subfactor. One can define the
$\a$-induced morphisms $\a_\la^\pm\in\hbox{End}(M)$ by the
Longo-Rehren formula \cite{LR}: $\a_\la^\pm= \bar{\iota}^{\
-1}\circ\hbox{Ad}(\varepsilon_{(\la,\theta)}^\pm)
\circ\la\circ\bar{\iota}$, where $\varepsilon$ is the braiding, so
that $\a_\la^\pm$ extends $\la$ in $\NXN$,
$\a_\la^\pm\iota=\iota\la$. We can define the positive integral
matrix $Z_{\la,\mu}=\lan\a_\la^+,\a_\mu^- \ran,$ normalised at the
vacuum $Z_{0,0}=1$ if $M$ is a factor, sometimes denoted by
$Z_{N\subset M}$ when we want to emphasize the inclusion from
which it was constructed. By \cite{BE4,E2}, the matrix
$Z_{N\subset M}$ commutes with the modular $S$- and $T$-matrices
for subfactors which also holds for inclusions under the
decomposition into normalised ones (cf.\ Theorem \ref{nice}).
Therefore $Z$ is a modular invariant. Now we use $\a$-induction
and the inclusion map $\iota$ to construct finite systems whose
general theory has been developed in \cite{BEK1,BEK2}. Let us
choose representative endomorphisms of each irreducible subsector
of sectors of the form $[\iota\la\bar{\iota}]$, $\la\in\NXN$. Any
subsector of $[\a_\la^+\a_\mu^-]$ is automatically a subsector of
$[\iota\nu\bar{\iota}]$ for some $\nu$ in $\NXN$ and since we
assume the non-degeneracy of the braiding the converse also holds
true \cite{BEK1}. This set of sectors yields a system $\MXM$ of
sectors in general non-commutative (the original sectors from the
system $\NXN$ is commutative since it is braided). We define in a
similar fashion the chiral systems $\MXMpm$ to be the subsystems
of $\beta\in\MXM$ such that $[\beta]$ is an irreducible subsector
of $[\a_\la^\pm]$. The  neutral or ambichiral system $\MXMo$ is
defined as the intersection $\MXMm\cap\MXMp$, so that we obtain
$\MXMo\subset\MXMpm\subset\MXM$ (see e.g.\ \cite{BE5}).

In the braided subfactor case, their global indices (the sum of
the squares of the quantum dimensions) are denoted by
$\omega_0,\omega_\pm$ and $\omega$, and are completely encoded in
the modular invariant $Z$, namely \cite{BEK2, BE4},
\begin{eqnarray}
\omega_\pm={\omega\over \sum_{\la}d_\la Z_{\la,0}}=
{\omega\over \sum_\la d_\la Z_{0,\la}},\quad \omega_0=
\omega_\pm^2/ \omega.
\end{eqnarray}
To help find the irreducible sectors in each of the above system
one has the relation \cite[\eerf{33}]{BEK2}
\begin{eqnarray}\label{inequ1}
\lan\a_\la^\pm,\a_\mu^\pm\ran\leq\lan\theta\la,\mu\ran.
\end{eqnarray}
Finally, the neutral system $\MXMo$ inherits a
non-degenerate braiding (therefore their fusion rules are
commutative) by \cite{BE4}
whose modular matrices are denoted by $S^{\mathrm{ext}}$ and
$T^{\mathrm{ext}}$.
We can recover the matrix $Z$ from the branching coefficients
($b_{\tau,\la}^\pm=\lan\tau,\a_\la^\pm\ran$) as follows:

\begin{eqnarray}
Z_{\la,\mu}=\sum_{\tau\in\MXMo} b_{\tau,\la}^+ b_{\tau,\mu}^-.
\end{eqnarray}
Moreover we have the following intertwining properties of the
branching coefficient matrices between original and extended
 $S$- and $T$-matrices \cite[Theorem 6.5]{BE4}:
\begin{eqnarray}
S^{\mathrm{ext}}b^\pm=b^{\pm}S,\quad T^{\mathrm{ext}}b^\pm=b^{\pm}T.
\end{eqnarray}
The modular invariant $Z$ is a permutation matrix if and only if
$\MXMo=\MXM$ \cite{BEK1,BEK2}. The modular invariant $Z$ produced by a
subfactor $N\subset M$ is said to be of type I if
$Z_{0,\la}=\lan\theta,\la\ran$ for all $\la\in\NXN$ \cite{BEK2,BE4} (in
particular $Z$ is symmetric). In this situation we say that the
chiral locality holds for the subfactor $N\subset M$. Note
however that there may be other braided subfactors producing the
same modular matrix $Z$. In other words, the chiral locality
holds if the dual canonical endomorphism is visible in the vacuum
row (and hence column), so $[\theta]=\oplus Z_{0,\la}[\la]$.
In the
presence of chiral locality, we can usually recover the canonical
endomorphism from the full $M$-$M$ system \cite{BE3} by computing
the dimensions
$\lan\a_\la^+\a_{\bar{\mu}}^-,\gamma\ran=\lan\a_\la^+,\a_\mu^+\ran$.

If a braided subfactor $N\subset M$ produces a
modular invariant $Z$, then there are intermediate subfactors
$N\subset M_\pm\subset M$, such that $N\subset M_\pm$ produce
symmetric modular invariants $Z_\pm$ (type I parents) \cite{BE4}.
Moreover
we have $Z_{\la,0}=Z_{\la,0}^+$ and $Z_{0,\la}=Z_{0,\la}^-$ and
the canonical endomorphisms of $N\subset M_\pm$ are visible from
the vacuum row and column of $Z$: $[\theta_+]=\oplus_\la
Z_{\la,0}[\la]$ and $[\theta_-]=\oplus_\la Z_{0,\la}[\la]$, which
incidentally implies Eq.\ (\ref{inequ1}).

\begin{definition}
A non-negative integer representation (nimrep, for short
\cite{BEK2}) of a given modular data is an assignment of a matrix
$G_{[\la]}$ (or $G_{\la}$ for short) to each primary field $\la$,
with non-negative integer entries, preserving the fusion rules,
i.e.
$$G_\la\cdot G_\mu=
\sum_\eta N_{\la,\mu}^\eta G_\eta, \qquad G_0=\mathrm{id}, \hbox{
and } G_{\bar{\la}}=G_\la^t, \hbox{ for all } \la,\mu$$ where
$N_{\la,\mu}^\eta$ are the Verlinde fusion integers from the
original modular data.
\end{definition}

By a standard argument (see e.g. \cite{Ve} or \cite[Page
425]{EK}), we can simultaneously diagonalize all $G_\la$'s, and
moreover the eigenvalues of $G_\la$ are $S_{\mu\la}/S_{\mu 0}$ for
$\mu$ running in some multi-set (possibly with multiplicities).
The cardinality of this set is the dimension of the nimrep. Two
nimreps $G^\prime$ and $G^{\prime\prime}$ are equivalent if there
is a permutation matrix $P$ such that $PG_\la^\prime
P^{-1}=G_\la^{\prime\prime}$ for all $\la$. There is a natural
notion of direct sum of nimreps \cite{Gan2}.
The
{\it regular nimrep} is obtained, by setting $G_\la:=N_\la$ where $N_\la$
is the Verlinde fusion matrix associated to the primary
field $\la\in\NXN$.
By the exponents of a modular invariant $Z$,
we mean the multi-set Exp consisting of $Z_{\la \la}$
copies of $\la\in\NXN$.

RCFT is thought to require that each physical modular
invariant has at least one nimrep such that their
exponents coincide.
A modular invariant is {\it nimble} if it has a matching
nimrep. A nimrep hereafter is meant to be compatible
with some modular invariant (i.e., its spectrum and
dimension are dictated by the diagonal part of the
modular invariant). Hence we discard (although some
of the results below remain valid) nimreps with no
consistent modular invariant as being spurious nimreps.
A {\it nimless} modular invariant is a modular invariant
without a matching nimrep.
We regard the category of nimless modular invariants
the least interesting of all.

Following closely Gannon's paper \cite{Gan2}, let us fix a nimrep
$G$. Then by Perron-Frobenius theory $S_{\la 0}/S_{00}$ is the
norm of the matrix $G_\la$ for every $\la\in\NXN$. Clearly, all
$G_\la$'s are symmetric if all the exponents $\mu$ are self
conjugate $\mu=\bar{\mu}$. More generally, $G_\la=G_\eta$ if and
only if $S_{\la\mu}= S_{\eta\mu}$ for all $\mu\in\hbox{Exp}.$ We
conveniently enumerate some properties \cite{Gan2} of nimreps
that we use in the sequel.

\begin{enumerate}
\item For a given modular data, there are only finitely
many indecomposable inequivalent nimreps.

\item For every primary field $\la$, the norm of every
connected component of $G_\la$ is $S_{\la 0}/S_{00}$.
Moreover, ${[G_\la]}_{i,j}\leq S_{\la 0}/S_{00}$.

\item The spectrum of every matrix
$G_\la$ is $\{S_{\la\mu}/S_{0\mu}: \mu \in \hbox{Exp}\}$.

\item The number of indecomposable components of $G_\la$ is
the number of exponents $\mu$ such that $S_{\la\mu}/S_{0\mu}=
S_{\la 0}/S_{00}$.

\item No column or row of each $G_\la$ can be identically
equal to zero.
\end{enumerate}

We finish this section with a well posed problem.\\
{\bf Problem.} Classify the nimreps of a given modular
data.

\subsection{Induced nimreps for sufferable modular invariants}
\label{mixing}

\begin{figure}[htb]
\begin{center}
\unitlength 0.6mm

\begin{picture}(75,70)
\thinlines

\put(-10,30){\makebox(0,0){{\tiny $\bullet$}}}
\put(45,10){\makebox(0,0){{\tiny $\bullet$}}}
\put(45,50){\makebox(0,0){{\tiny $\bullet$}}}
\path(-10,30)(45,10)(45,50)(-10,30)

\put(-15,30){\ellipse{10}{3}}
\put(50,10){\ellipse{10}{3}}
\put(50,50){\ellipse{10}{3}}

\put(-10,27){\makebox(0,0){{\tiny $N$}}}
\put(48,14){\makebox(0,0){{\tiny $M_a$}}}
\put(46,54){\makebox(0,0){{\tiny $M_b$}}}

\put(60,30){\makebox(0,0){{\tiny
$\matrix{\#M_a$-$M_b\cr | |\cr\Tr(Z_a(Z_b)^t)}$}}}
\put(60,5){\makebox(0,0){{\tiny $\#M_a$-$M_a$=$\Tr(Z_a(Z_a)^t)$}}}
\put(11,44){\makebox(0,0){{\tiny $\#N$-$M_b$=$\Tr(Z_b)$}}}

\end{picture}
\end{center}
\caption{Irreducible $M_a$-$M_b$ sectors from the subfactors
$N\subset M_a$ and $N\subset M_b$ with modular invariants $Z_a$
and $Z_b$ } \label{MaMb}
\end{figure}
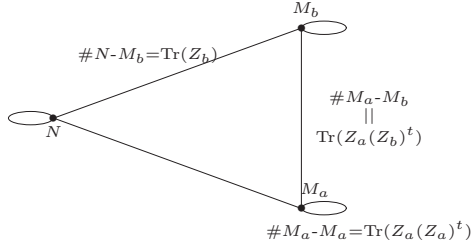

Let $N\subset M_a$ and $N\subset M_b$ be two braided subfactors,
with $\iota_a$ and $\iota_b$ denoting the inclusion maps
respectively, whose modular invariants of the given modular data
we denote by $Z_a$ and $Z_b$ respectively. Let $\MaXMb$ denote
the irreducible representative endomorphisms from the subsectors
of $[\iota_a\la\bar{\iota}_b]$ with $\la$ running over $\NXN$.
Then by extending ideas from \cite{BEK1} we can show \cite{E1}
that $\MaXMb$ under the left action $\MaXMa$ and right action of
$\MbXMb$ is isomorphic to
\begin{eqnarray}
\bigoplus_{\la,\mu\in\NXN}H_{\la,\mu}^a\otimes\overline{H}_{\la,\mu}^b
\label{iso}
\end{eqnarray}
where
$H_{\la,\mu}^c\subseteq\oplus_{x\in\NXMc}\mathrm{Hom}(\la\bar{\mu},x\bar{x})$
is a certain Hilbert space of dimension $[Z_c]_{\la,\mu}$,
$c=a,b$. In particular $\#\MaXMb=\Tr(Z_a(Z_b)^t)$. Moreover, we
have an induced nimrep in a stronger sense as follows (see
\cite{E1}).
The matrices
\begin{eqnarray}
\Gamma_{\la,\mu;\beta}^{\beta^\prime}=
\lan\beta\a_\la^+\a_\mu^-,\beta^\prime\ran,\qquad
\beta,\beta^\prime\in\MaXMb
\end{eqnarray}
form a nimrep (we can use either inductions $\a^{(a)}$ or $\a^{(b)}$ from
the braided subfactors $N\subset M_a$ or $N\subset M_b$
respectively). Such a $Z_a$-$Z_b$ nimrep provides us with further
constraints for a modular invariant to be produced by a braided
subfactor. The eigenvalues of $\Gamma_{\la,\mu}$ arise from
$\{\chi_\zeta(\la)\chi_\eta(\mu)\}$, where
$\chi_\la(\nu)=S_{\la,\nu}/S_{\la,0}$, with multiplicities
\begin{eqnarray}
\hbox{mult}(\chi_\zeta(\la)\chi_\eta(\mu))=
[Z_a]_{\zeta,\eta}[Z_b]_{\zeta,\eta}
\end{eqnarray}
This is a generalization of the classification result
\cite[Theorem 4.16]{BEK2}. We display in \fig{MaMb} the number of
irreducible sectors; each vertex is labelled with an algebra $N$,
$M_a$ or $M_b$ and a line between two algebras represents the
number of irreducible sectors.

Hence if $M_a=M_b$, then $\#\MXM=\Tr(ZZ^t)$ and moreover we can
see in particular that $\MXM$ is commutative if and only if
$Z_{\la,\mu}\leq 1$ for all $\la,\mu$ (so recovering \cite[Theorem
6.8]{BEK1}). Also $b_{\tau,\la}^\pm\leq 1$ for all $\tau, \la$ if
and only if $\MXMpm$ is commutative (as part of
\cite[\eerf{4.11}]{BEK2}). Furthermore, if $M_b=N$ we get
$\Tr(Z_{N\subset M})$ irreducible $M$-$N$ sectors by decomposing
$[\iota \la]$, for $\la$ running in the braided system $\NXN$,
into irreducible sectors and we recover the nimrep constructed in
\cite[Page 768]{BEK2} or \cite[Page 10]{BE6}:
$$G_{\la;a}^b:=\lan \a_\la^\pm a, b\ran=\lan a\la ,b\ran,\quad a,b\in\MXN.$$
As an application, these matrices $G_\la$ are used together with
the fact that $Z:=Z_{N\subset M}$ commutes with $S$ to get the
following curious identity \cite[Proposition 3.3]{E1}:
\begin{equation}
\bigoplus_{a\in\MXN}[a\bar{a}]=\bigoplus_{\la,\mu\in\NXN}Z_{\la,\mu}
[\la\bar{\mu}].
\end{equation}

The matrices $G_\la$ do not depend on the $\pm$ induction choice.
We can consider $G_\la$ as the adjacency matrix of the fusion
graph of $[\a_\la^\pm]$ on the $\MXN$ sectors via left
multiplication. Then the set of matrices $\{G_\la\}$ yield a
nimrep of the underlying modular data \cite{BEK2} whose exponents
match with the exponents of the modular invariant
$Z_{N\subset M}$.
\begin{definition}
A modular invariant $Z$ is  said to be sufferable if it arises
from a braided inclusion $N\subset M$ through the process of
$\a$-induction $[Z]_{\la,\mu}=\lan\a_\la^+,\a_\mu^-\ran$. It is
called insufferable otherwise.
\end{definition}
Therefore every sufferable modular invariant has a compatible
nimrep.
We are also interested in tackling the following well posed problem.\\
{\bf Problem.} For a given modular data, classify their
sufferable modular invariants.
We
assemble the 3 categories (from the ones that encode a very deep
and rich structure to poorer ones) of modular invariants that are
interesting to classify for a fixed RCFT or modular data:
sufferable modular invariants, nimble but insufferable modular
invariants and nimless modular invariants.

\section{Q-systems for inclusions}
\label{frobsub}

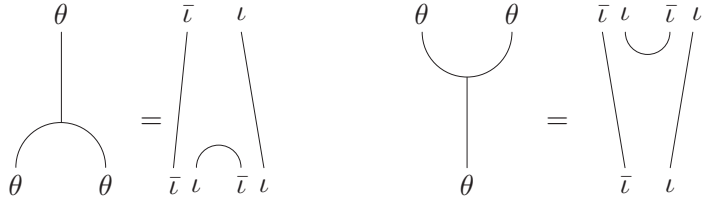
\begin{figure}[htb]
\begin{center}
\unitlength 0.6mm
\begin{picture}(100,55)

\path(-20,20)(-20,40) \put(-20,10){\arc{20}{3.142}{0}}
\put(-20,44){\makebox(0,0){$\theta$}}
\put(-30,6){\makebox(0,0){$\theta$}}
\put(-10,6){\makebox(0,0){$\theta$}}
\put(0,20){\makebox(0,0){$=$}}

\path(5,10)(8,40) \path(25,10)(20,40)
\put(15,10){\arc{10}{3.142}{0}}
\put(5,6){\makebox(0,0){$\bar{\iota}$}}
\put(8,44){\makebox(0,0){$\bar{\iota}$}}
\put(10,6){\makebox(0,0){$\iota$}}
\put(20,6){\makebox(0,0){$\bar{\iota}$}}
\put(25,6){\makebox(0,0){$\iota$}}
\put(20,44){\makebox(0,0){$\iota$}}

\put(70,40){\arc{20}{0}{3.142}} \path(70,30)(70,10)
\put(90,20){\makebox(0,0){$=$}} \put(110,40){\arc{10}{0}{3.142}}
\path(100,40)(105,10) \path(115,10)(120,40)

\put(70,6){\makebox(0,0){$\theta$}}
\put(80,44){\makebox(0,0){$\theta$}}
\put(60,44){\makebox(0,0){$\theta$}}

\put(105,6){\makebox(0,0){$\bar{\iota}$}}
\put(115,6){\makebox(0,0){$\iota$}}
\put(100,44){\makebox(0,0){$\bar{\iota}$}}
\put(115,44){\makebox(0,0){$\bar{\iota}$}}
\put(105,44){\makebox(0,0){$\iota$}}
\put(121,44){\makebox(0,0){$\iota$}}

\end{picture}
\end{center}
\caption{Diagrammatic representation of $\sqrt{d(\iota)}w$ and
$\sqrt{d(\iota)}w^\ast$, respectively} \label{w}
\end{figure}

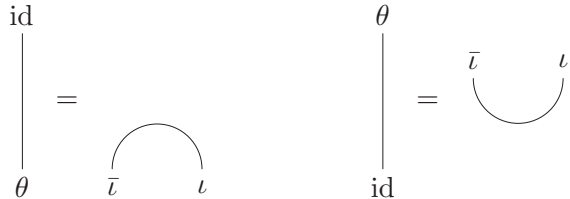
\begin{figure}[htb]
\begin{center}
\unitlength 0.6mm
\begin{picture}(100,55)

\path(-20,30)(-20,0) \put(-20,34){\makebox(0,0){$\hbox{id}$}}
\put(-20,-4){\makebox(0,0){$\theta$}}
\put(-10,15){\makebox(0,0){$=$}} \put(10,0){\arc{20}{3.142}{0}}
\put(0,-4){\makebox(0,0){$\bar{\iota}$}}
\put(20,-4){\makebox(0,0){$\iota$}}

\path(60,30)(60,0) \put(60,34){\makebox(0,0){$\theta$}}
\put(60,-4){\makebox(0,0){$\hbox{id}$}}
\put(70,15){\makebox(0,0){$=$}} \put(90,20){\arc{20}{0}{3.142}}
\put(80,24){\makebox(0,0){$\bar{\iota}$}}
\put(100,24){\makebox(0,0){$\iota$}}

\end{picture}
\end{center}
\caption{Diagrammatic representation of $\sqrt{d(\iota)}v$ and
$\sqrt{d(\iota)}v^\ast$, respectively} \label{v}
\end{figure}

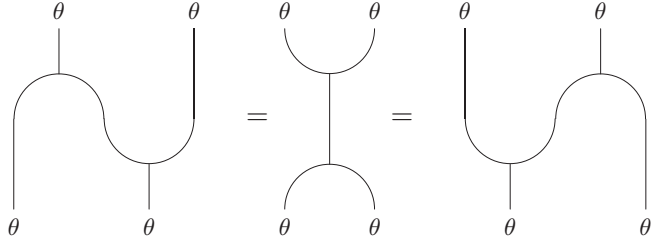
\begin{figure}[htb]
\begin{center}
\unitlength 0.6mm
\begin{picture}(170,40)
\thinlines \path(15,20)(15,0) \put(25,20){\arc{20}{3.142}{0}}
\put(45,20){\arc{20}{0}{3.142}}

\path(25,30)(25,40) \path(45,10)(45,0)

\put(55,40){\line(0,-1){20}} \put(69,20){\makebox(0,0){$=$}}

\put(85,0){\arc{20}{3.142}{0}} \put(85,40){\arc{20}{0}{3.142}}
\path(85,30)(85,10)

\put(101,20){\makebox(0,0){$=$}} \put(115,40){\line(0,-1){20}}
\put(125,20){\arc{20}{0}{3.142}} \put(145,20){\arc{20}{3.142}{0}}
\path(155,20)(155,0) \path(125,10)(125,0) \path(145,30)(145,40)

\put(15,-4){\makebox(0,0){{\footnotesize $\theta$}}}
\put(45,-4){\makebox(0,0){{\footnotesize $\theta$}}}
\put(75,-4){\makebox(0,0){{\footnotesize $\theta$}}}
\put(95,-4){\makebox(0,0){{\footnotesize $\theta$}}}
\put(125,-4){\makebox(0,0){{\footnotesize $\theta$}}}
\put(155,-4){\makebox(0,0){{\footnotesize $\theta$}}}

\put(25,44){\makebox(0,0){{\footnotesize $\theta$}}}
\put(55,44){\makebox(0,0){{\footnotesize $\theta$}}}
\put(75,44){\makebox(0,0){{\footnotesize$\theta$}}}
\put(95,44){\makebox(0,0){{\footnotesize $\theta$}}}
\put(115,44){\makebox(0,0){{\footnotesize $\theta$}}}
\put(145,44){\makebox(0,0){{\footnotesize $\theta$}}}

\end{picture}
\end{center}
\caption{Q-system relation
$\theta(w^\ast)w=ww^\ast=w^\ast\theta(w)$}
\label{md}
\end{figure}

\begin{figure}[htb]
\begin{center}
\unitlength 0.6mm
\begin{picture}(80,45)

\put(0,15){\arc{20}{3.142}{0}}
\put(-10,10){\arc{10}{3.142}{0}}

\path(0,25)(0,35)
\path(-15,10)(-15,0)
\path(-5,10)(-5,0)
\path(10,15)(10,0)

\put(25,15){\makebox(0,0){{\footnotesize $=$}}}

\put(45,15){\arc{20}{3.142}{0}}
\put(55,10){\arc{10}{3.142}{0}}
\path(45,25)(45,35)
\path(50,10)(50,0)
\path(60,10)(60,0)
\path(35,15)(35,0)

\put(-15,-4){\makebox(0,0){{\footnotesize $\theta$}}}
\put(-5,-4){\makebox(0,0){{\footnotesize $\theta$}}}
\put(10,-4){\makebox(0,0){{\footnotesize $\theta$}}}
\put(35,-4){\makebox(0,0){{\footnotesize $\theta$}}}
\put(50,-4){\makebox(0,0){{\footnotesize $\theta$}}}
\put(60,-4){\makebox(0,0){{\footnotesize $\theta$}}}

\put(0,39){\makebox(0,0){{\footnotesize $\theta$}}}
\put(45,39){\makebox(0,0){{\footnotesize $\theta$}}}

\end{picture}
\end{center}
\caption{Q-system relation $d(\iota)w^2=d(\iota)\theta(w)w$}
\label{qsystem1}
\end{figure}
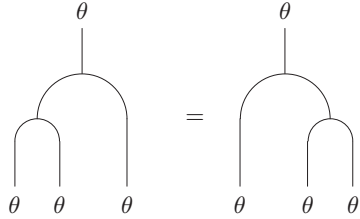

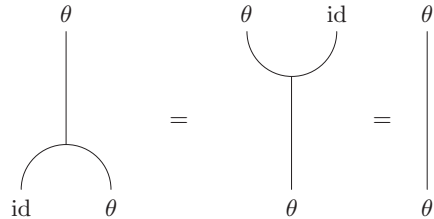
\begin{figure}[htb]
\begin{center}
\unitlength 0.6mm
\begin{picture}(80,47)

\put(0,0){\arc{20}{3.142}{0}}
\path(0,10)(0,35)

\put(25,15){\makebox(0,0){{\footnotesize $=$}}}

\put(50,35){\arc{20}{0}{3.142}}
\path(50,25)(50,0)

\put(-10,-4){\makebox(0,0){{\footnotesize ${\rm id}$}}}
\put(10,-4){\makebox(0,0){{\footnotesize $\theta$}}}
\put(50,-4){\makebox(0,0){{\footnotesize $\theta$}}}

\put(0,39){\makebox(0,0){{\footnotesize $\theta$}}}
\put(40,39){\makebox(0,0){{\footnotesize $\theta$}}}
\put(60,39){\makebox(0,0){{\footnotesize ${\rm id}$}}}

\put(70,15){\makebox(0,0){{\footnotesize $=$}}}

\path(80,0)(80,35)
\put(80,-4){\makebox(0,0){{\footnotesize $\theta$}}}
\put(80,39){\makebox(0,0){{\footnotesize $\theta$}}}

\end{picture}
\end{center}
\caption{Q-system relation $d(\iota)v^\ast w=d(\iota)w^\ast
\theta(v)={\rm id}$} \label{qsystem2}
\end{figure}

Let $N$ be a factor and  $\iota:
N\to M$ be an inclusion in a von Neumann algebra $M$,
with $\bar{\iota}: M\to N$ a conjugate
endomorphism of $\iota$. Then since $ \bar{\iota}\iota$,
$\iota\bar{\iota}$ both contain the identity $\hbox{id}_M$,
$\hbox{id}_N$ respectively, there are intertwining isometries,
$v$ and $w_1$, in Hom$(\hbox{id}_N, \bar{\iota}\iota)$,
Hom$(\hbox{id}_M, \iota\bar{\iota})$ respectively
\cite{longo2,iz2}. Then $w=\bar{\iota}(w_1)$ is an isometry in
Hom$(\theta,\theta^2)$ where $\theta=\bar{\iota}\iota$ is the
dual canonical endomorphism which satisfy \cite{longo}
\begin{eqnarray}
w^\ast\theta(w)=ww^\ast, \quad w^2=\theta(w)w,\quad v^\ast
w=w^\ast\theta(v)=1/d \label{qsystem}
\end{eqnarray}
with $d=d(\iota)$. Hence $d^2=[M:N]$.
\begin{figure}[htb]
\begin{center}
\unitlength 0.6mm
\begin{picture}(30,35)
\thinlines \put(-20,10){\dashbox{2}(20,10){$t$}}
\put(-10,30){\vector(0,-1){10}} \put(-10,10){\vector(0,-1){10}}
\put(-5,25){\makebox(0,0){$\rho$}}
\put(-5,5){\makebox(0,0){$\sigma$}}

\put(10,15){\makebox(0,0){$=$}}

\put(20,10){\dashbox{2}(20,10){$t$}} \path(30,30)(30,20)
\path(30,10)(30,0) \put(35,25){\makebox(0,0){$\rho$}}
\put(35,5){\makebox(0,0){$\sigma$}}

\put(70,10){\dashbox{2}(20,10){$t^\ast$}} \path(80,30)(80,20)
\path(80,10)(80,0) \put(85,25){\makebox(0,0){$\sigma$}}
\put(85,5){\makebox(0,0){$\rho$}}

\end{picture}
\end{center}
\caption{Intertwiners $t\in\hbox{Hom}(\rho,\sigma)$ and
$t^\ast\in\hbox{Hom}(\sigma,\rho)$} \label{trs}
\end{figure}
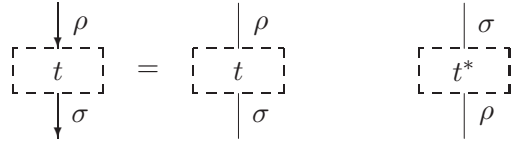
It is convenient to represent intertwiners graphically
\cite{BEK1}. Our convention is that an element $t$ in the
intertwiner space Hom$(\rho,\sigma)$ is written as in \fig{trs}.
With this convention, we write the isometries $w$ and $v$ and as
in \fig{w} and \fig{v}, respectively. Then the
relations of \erf{qsystem} can be displayed graphically as in
Figs.\ \ref{md}, \ref{qsystem1} and \ref{qsystem2}.

\begin{figure}[htb]
\begin{center}
\unitlength 0.6mm
\begin{picture}(100,45)

\put(0,20){\arc{8}{0}{3.143}}
\put(0,15){\arc{10}{0}{3.143}}

\put(9,15){\arc{8}{3.143}{0}}
\put(9,20){\arc{10}{3.143}{0}}

\path(-4,20)(-8,35)\path(-5,15)(-11,35)
\path(13,15)(15,0)\path(14,20)(18,0)

\put(15,-4){\makebox(0,0){{\footnotesize $\bar{\iota}$}}}
\put(19,-4){\makebox(0,0){{\footnotesize $\iota$}}}
\put(-11,39){\makebox(0,0){{\footnotesize $\bar{\iota}$}}}
\put(-7,39){\makebox(0,0){{\footnotesize $\iota$}}}

\put(35,15){\makebox(0,0){{\footnotesize $=$}}}

\put(55,15){\arc{8}{3.143}{0}}
\put(55,20){\arc{10}{3.143}{0}}

\put(64,20){\arc{8}{0}{3.143}}
\put(64,15){\arc{10}{0}{3.143}}

\path(51,15)(50,0)\path(50,20)(47,0)
\path(69,15)(72,35)\path(68,20)(69,35)

\put(47,-4){\makebox(0,0){{\footnotesize $\bar{\iota}$}}}
\put(51,-4){\makebox(0,0){{\footnotesize $\iota$}}}
\put(69,39){\makebox(0,0){{\footnotesize $\bar{\iota}$}}}
\put(73,39){\makebox(0,0){{\footnotesize $\iota$}}}

\end{picture}
\end{center}
\caption{Symmetric Q-system}
\label{symmetric}
\end{figure}
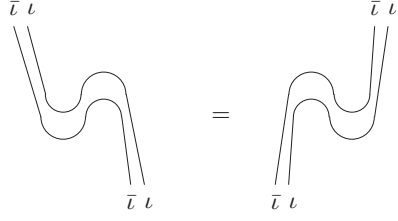

The system $\Theta=(\theta,v,w)$ is called a Q-system by Longo
\cite{longo} which characterises precisely which endomorphisms
can arise as dual canonical endomorphisms for $N\subset M$.
 More
precisely, Longo worked with equivalent notion of characterising
canonical endomorphisms. It is well known that the dual sector
$[\theta]$ does not determine $M$ uniquely up to inner conjugacy.
This is an ${\rm H}^2$ cohomological obstruction that has been
studied in \cite{iz1,KI} where the following definition is
proposed. Two Q-systems $(\theta_1,v_1,w_1)$,
$(\theta_2,v_2,w_2)$ are equivalent if there exists a unitary
$u\in\hbox{Hom}(\theta_1,\theta_2)$ such that $v_2=uv_1$,
$w_2=u\theta_1(u)w_1u^\ast$. In particular we obtain
$\theta_2(\cdot)=u\theta_1(\cdot)u^\ast$ thus
$[\theta_1]=[\theta_2]$. The above relations \erf{qsystem} mean
that a Q-system is a Frobenius algebra
$A=(\theta,m,e,\Delta,\epsilon)$ where $e\in \hbox{ Hom}({\bf
1},\theta)$, $m\in\hbox{ Hom}(\theta^2, \theta)$,
$\epsilon\in\hbox{ Hom}(\theta, {\bf 1})$, $\Delta\in \hbox{
Hom}(\theta, \theta^2)$ such that $(\theta,m,e)$ is an algebra,
$(\theta,\Delta,\epsilon)$ is a co-algebra with the algebraic and
co-algebraic structure related by
$$(\hbox{id}_\Theta\otimes m)\circ(\Delta\otimes \hbox{id}_\Theta)=
\Delta\circ m=(m\otimes
\hbox{id}_\Theta)\circ(\hbox{id}_\Theta\otimes\Delta).$$ This is
in as \fig{md} where $m$ and $\Delta$ are identified with $w^\ast$
and $w$ respectively. A Frobenius algebra $A$ is said to be
special \cite[Definition 3.4]{frs} if there are non-zero constants
$\beta_{\bf 1}$ and $\beta_A$ such that $\varepsilon\circ
e=\beta_1\hbox{id}_{\bf 1}$, $m\circ\Delta=\beta_A\hbox{id}_A$.
For a Q-system $\Theta=(\theta,v,w)$, $v^\ast v={\bf 1}$ and
$w^\ast w={\bf 1}$, implying that $\Theta$ is a special
$\ast$-Frobenius algebra (with $\beta_{\bf
1}=\beta_\Theta=d(\iota)$, see Figs.\ \ref{w} and \ref{v}). A
Q-system $\Theta$ is automatically a {\it symmetric} Frobenius
algebra \cite[Definition 3.4]{frs}, as in \fig{symmetric} since
${\rm id} = d^2 v^\ast
w^\ast\theta(w)\theta(v)=d^2\theta(v^\ast)\theta(w^\ast)wv.$

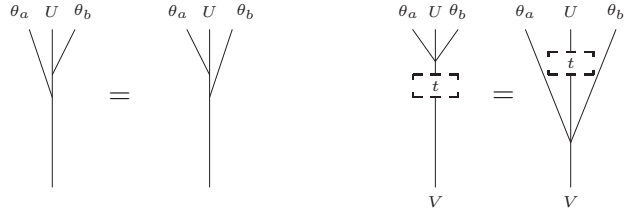
\begin{figure}[htb]
\begin{center}
\unitlength 0.6mm
\begin{picture}(80,45)


\path(-20,-5)(-20,30) \path(-20,15)(-25,30) \path(-20,20)(-15,30)

\put(-20,34){\makebox(0,0){{\tiny $U$}}}
\put(-27,34){\makebox(0,0){{\tiny $\theta_a$}}}
\put(-13,34){\makebox(0,0){{\tiny $\theta_b$}}}

\put(-5,15){\makebox(0,0){$=$}}

\path(15,-5)(15,30)
\path(15,20)(10,30)
\path(15,15)(20,30)

\put(15,34){\makebox(0,0){{\tiny $U$}}}
\put(8,34){\makebox(0,0){{\tiny $\theta_a$}}}
\put(23,34){\makebox(0,0){{\tiny $\theta_b$}}}

\put(60,15){\dashbox{2}(10,5){{{\tiny $t$}}}} \path(65,15)(65,-5)
\path(65,20)(65,30) \path(65,23)(60,30) \path(65,23)(70,30)
\put(80,15){\makebox(0,0){$=$}}

\put(90,20){\dashbox{2}(10,5){{{\tiny $t$}}}} \path(95,-5)(95,20)
\path(95,25)(95,30) \path(95,5)(85,30) \path(95,5)(105,30)

\put(65,34){\makebox(0,0){{\tiny $U$}}}
\put(59,34){\makebox(0,0){{\tiny $\theta_a$}}}
\put(70,34){\makebox(0,0){{\tiny $\theta_b$}}}
\put(95,34){\makebox(0,0){{\tiny $U$}}}
\put(85,34){\makebox(0,0){{\tiny $\theta_a$}}}
\put(105,34){\makebox(0,0){{\tiny $\theta_b$}}}

\put(65,-8){\makebox(0,0){{\tiny $V$}}}
\put(95,-8){\makebox(0,0){{\tiny $V$}}}

\end{picture}
\end{center}
\caption{$\Theta_a$-$\Theta_b$ bimodule}
\label{AB}
\end{figure}

\begin{figure}[htb]
\begin{center}

\unitlength 0.6mm
\begin{picture}(75,60)

\thicklines \path(-20,25)(-30,50)

\thinlines \path(-20,50)(-20,5)

\put(-20,25){\makebox(0,0){{\tiny $\bullet$}}}
\put(-30,54){\makebox(0,0){{\tiny $\theta_a$}}}
\put(-20,54){\makebox(0,0){{\tiny $U$}}}
\put(-25,-5){\makebox(0,0){{\tiny (left) $\Theta_a$-module}}}

\path(40,50)(40,5) \thicklines \path(40,25)(30,50)
\path(40,20)(47,50)

\thinlines \put(40,-5){\makebox(0,0){{\tiny $\Theta_a$-$\Theta_b$
bimodule}}} \put(30,54){\makebox(0,0){{\tiny $\theta_a$}}}
\put(47,54){\makebox(0,0){{\tiny $\theta_b$}}}
\put(40,54){\makebox(0,0){{\tiny $U$}}}
\put(40,25){\makebox(0,0){{\tiny $\bullet$}}}
\put(40,20){\makebox(0,0){{\tiny $\bullet$}}}

\put(50,25){\vector(1,0){8}}

\path(75,50)(75,5)

\thicklines \path(75,25)(65,50)
\path(75,20)(78,23)(78,25)(77,26.5)(76.5,28)

 \path(70,50)(74,34)



\thinlines

\put(75,25){\makebox(0,0){{\tiny $\bullet$}}}
\put(75,20){\makebox(0,0){{\tiny $\bullet$}}}

\put(85,25){\makebox(0,0){{\tiny $\equiv$}}}

\thicklines \path(105,25)(95,50)

\thinlines \path(105,50)(105,5) \put(105,25){\makebox(0,0)
{{\tiny $\bullet$}}}

\put(93,54){\makebox(0,0){{\tiny $\theta_a\theta_b$}}}

\put(108,54){\makebox(0,0){{\tiny $U$}}}

\put(105,-5){\makebox(0,0){{\tiny
$\Theta_a$$\otimes$${\Theta_b}^{\rm{opp}}$-module}}}

\end{picture}
\end{center}
\caption{Notation for $\Theta_a$-modules and
$\Theta_a$-$\Theta_b$ bimodules} \label{Amod}
\end{figure}
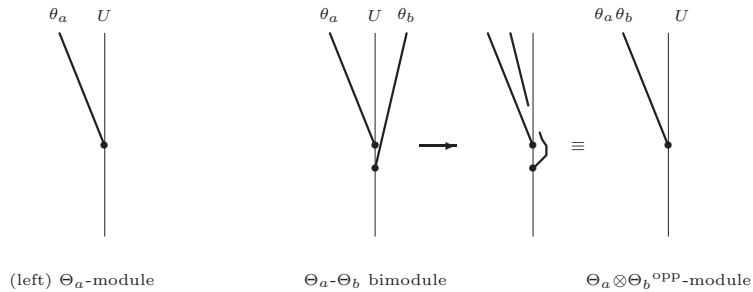

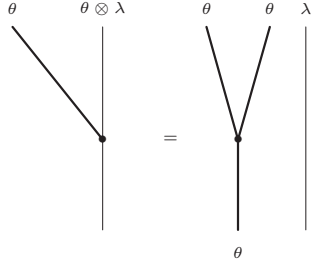
\begin{figure}[htb]
\begin{center}
\unitlength 0.6mm
\begin{picture}(75,60)

\thicklines \path(20,25)(0,50) \thinlines \path(20,50)(20,5)
\put(20,25){\makebox(0,0){{\tiny $\bullet$}}}
\put(0,54){\makebox(0,0){{\tiny $\theta$}}}
\put(20,54){\makebox(0,0){{\tiny $\theta\otimes\la$}}}

\thinlines \put(35,25){\makebox(0,0){{\tiny $=$}}}

\thicklines

\path(50,5)(50,25) \path(43,50)(50,25)(57,50)

 \thinlines \path(65,50)(65,5)
\put(50,25){\makebox(0,0){{\tiny $\bullet$}}}

\put(57,54){\makebox(0,0){{\tiny $\theta$}}}
\put(43,54){\makebox(0,0){{\tiny $\theta$}}}
\put(50,0){\makebox(0,0){{\tiny $\theta$}}}
\put(65,54){\makebox(0,0){{\tiny $\la$}}}

\end{picture}
\end{center}
\caption{Induced (left) $\Theta$-module module Ind$_\Theta(\la)$
for $\la\in\NXN$} \label{induced}
\end{figure}

For a Q-system $\Theta=(\theta, v,w)$, a (left) $\Theta$-module
\cite{ost1} (see \fig{Amod}) is a pair $(U,\rho)$, where $U\in\Sigma(\NXN)$
and $\rho\in\hbox{Hom}(\theta\otimes U,U)$ such that
\begin{eqnarray}
\rho\circ(m\otimes\hbox{id})=\rho\circ(\hbox{id}\otimes\rho),\quad
\rho\circ(e\otimes\hbox{id})= \hbox{id}.
\end{eqnarray}
In particular, we can define the induced $\Theta$-modules
Ind$_\Theta(U)=\theta\otimes U$ for any $U\in\Sigma(\NXN)$, and it
is the case that any simple $\Theta$-module is a submodule of
Ind$_\Theta(\la)$, for some $\la\in\NXN$ (see \fig{induced}). A
$\Theta_a$-$\Theta_b$ bimodule is a triple
$(U,\rho_{a},\rho_{b})$ with $U\in\Sigma(\NXN)$,
$\rho_{a}\in\hbox{Hom}(\theta_a\otimes\nobreak U,U)$ and
$\rho_{b}\in\hbox{Hom}(U\otimes \theta_b,U)$ such that
$(U,\rho_{a})$ is a left $\Theta_a$-module, $(U,\rho_{b})$ a
right $\Theta_b$-module so that
$$\rho_{a}\circ(\hbox{id}_{\Theta_a}\otimes\rho_{b})=
\rho_{b}\circ(\rho_{a} \otimes\hbox{id}_{\Theta_b})$$ as LHS of
\fig{AB}. A map $t:U\to V$ is an $\Theta_a$-$\Theta_b$
intertwiner if the condition on RHS of \fig{AB} holds.

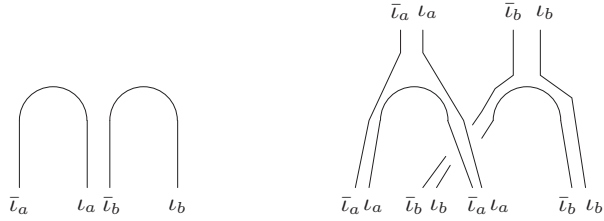
\begin{figure}[htb]
\begin{center}
\unitlength 0.6mm
\begin{picture}(90,45)

\put(-10,10){\arc{15}{3.143}{0}}
\put(-30,10){\arc{15}{3.143}{0}}
\path(-17.5,10)(-17.5,-5)\path(-2.5,10)(-2.5,-5)
\path(-22.5,10)(-22.5,-5)\path(-37.5,10)(-37.5,-5)

\put(-37.5,-9){\makebox(0,0){{\footnotesize $\bar{\iota}_a$}}}
\put(-22.5,-9){\makebox(0,0){{\footnotesize $\iota_a$}}}
\put(-17,-9){\makebox(0,0){{\footnotesize $\bar{\iota}_b$}}}
\put(-2.5,-9){\makebox(0,0){{\footnotesize $\iota_b$}}}

\put(50,10){\arc{15}{3.143}{0}}
\put(75,10){\arc{15}{3.143}{0}}
\path(57.5,10)(63,-5)\path(65,-5)(61,10)(52,25)(52,30)
\path(42.5,10)(40,-5)\path(37,-5)(40,10)(47,25)(47,30)

\path(82.5,10)(85,-5)
\path(67.5,10)(65,6)\path(58,0)(55,-5)
\path(52,-5)(57,2)
\path(63,9)(65,12)(68,17)(72,20)(72,30)
\path(88,-5)(85,15)(78,20)(78,30)

\put(36,-9){\makebox(0,0){{\footnotesize $\bar{\iota}_a$}}}
\put(41,-9){\makebox(0,0){{\footnotesize $\iota_a$}}}
\put(50,-9){\makebox(0,0){{\footnotesize $\bar{\iota}_b$}}}
\put(55.5,-9){\makebox(0,0){{\footnotesize $\iota_b$}}}

\put(64,-9){\makebox(0,0){{\footnotesize $\bar{\iota}_a$}}}
\put(69,-9){\makebox(0,0){{\footnotesize $\iota_a$}}}
\put(84,-9){\makebox(0,0){{\footnotesize $\bar{\iota}_b$}}}
\put(90,-9){\makebox(0,0){{\footnotesize $\iota_b$}}}

\put(47,34){\makebox(0,0){{\footnotesize $\bar{\iota}_a$}}}
\put(53,34){\makebox(0,0){{\footnotesize $\iota_a$}}}
\put(72,34){\makebox(0,0){{\footnotesize $\bar{\iota}_b$}}}
\put(79,34){\makebox(0,0){{\footnotesize $\iota_b$}}}

\end{picture}
\end{center}
\caption{Braided product of Q-systems $\Theta_{ab}=\Big(\theta_a\theta_b,
\theta_a(v_b)v_a,\theta_a(\epsilon{(\theta_a,\theta_b)})
\theta_a^2(w_b)w_a\Big)$}
\label{thetaproduct}
\end{figure}

Given Q-systems $\Theta_a$ and $\Theta_b$, the simple
$\Theta_a$-$\Theta_b$ bimodules can be identified with the simple
$\Theta_a\otimes {\Theta_b}^{\rm{opp}}$-modules as in \fig{Amod}
where ${\Theta_b}^{\rm{opp}}$ is the opposite algebra \cite[Remark
12]{ost1}. Recall
that if $\Theta=(\theta,v,w)$ is a Q-system on a braided factor
$N$, then we can define an associated opposite Q-system
$(\theta,v,\varepsilon_{(\theta,\theta)} w)$ denoted by
$\Theta^{\rm{opp}}$ as in \fig{opposite}. The graphical
representation of
$\epsilon_{(\la,\mu)}\equiv\epsilon_{(\la,\mu)}^+$ is in
\fig{epsilonpm}. We remark that if $\Theta$ and
${\Theta}^{\rm{opp}}$ are equivalent then the modular invariant
$Z$ is symmetric $Z$=$Z^t$, but by \cite{BE4}, there are
Q-systems producing non-symmetric modular invariants, thus we may
have $\Theta\not\simeq{\Theta}^{\rm{opp}}$.
The $\Phi$ product $\Theta_{a} \otimes {\Theta_{b}}= 
\Theta_{{ab}}$ of two Q-systems
$\Theta_a=(\theta_a, v_a, w_a)$  and $\Theta_b=(\theta_b,v_b,w_b)$ is
 $(\theta_a\theta_b,  \theta_a(v_b)v_a,
\theta_a(\epsilon(\theta_a,\theta_b))\theta_a^2(w_b)w_a)$
as in \fig{thetaproduct}. For completion, we define here the
direct sum $\Theta_{a}\oplus\Theta_{b}$ of Q-systems,
whose associated braided inclusion is $N\subset
M_{a}\oplus M_b$, if the inclusions $N\subset
M_a, N\subset M_b$ represent $\Theta_{a}$, and $\Theta_{b}$
respectively.
Let
$s_a,s_b\in N$ be Cuntz generators, i.e.\ $s_i^\ast s_j=
\delta_{i,j}{\bf 1}$ and
$s_as_a^\ast+s_bs_b^\ast={\bf 1}$. Define now
\begin{eqnarray*}
\theta(n)&=&s_a\theta_a(n)s_a^\ast+s_b\theta_b(n)s_b^\ast,\quad n\in N\\
v&=& (\sqrt{d_a}s_av_a+\sqrt{d_b}s_bv_b)/\sqrt{d(\theta)},\\
w&=&\theta(s_a)s_aw_as_a^\ast+\theta(s_b)s_bw_bs_b^\ast.
\end{eqnarray*}
where $d(\theta)=d(\theta_a)+d(\theta_b)$,
and $d_a=d(\theta_a), d_b=d(\theta_b)$.

\begin{figure}[htb]
\begin{center}
\unitlength 0.6mm
\begin{picture}(100,40)

\path(10,30)(30,0) \path(30,30)(21,16) \path(17,10)(10,0)

\put(10,34){\makebox(0,0){{\footnotesize $\la$}}}
\put(30,34){\makebox(0,0){{\footnotesize $\mu$}}}
\put(30,-4){\makebox(0,0){{\footnotesize $\la$}}}
\put(10,-4){\makebox(0,0){{\footnotesize $\mu$}}}


\path(70,30)(78,18) \path(83,11)(90,0) \path(90,30)(70,0)

\put(70,34){\makebox(0,0){{\footnotesize $\mu$}}}
\put(90,34){\makebox(0,0){{\footnotesize $\la$}}}
\put(90,-4){\makebox(0,0){{\footnotesize $\mu$}}}
\put(70,-4){\makebox(0,0){{\footnotesize $\la$}}}

\end{picture}
\end{center}
\caption{The braiding operators $\epsilon_{(\la,\mu)}^+$ and
$\epsilon_{(\la,\mu)}^-:=\epsilon_{(\mu,\la)}^{+\quad \ast}$ as
over- and undercrossings} \label{epsilonpm}
\end{figure}
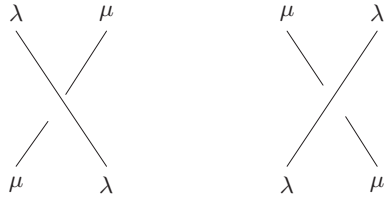

\begin{figure}[htb]
\begin{center}
\unitlength 0.6mm
\begin{picture}(60,60)

\path(20,30)(20,50) \put(20,20){\arc{20}{3.142}{0}}
\put(20,54){\makebox(0,0){${\footnotesize \theta}$}}
\path(10,20)(30,-5) \path(30,20)(23,10) \path(16,3)(10,-5)
\put(30,-10){\makebox(0,0){${\footnotesize \theta}$}}
\put(10,-10){\makebox(0,0){${\footnotesize \theta}$}}

\end{picture}
\end{center}
\caption{The isometry $\epsilon_{(\theta,\theta)}w$}
\label{opposite}
\end{figure}
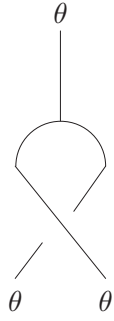

\begin{figure}[htb]
\begin{center}
\unitlength 0.6mm
\begin{picture}(75,60)

\thicklines \path(0,25)(-20,50)

\thinlines

\path(0,50)(0,5) \put(0,25){\makebox(0,0){{\tiny $\bullet$}}}
\put(-20,54){\makebox(0,0){{\tiny $\theta$}}}
\put(0,54){\makebox(0,0){{\tiny $U$}}}

\put(20,30){\makebox(2,3){{$\equiv$}}}

\path(40,25)(20,50) \path(40,25)(25,50) \path(40,50)(40,5)
\put(40,25){\makebox(0,0){{\tiny $\bullet$}}}
\put(20,54){\makebox(0,0){{\tiny $\bar{\iota}$}}}
\put(25,54){\makebox(0,0){{\tiny $\iota$}}}
\put(40,54){\makebox(0,0){{\tiny $U$}}}

\path(80,50)(80,5) \put(85,25){\vector(1,0){8}} \path(95,50)(95,5)
\path(100,50)(100,5) \path(106,50)(106,5)
 \put(80,54){\makebox(0,0){{\tiny $\beta$}}}
\put(90,27){\makebox(0,0){{\tiny $\Phi$}}}
\put(95,54){\makebox(0,0){{\tiny $\bar{\iota}_a$}}}
\put(101,54){\makebox(0,0){{\tiny $\beta$}}}
\put(107,54){\makebox(0,0){{\tiny $\iota_b$}}}

\end{picture}
\end{center}
\caption{$\Theta$ modules and the map $\Phi$}
 \label{corr}
\end{figure}
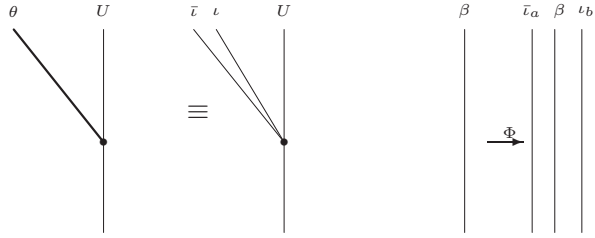

\begin{figure}[htb]
\begin{center}
\unitlength 0.6mm
\begin{picture}(110,60)

\thinlines \path(-25,50)(-25,30) \path(-25,20)(-25,0)
\put(-30,20){\dashbox{5}(10,10){\tiny $t$}}
\put(-25,54){\makebox(0,0){{\tiny $\beta$}}}
\put(-25,-2){\makebox(0,0){{\tiny $\beta^\prime$}}}
\put(-15,25){\vector(1,0){8}} \put(-12,27){\makebox(0,0){{\tiny
$\Phi$}}} \path(0,50)(0,0) \path(10,50)(10,30) \path(10,20)(10,0)
\put(5,20){\dashbox{5}(10,10){\tiny $t$}}
\put(10,54){\makebox(0,0){{\tiny $\beta$}}}
\put(10,-2){\makebox(0,0){{\tiny $\beta^\prime$}}}
\put(0,54){\makebox(0,0){{\tiny $\bar{\iota}_a$}}}

\path(20,0)(20,50) \put(20,54){\makebox(0,0){{\tiny $\iota_b$}}}
\put(45,25){\makebox(0,0){{\tiny $\in{\rm
Hom}(\bar{\iota}_a\beta\iota_b,\bar{\iota}_a\beta^\prime\iota_b$)}}}
\put(75,25){\makebox(0,0){{\tiny $\equiv$}}}

\path(86,50)(86,30) \path(83.5,50)(83.5,30) \path(86,20)(86,0)
\path(83.5,20)(83.5,0) \put(80,20){\dashbox{5}(12,10){{\tiny
$t^\prime$}}} \put(87,54){\makebox(0,0){{\tiny $\beta$}}}
\put(87,-2){\makebox(0,0){{\tiny $\beta^\prime$}}}
\put(82.5,54){\makebox(0,0){{\tiny $\bar{\iota}_a$}}}
\put(82.5,-2){\makebox(0,0){{\tiny $\bar{\iota}_a$}}}
\path(90,0)(90,20) \path(90,30)(90,50)
\put(92,54){\makebox(0,0){{\tiny $\iota_b$}}}
\put(92,-2){\makebox(0,0){{\tiny $\iota_b$}}}

\put(100,25){\makebox(0,0){{\tiny $=$}}}

\path(110,50)(110,30) \path(110,20)(110,0)
\put(105,20){\dashbox{5}(10,10){\tiny $t^\prime$}}
\put(110,54){\makebox(0,0){{\tiny $\bar{\iota}_a\beta\iota_b$}}}
\put(110,-2){\makebox(0,0){{\tiny
$\bar{\iota}_a\beta^\prime\iota_b$}}}

\end{picture}
\end{center}
\caption{From $M_a$-$M_b$ intertwiners to
$\Theta_a$-$\Theta_b$-intertwiners} \label{mapping}
\end{figure}
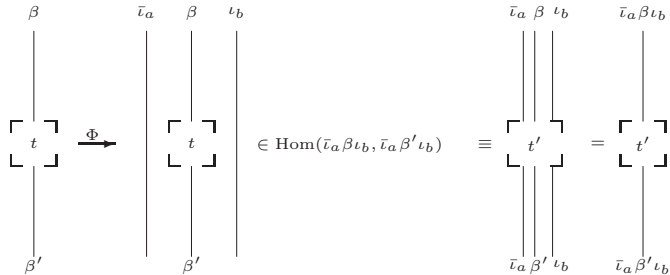

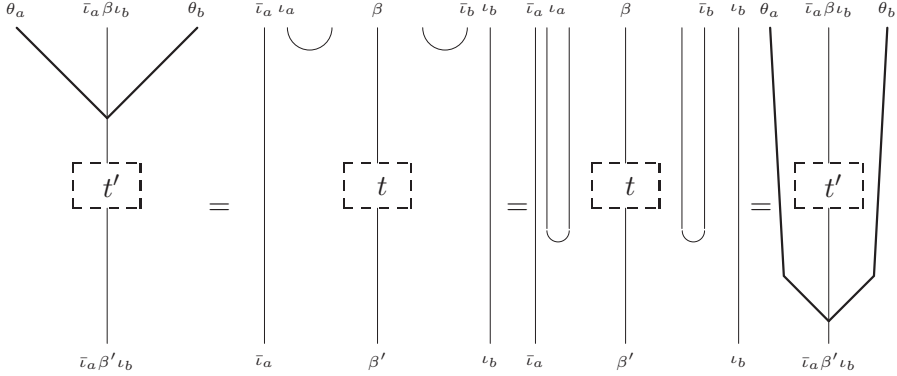
\begin{figure}[htb]
\begin{center}
\unitlength 0.6mm
\begin{picture}(100,90)
\thinlines

\put(-37.5,30){\dashbox{2}(15,10){{ $t^\prime$}}}
\path(-30,0)(-30,30)

\put(-30,-4){\makebox(0,0){{\tiny
$\bar{\iota}_a\beta^\prime\iota_b$}}}

\path(-30,40)(-30,70)

\put(-30,74){\makebox(0,0){{\tiny $\bar{\iota}_a\beta\iota_b$}}}

\thicklines

\path(-30,50)(-50,70) \path(-30,50)(-10,70)
\put(-50,74){\makebox(0,0){{\tiny $\theta_a$}}}
\put(-10,74){\makebox(0,0){{\tiny $\theta_b$}}}

\thinlines

\put(-5,30){\makebox(0,0){{$=$}}}

\put(22.5,30){\dashbox{2}(15,10){{ $t$}}} \path(30,0)(30,30)
\path(30,40)(30,70)

\path(5,0)(5,70) \put(15,70){\arc{10}{0}{3.145}}
\put(45,70){\arc{10}{0}{3.145}}

\path(55,0)(55,70)

\put(5,-4){\makebox(0,0){{\tiny $\bar{\iota}_a$}}}
\put(5,74){\makebox(0,0){{\tiny $\bar{\iota}_a$}}}
\put(55,-4){\makebox(0,0){{\tiny $\iota_b$}}}
\put(55,74){\makebox(0,0){{\tiny $\iota_b$}}}
\put(30,-4){\makebox(0,0){{\tiny $\beta^\prime$}}}
\put(30,74){\makebox(0,0){{\tiny $\beta$}}}

\put(10,74){\makebox(0,0){{\tiny $\iota_a$}}}
\put(50,74){\makebox(0,0){{\tiny $\bar{\iota}_b$}}}

\put(61,30){\makebox(0,0){{$=$}}}

\put(77.5,30){\dashbox{2}(15,10){{ $t$}}} \path(85,0)(85,30)
\path(85,40)(85,70)

\path(65,0)(65,70)

\put(70,25){\arc{5}{0}{3.145}} \put(100,25){\arc{5}{0}{3.145}}
\path(67.5,25.5)(67.5,70) \path(72.5,25.5)(72.5,70)

\path(97.5,25.5)(97.5,70)

\path(102.5,25.5)(102.5,70)

\path(110,0)(110,70)

\put(65,-4){\makebox(0,0){{\tiny $\bar{\iota}_a$}}}
\put(65,74){\makebox(0,0){{\tiny $\bar{\iota}_a$}}}
\put(110,-4){\makebox(0,0){{\tiny $\iota_b$}}}
\put(110,74){\makebox(0,0){{\tiny $\iota_b$}}}
\put(85,-4){\makebox(0,0){{\tiny $\beta^\prime$}}}
\put(85,74){\makebox(0,0){{\tiny $\beta$}}}

\put(70,74){\makebox(0,0){{\tiny $\iota_a$}}}
\put(103,74){\makebox(0,0){{\tiny $\bar{\iota}_b$}}}

\put(115,30){\makebox(0,0){{$=$}}}

\put(122.5,30){\dashbox{2}(15,10){{ $t^\prime$}}}
\path(130,0)(130,30) \path(130,40)(130,70)
\put(130,74){\makebox(0,0){{\tiny $\bar{\iota}_a\beta\iota_b$}}}
\put(130,-4){\makebox(0,0){{\tiny
$\bar{\iota}_a\beta^\prime\iota_b$}}}

\thicklines

\path(130,5)(120,15)(117,70) \path(130,5)(140,15)(143,70)


\put(117,74){\makebox(0,0){{\tiny $\theta_a$}}}
\put(143,74){\makebox(0,0){{\tiny $\theta_b$}}}

\end{picture}
\end{center}
\caption{$\Phi$ is an $\Theta_a$-$\Theta_b$ morphism}

\label{ainter}
\end{figure}

\begin{figure}[htb]
\begin{center}
\unitlength 0.6mm
\begin{picture}(100,90)
\thinlines

\put(-40,30){\dashbox{2}(20,10){{ $t$}}}

\path(-40,0)(-40,30) \path(-30,0)(-30,30) \path(-20,0)(-20,30)
\path(-40,40)(-44,70) \path(-20,40)(-18,70)\path(-30,40)(-30,70)

\put(-37,70){\arc{5}{0}{3.145}} \put(-25,70){\arc{5}{0}{3.145}}
\put(-40,-4){\makebox(0,0){{\tiny $\bar{\iota}_a$}}}
\put(-30,-4){\makebox(0,0){{\tiny $\beta^\prime$}}}
\put(-20,-4){\makebox(0,0){{\tiny $\iota_b$}}}

\put(-44,74){\makebox(0,0){{\tiny $\bar{\iota}_a$}}}
\put(-37,74){\makebox(0,0){{\tiny $\iota_a$}}}
\put(-30,74){\makebox(0,0){{\tiny $\beta$}}}
\put(-15,74){\makebox(0,0){{\tiny $\iota_b$}}}
\put(-22,74){\makebox(0,0){{\tiny $\bar{\iota}_b$}}}

\put(-12,35){\makebox(0,0){{\tiny $=$}}}

\put(5,30){\dashbox{2}(20,10){{ $t$}}}
\path(0,0)(0,70)\path(30,0)(30,70)
\path(15,0)(15,30)\path(15,40)(15,70)

\put(5,30){\arc{5}{0}{3.145}} \put(25,30){\arc{5}{0}{3.145}}
\path(7.5,40)(7.5,70)

\path(21,40)(21,70)

\path(2.5,30)(2.5,70) \path(27.5,30)(27.5,70)

\put(15,74){\makebox(0,0){{\tiny $\beta$}}}
\put(15,-4){\makebox(0,0){{\tiny $\beta^\prime$}}}
\put(0,-4){\makebox(0,0){{\tiny $\bar{\iota}_a$}}}
\put(30,-4){\makebox(0,0){{\tiny $\iota_b$}}}

\put(9,74){\makebox(0,0){{\tiny $\iota_a$}}}
\put(21,74){\makebox(0,0){{\tiny $\bar{\iota}_b$}}}

\thicklines \put(40,35){\makebox(0,0){{\bf {\tiny
$\Longrightarrow$}}}}

\thinlines

\put(55,30){\dashbox{2}(20,10){{ $t$}}} \path(55,0)(55,30)
\path(55,40)(53,70)

\path(65,0)(65,30) \path(65,40)(65,70) \path(75,0)(75,30)
\path(75,40)(78,70)

\put(60,65){\circle{6}} \put(70,65){\circle{6}}

\put(55,-4){\makebox(0,0){{\tiny $\bar{\iota}_a$}}}
\put(65,-4){\makebox(0,0){{\tiny $\beta^\prime$}}}
\put(75,-4){\makebox(0,0){{\tiny $\iota_b$}}}

\put(52,74){\makebox(0,0){{\tiny $\bar{\iota}_a$}}}
\put(58,70){\makebox(0,0){{\tiny $\iota_a$}}}
\put(65,74){\makebox(0,0){{\tiny $\beta$}}}
\put(70,70){\makebox(0,0){{\tiny $\iota_b$}}}
\put(79,74){\makebox(0,0){{\tiny $\iota_b$}}}

\put(85,35){\makebox(0,0){{\tiny $=$}}}
\put(100,30){\dashbox{2}(20,10){{ $t^\prime$}}} \path(90,0)(90,70)
\path(110,0)(110,30) \path(110,40)(110,70) \path(130,0)(130,70)
\put(100,30){\arc{5}{0}{3.145}} \put(120,30){\arc{5}{0}{3.145}}
\put(100,40){\arc{5}{3.145}{0}} \put(120,40){\arc{5}{3.145}{0}}
\path(97.5,30)(97.5,40) \path(122.5,30)(122.5,40)

\put(90,-4){\makebox(0,0){{\tiny $\bar{\iota}_a$}}}
\put(110,-4){\makebox(0,0){{\tiny $\beta^\prime$}}}
\put(110,74){\makebox(0,0){{\tiny $\beta$}}}
\put(130,-4){\makebox(0,0){{\tiny $\iota_b$}}}

\put(100,44){\makebox(0,0){{\tiny $\iota_a$}}}
\put(120,44){\makebox(0,0){{\tiny $\iota_b$}}}

\end{picture}
\end{center}
\caption{Proving that $t={\bf 1}_{\bar{\iota}_a} \otimes
t^{\prime\prime}\otimes {\bf 1}_{\iota_b}$} \label{onto}
\end{figure}

\begin{lemma}
Let $\Theta_a$ and $\Theta_b$ be Q-systems associated with the
braided subfactors $N\subset M_a$ and $N\subset M_b$,
respectively. Then we can identify the category of
$\Theta_a$-$\Theta_b$ bimodules with the category of $M_a$-$M_b$
sectors $\MaXMb$. \lablth{categ}
\end{lemma}
\bproof Denote by $\iota_a$ and $\iota_b$ the inclusion maps of
the subfactors $N\subset M_a$ and $N\subset M_b$ respectively.
Every irreducible $\beta$ in $\Sigma(\MaXMb)$ arises from the
decomposition $\iota_a\la\bar{\iota}_b$ with $\la$ in $\NXN$.
Define now $\Phi:\Sigma(\MaXMb)\to \Theta_a$-$\Theta_b$-bimodules
by $\Phi(\beta)=\bar{\iota}_a\beta\iota_b$ for $\beta\in\MaXMb$
(see RHS of \fig{corr}). In particular,
$\Phi(\iota_a\la\bar{\iota}_b)=\theta_a\la\theta_b$. If
$\beta,\beta^\prime\in\Sigma(\MaXMb)$, then we map an intertwiner
$t\in\hbox{Hom}(\beta,\beta^\prime)$ to $\Phi(t)=\bar{\iota}_a t
\iota_b\in\hbox{Hom}(\bar{\iota}_a\beta\iota_b,
\bar{\iota}_a\beta^\prime\iota_b)=
\hbox{Hom}(\Phi(\beta),\Phi(\beta^\prime)$, see RHS of
\fig{mapping}. Then $\Phi$ is an $\Theta_a$-$\Theta_b$ morphism
as in \fig{ainter}. Let us finally prove that $\Phi$ is injective.
Suppose that
$\bar{\iota}_a\beta\iota_b\simeq\bar{\iota}_a\beta^\prime\iota_b$
as $\Theta_a$-$\Theta_b$ bimodules. Then we prove in \fig{onto}
that $t={\bf 1}_{\bar{\iota}_a}\otimes t^{\prime\prime}
\otimes{\bf 1}_{\iota_b}$ with
$t^{\prime\prime}\in\hbox{Hom}(\beta,\beta^\prime)$. If $t$ is an
isomorphism so is $t^{\prime\prime}$, therefore $\beta\simeq\beta^\prime$.
We have $\#\MaXMb=\Tr(Z_aZ_b^t)$ by \cite{E1}\footnote{This was shown
as in  \sect{mixing}
for subfactors and normalised invariants, but the results extend
to inclusions using central decompositions as in \sect{frsuff},
particularly as in Lemma \ref{sumsmod}. Similarly for the generating
property for $M$-$M$ sectors by $\a$-induced ones.}, so since $\Phi$ is
injective $\#\MaXMb\leq\#$(irreducible $\Theta_a$-$\Theta_b$
bimodules). On the other hand, by \cite[Proposition 5.16]{frs}\footnote{
We do not need the full strength of \cite[Theorem 5.18]{frs}},
$\Tr(Z_aZ_b^t)\geq \#(\hbox{irreducible}\ \Theta_a\otimes
{\Theta_b}^{{\rm opp}}- \hbox{modules})=\#$(irreducible
$\Theta_a$-$\Theta_b$ bimodules). So $\#\MaXMb$=\#(irreducible
$\Theta_a$-$\Theta_b$ bimodules). Therefore $\Phi$ is surjective.
\eproof

There are  particular cases from Lemma \ref{categ}: namely,
when $N=M_a$ and $M_a=M_b$. In the first case the lemma implies
that for the Q-system $\Theta_b=:\Theta$ with $N\subset M$ the
associated braided subfactor, we can identify the category of
$\Theta$-modules with the category $\Sigma(\MXN)$ compatible with
the left $\NXN$ and right $\MXM$ actions. The irreducible $N$-$M$
sectors $\beta$ in $\MXN$ arise from the decompositon of
$\iota\la$, for $\la\in\NXN$. (The map $\Phi:\Sigma(\MXN)\to$
$\Theta$-modules is given by $\Phi(\beta)=\bar{\iota}\beta$ for
$\beta\in\MXN$. In particular
$\iota\la\mapsto\bar{\iota}\iota\la=\hbox{Ind}_\Theta(\la)$.)
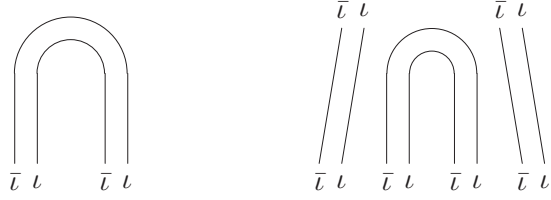
\begin{figure}[htb]
\begin{center}
\unitlength 0.6mm
\begin{picture}(100,45)

\put(-20,20){\arc{15}{3.143}{0}} \put(-20,20){\arc{25}{3.143}{0}}
\path(-27.5,20)(-27.5,0)\path(-12.5,20)(-12.5,0)
\path(-32.5,20)(-32.5,0)
\path(-7.5,20)(-7.5,0)

\put(-32.5,-4){\makebox(0,0){$\bar{\iota}$}}
\put(-27.5,-4){\makebox(0,0){$\iota$}}
\put(-12.5,-4){\makebox(0,0){$\bar{\iota}$}}
\put(-7.6,-4){\makebox(0,0){$\iota$}}

\put(60,20){\arc{10}{3.143}{0}} \put(60,20){\arc{20}{3.143}{0}}
\path(65,20)(65,0)\path(55,20)(55,0)\path(70,20)(70,0)
\path(50,20)(50,0)

\path(45,30)(40,0)\path(40,30)(35,0)
\path(75,30)(80,0)\path(80,30)(85,0)

\put(35,-4){\makebox(0,0){$\bar{\iota}$}}
\put(40,-4){\makebox(0,0){$\iota$}}
\put(50,-4){\makebox(0,0){$\bar{\iota}$}}
\put(55,-4){\makebox(0,0){$\iota$}}
\put(65,-4){\makebox(0,0){$\bar{\iota}$}}
\put(70,-4){\makebox(0,0){$\iota$}}
\put(80,-4){\makebox(0,0){$\bar{\iota}$}}
\put(85,-4){\makebox(0,0){$\iota$}}

\put(40,34){\makebox(0,0){$\bar{\iota}$}}
\put(45,34){\makebox(0,0){$\iota$}}
\put(75,34){\makebox(0,0){$\bar{\iota}$}}
\put(80,34){\makebox(0,0){$\iota$}}

\end{picture}
\end{center}
\caption{The Jones product from $N\subset M_1$}
\label{jonesproduct}
\end{figure}
Combining with previous remarks,
this means that we can identify the $N$-$M_{ab^{opp}}$
system ${\NXM}_{{a{b}^{opp}}}$  with the $M_{a}$-$M_{b}$ system
$\MaXMb$.

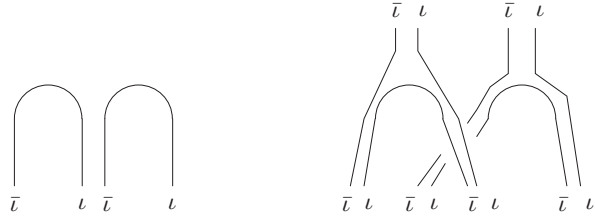
\begin{figure}[htb]
\begin{center}
\unitlength 0.6mm
\begin{picture}(90,45)

\put(-10,10){\arc{15}{3.143}{0}}
\put(-30,10){\arc{15}{3.143}{0}}
\path(-17.5,10)(-17.5,-5)\path(-2.5,10)(-2.5,-5)
\path(-22.5,10)(-22.5,-5)\path(-37.5,10)(-37.5,-5)

\put(-37.5,-9){\makebox(0,0){{\footnotesize $\bar{\iota}$}}}
\put(-22.5,-9){\makebox(0,0){{\footnotesize $\iota$}}}
\put(-17,-9){\makebox(0,0){{\footnotesize $\bar{\iota}$}}}
\put(-2.5,-9){\makebox(0,0){{\footnotesize $\iota$}}}

\put(50,10){\arc{15}{3.143}{0}}
\put(75,10){\arc{15}{3.143}{0}}
\path(57.5,10)(63,-5)\path(65,-5)(61,10)(52,25)(52,30)
\path(42.5,10)(40,-5)\path(37,-5)(40,10)(47,25)(47,30)

\path(82.5,10)(85,-5)
\path(67.5,10)(65,6)\path(58,0)(55,-5)
\path(52,-5)(57,2)
\path(63,9)(65,12)(68,17)(72,20)(72,30)
\path(88,-5)(85,15)(78,20)(78,30)

\put(36,-9){\makebox(0,0){{\footnotesize $\bar{\iota}$}}}
\put(41,-9){\makebox(0,0){{\footnotesize $\iota$}}}
\put(50,-9){\makebox(0,0){{\footnotesize $\bar{\iota}$}}}
\put(55.5,-9){\makebox(0,0){{\footnotesize $\iota$}}}

\put(64,-9){\makebox(0,0){{\footnotesize $\bar{\iota}$}}}
\put(69,-9){\makebox(0,0){{\footnotesize $\iota$}}}
\put(84,-9){\makebox(0,0){{\footnotesize $\bar{\iota}$}}}
\put(90,-9){\makebox(0,0){{\footnotesize $\iota$}}}

\put(47,34){\makebox(0,0){{\footnotesize $\bar{\iota}$}}}
\put(53,34){\makebox(0,0){{\footnotesize $\iota$}}}
\put(72,34){\makebox(0,0){{\footnotesize $\bar{\iota}$}}}
\put(79,34){\makebox(0,0){{\footnotesize $\iota$}}}

\end{picture}
\end{center}
\caption{The braided product}
\label{braidproduct}
\end{figure}

The second interesting case $\Theta_a$=$\Theta_b$ is a Q-system
from a braided subfactor $N\subset M$. Let us consider how the
product $\Theta \otimes \Theta^{opp}$ is related or not to the
Jones basic construction \cite{J1}. If $\Theta=(\theta,v,w)$ is
Q-system for a braided subfactor $N\subset M$ then we can
understand a Q-system for the Jones basic construction $N\subset
M\subset M_1$. If $M\subset M_1$ is the Jones basic construction
for $N\subset M$ and $\iota_1:M\subset M_1$ is the canonical
inclusion map, we can naturally identify the dual canonical
endomorphism $\bar{\iota}_1\iota_1$ of $M\subset M_1$ with the
canonical endomorphism $\iota\bar{\iota}$ of $N\subset M$. Thus
the dual canonical endomorphism for $N\subset M_1$ is
$\overline{\iota_1\iota}\iota_1\iota=
\bar{\iota}\bar{\iota}_1\iota_1\iota=
\bar{\iota}\iota\bar{\iota}\iota=\theta^2$. The corresponding
isometries are $V=wv$ and $W=\theta(wv)$, i.e.\ we have a new
Q-system $\Theta_J$ associated with the Jones basic construction
$\Theta_J=\Big(\theta^2,wv,\theta(wv)\Big)$, see
\fig{jonesproduct}. On the other hand we also have another
Q-system $\Theta _\Phi$ associated with $\theta^2$ from the $\Phi$
product: namely, $\Theta_\Phi=\Big(\theta^2,
\theta(v)v,\theta(\epsilon_{(\theta,\theta)})\theta^2(w)w\Big)$
where $\epsilon$ is the braiding operator, see \fig{braidproduct}.
Denote by $N\subset M_{J}$ and $N\subset M_{\Phi}$ the braided
inclusions from $\Theta_J$ and $\Theta_\Phi$, respectively.

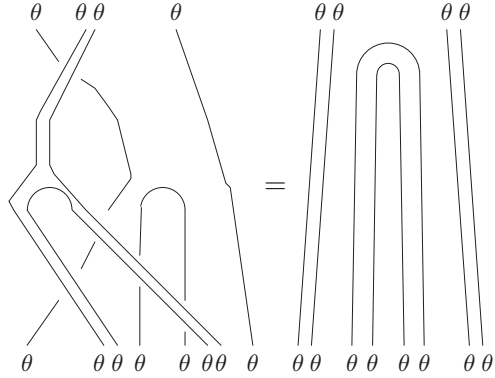
\begin{figure}[htb]
\begin{center}
\unitlength 0.6mm
\begin{picture}(120,80)

\put(25,30){\arc{10}{3.14}{0}} \put(50,30){\arc{10}{3.3}{0}}
\path(20,30)(40,0) \path(30,30)(60,0)

\path(63,0)(33,30)(26,38)(25,40)(25,50)(26,52)(35,70)
\path(37,0)(17,30)(16,32)(22,40)(22,50)(23,52)(33,70)

\path(45.3,30.6)(45,20) \path(45,13)(45,0) \path(55,30)(55,10)
\path(55,3)(55,0)

\path(20,0)(27,10) \path(32,17)(35,23)
\path(38,30)(43,37)(43,38)(40,50)(38,53)(35,57)(32,59)
\path(27,63)(22,70)

\path(70,0)(65,35)(64,36)(60,50)(53,70)

\put(36,74){\makebox(0,0){{\footnotesize $\theta$}}}
\put(32,74){\makebox(0,0){{\footnotesize $\theta$}}}
\put(22,74){\makebox(0,0){{\footnotesize $\theta$}}}
\put(53,74){\makebox(0,0){{\footnotesize $\theta$}}}
\put(63,-4){\makebox(0,0){{\footnotesize $\theta$}}}
\put(36,-4){\makebox(0,0){{\footnotesize $\theta$}}}
\put(45,-4){\makebox(0,0){{\footnotesize $\theta$}}}
\put(55,-4){\makebox(0,0){{\footnotesize $\theta$}}}
\put(20,-4){\makebox(0,0){{\footnotesize $\theta$}}}
\put(70,-4){\makebox(0,0){{\footnotesize $\theta$}}}
\put(40,-4){\makebox(0,0){{\footnotesize $\theta$}}}
\put(60,-4){\makebox(0,0){{\footnotesize $\theta$}}}

\put(75,35){\makebox(0,0){$=$}}

\path(80,0)(85,70) \path(83,0)(88,70)
\put(100,60){\arc{14}{3.14}{0}} \put(100,60){\arc{5}{3.14}{0}}
\path(93,60)(92,0) \path(107,60)(108,0) \path(102.5,60)(103.5,0)
\path(97.5,60)(96.5,0)

\path(113,70)(118,0)%
\path(116,70)(121,0)

\put(85,74){\makebox(0,0){{\footnotesize $\theta$}}}
\put(89,74){\makebox(0,0){{\footnotesize $\theta$}}}
\put(113,74){\makebox(0,0){{\footnotesize $\theta$}}}
\put(117,74){\makebox(0,0){{\footnotesize $\theta$}}}
\put(80,-4){\makebox(0,0){{\footnotesize $\theta$}}}
\put(84,-4){\makebox(0,0){{\footnotesize $\theta$}}}
\put(92,-4){\makebox(0,0){{\footnotesize $\theta$}}}
\put(108,-4){\makebox(0,0){{\footnotesize $\theta$}}}
\put(103.5,-4){\makebox(0,0){{\footnotesize $\theta$}}}
\put(96.5,-4){\makebox(0,0){{\footnotesize $\theta$}}}
\put(118,-4){\makebox(0,0){{\footnotesize $\theta$}}}
\put(122,-4){\makebox(0,0){{\footnotesize $\theta$}}}



\end{picture}
\end{center}
\caption{Equivalent Q-systems
$u\theta^4(u)\theta^2(\epsilon{(\theta^2,\theta^2)})\theta^4(\theta(wv))
\theta(wv)u^\ast=\theta^2(\theta(wv)wv)$} \label{equivalente}
\end{figure}

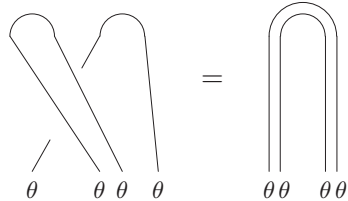
\begin{figure}[htb]
\begin{center}
\unitlength 0.6mm
\begin{picture}(100,55)


\put(0,30){\arc{10}{3.14}{0}} \put(20,30){\arc{10}{3.14}{0}}
\path(-5,30)(15,0) \path(5,30)(20,0) \path(25,30)(28,0)

\path(15,30)(11,23) \path(4,7)(0,0)

\put(40,20){\makebox(0,0){$=$}}

\put(60,30){\arc{10}{3.14}{0}} \put(60,30){\arc{15}{3.14}{0}}
\path(65,30)(65,0) \path(55,30)(55,0) \path(67.5,30)(67.5,0)
\path(52.5,30)(52.5,0)

\put(15,-4){\makebox(0,0){{\footnotesize $\theta$}}}
\put(20,-4){\makebox(0,0){{\footnotesize $\theta$}}}
\put(28,-4){\makebox(0,0){{\footnotesize $\theta$}}}
\put(0,-4){\makebox(0,0){{\footnotesize $\theta$}}}
\put(65,-4){\makebox(0,0){{\footnotesize $\theta$}}}
\put(56,-4){\makebox(0,0){{\footnotesize $\theta$}}}
\put(68.5,-4){\makebox(0,0){{\footnotesize $\theta$}}}
\put(52.5,-4){\makebox(0,0){{\footnotesize $\theta$}}}

\end{picture}
\end{center}
\caption{Equivalent Q-systems $u\theta^2(wv)wv=\theta(wv)wv$}
\label{eqqui2}
\end{figure}

For each Q-system on $\theta^2$, we consider
$\theta^4\in\Sigma(\NXN)$ with two  different products arising
from the Jones or braided products. If we label
the strings $\iota$ and $\bar{\iota}$ by $\theta$ in Figs.\
\ref{jonesproduct} and \ref{braidproduct}, the two Q-systems on
$\theta^4$ are equivalent using the unitary
$u=\epsilon_{(\theta,\theta^2)}$ arising from the braiding, see
\fig{equivalente}, and then we pull the strings to obtain the
corresponding Jones co-product, see also \fig{eqqui2}. We
actually have four Q-systems on $\theta^4$. The first two arise
from the Q-system $\Theta_J$ by performing the Jones and
braided products, thus obtaining $\Theta_{JJ}$
and $\Theta_{J\Phi}$, respectively. The other two $\Theta_{\Phi
J}$ and $\Theta_{\Phi\Phi}$ arise similarly from the Q-system
$\Theta_\Phi$. Namely,
\begin{eqnarray*}
\Theta_{JJ}&=&\Big(\theta^4, \theta(wv)wv,\theta^2(\theta(wv)wv)\Big),\\
\Theta_{J\Phi}&=&\Big(\theta^4,\theta^2(wv)wv,
\theta^2(\epsilon_{(\theta^2,\theta^2)})\theta^4(\theta(wv))\theta(wv)\Big),\\
\Theta_{\Phi J}&=&\bigg(\theta^4,
\theta(\epsilon_{(\theta,\theta)})\theta^2(w)w\theta(v)v, \theta^2
\Big(\theta(\epsilon_{(\theta,\theta)})\theta^2(w)w\theta(v)v\Big)\bigg),\\
\Theta_{\Phi\Phi}&=&\Big(\theta^4,\theta^2(\theta(v)v)\theta(v)v,
\theta^2(\epsilon_{(\theta^2,\theta^2)})\theta^4
(\theta(\epsilon_{(\theta,\theta)})\theta^2(w)w)
\theta(\epsilon_{(\theta,\theta)})\theta^2(w)w \Big).
\end{eqnarray*}
>From the above Q-systems $\Theta_{JJ}$, $\Theta_{J\Phi}$,
$\Theta_{\Phi J}$ we denote by $N\subset M_{JJ}$, $N\subset
M_{J\Phi}$, $N\subset M_{\Phi J}$ and $N\subset M_{\Phi\Phi}$,
respectively, the associated braided inclusions.

\begin{figure}[htb]

\begin{center}
\unitlength 0.6mm
\begin{picture}(105,50)

\put(-20,30){\arc{10}{3.14}{0}} \put(0,30){\arc{10}{3.14}{0}}
\path(-25,30)(-25,0) \path(-15,30)(-5,0) \path(5,30)(4,0)
\path(-5,30)(-9,23) \path(-13,12)(-19,0)

\path(55,35)(55,0) \path(59,35)(59,0)
\put(70,30){\arc{10}{3.14}{0}} \put(90,30){\arc{10}{3.14}{0}}
\path(65,30)(65,0) \path(75,30)(85,0) \path(95,30)(94,0)
\path(85,30)(81,23) \path(77,12)(71,0) \path(100,35)(100,0)
\path(103,35)(103,0)

\put(-25,-4){\makebox(0,0){{\footnotesize $\theta$}}}
\put(-5,-4){\makebox(0,0){{\footnotesize $\theta$}}}
\put(4,-4){\makebox(0,0){{\footnotesize $\theta$}}}
\put(-19,-4){\makebox(0,0){{\footnotesize $\theta$}}}
\put(55,-4){\makebox(0,0){{\footnotesize $\theta$}}}
\put(59,-4){\makebox(0,0){{\footnotesize $\theta$}}}
\put(65,-4){\makebox(0,0){{\footnotesize $\theta$}}}
\put(85,-4){\makebox(0,0){{\footnotesize $\theta$}}}
\put(94,-4){\makebox(0,0){{\footnotesize $\theta$}}}
\put(71,-4){\makebox(0,0){{\footnotesize $\theta$}}}

\put(100,-4){\makebox(0,0){{\footnotesize $\theta$}}}
\put(104,-4){\makebox(0,0){{\footnotesize $\theta$}}}
\put(100,39){\makebox(0,0){{\footnotesize $\theta$}}}
\put(104,39){\makebox(0,0){{\footnotesize $\theta$}}}
\put(55,39){\makebox(0,0){{\footnotesize $\theta$}}}
\put(60,39){\makebox(0,0){{\footnotesize $\theta$}}}

\end{picture}
\end{center}
\caption{The isometries  $V$ and $W$ in $\Theta_{FJ}$}
\label{FJ}
\end{figure}
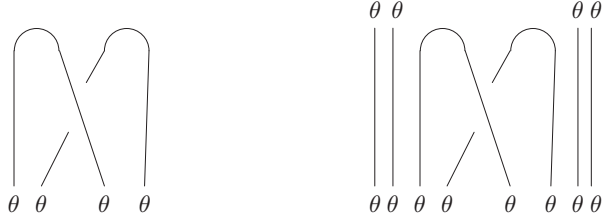

\begin{proposition}
Let $\theta=\bar{\iota}\iota$ be the dual canonical endomorphism
of a braided subfactor $\iota: N\subset M$. Then $\Theta_{JJ}$,
$\Theta_{\Phi J}$ and $\Theta_{J\Phi}$ are equivalent.
\lablth{thetaJF}
\end{proposition}
\bproof That $\Theta_{JJ}$ and $\Theta_{J\Phi}$ are equivalent
Q-systems has been proven above (cf.\ Figs.\ \ref{equivalente}
and \ref{eqqui2}). We can also prove that $\Theta_{\Phi J}$ and
$\Theta_{JJ}$ are equivalent using the unitary
$u=\epsilon_{(\theta,\theta)}$, see \fig{FJ}. \eproof

In the Q-system $\Theta_\Phi$, we can use the relative braiding
\cite{BEK1} and, e.g., use the unitary
$u=\epsilon_{(\iota,\theta)}$ to get the LHS of
\fig{equivalente}, but then we see that the labelling of this
picture is $\bar{\iota} \bar{\iota}\iota\iota\cdots$ whereas the
one for the Jones product is
$\bar{\iota}\iota\bar{\iota}\iota\cdots$ see \fig{jonesproduct},
and they do not match.

\begin{remark}{\rm
Let $\Theta=(\theta,v,w)$ be a Q-system on a braided system $\NXN$
such that $\theta=\la_0\oplus\la_\sigma$ and
$\la_\sigma^2=\la_0$. Thus we can endow $\theta^2$ with two
Q-systems $\Theta_J$ and $\Theta_\Phi$. We see that dim${\cal
Z}(\Theta_J)=1$ whereas dim${\cal Z}(\Theta_\Phi)$=1 if
$\epsilon_{(\sigma,\sigma)}=\kappa_\sigma$id with
$\kappa_\sigma=e^{2\pi ih_\sigma}\not=1$, and dim${\cal
Z}(\Theta_\Phi)$=2 if $\epsilon_{(\sigma,\sigma)}=\hbox{id}$.
In the ${\mathit{SU}}(2)$ modular invariants $Z_{D_n}$ associated
to the Dynkin diagrams $D_n$ \cite{CIZ}, the endomorphism
$\theta=\la_0\oplus\la_\sigma$ is a dual canonical endomorphism
\cite{BE3} producing $Z_{D_n}$. Moreover $\kappa_\sigma=1$ for
$D_{\rm odd}$ and $\kappa_\sigma\not=1$ for $D_{\rm even}$.
Therefore $\Theta_J$ and $\Theta_\Phi$ are inequivalent in the
$D_{\rm even}$ case. In fact
$M_\Phi$ is not a factor, dim$(M_\Phi^\prime\cap M_\Phi)$=$2$,
whereas $M_J$ is (always) a factor, cf.\ Lemma \ref{frobrecip}.
Note that $Z_{N\subset M_\Phi}=Z_{D_{\rm even}}^2$. Moreover,
$\Theta_{\Phi\Phi}$ and $\Theta_{JJ}$ are also not equivalent and
so $N\subset M_{\Phi\Phi}$ and $N\subset M_{JJ}$ are different
inclusions.
 }
\end{remark}

\subsection{Fusion of sufferable modular invariants}
\label{frsuff}

Let $\Theta=(\theta,v,w)$ be a Q-system on a factor $N$,
and $N\subset M$ the
associated inclusion. We will denote Hom(id,$\theta)$ by
$H(\Theta)$ \cite{longo2} and let ${\cal Z}(\Theta):=\{b\in
H(\Theta): w^\ast b=w^\ast\theta(b)\}$, see RHS of Fig.\
\ref{zatop}. If the multiplicity of ${\rm id}$ in $\theta$ is
one, the inclusion $N\subset M$ is irreducible $N^\prime \cap
M\simeq\bbC$ therefore $N\subset M$ is a subfactor. In case the
multiplicity of ${\rm id}$ in $\theta$ is not one, then by \cite{longo}
there is still an inclusion $N\subset M$ with $M$ a von Neumann
algebra with finite dimensional centre. In the following  we get a
condition in terms of data on
$N$, for which
$M$ becomes a factor, thereby yielding a subfactor $N\subset M$.

\begin{figure}[htb]
\begin{center}
\unitlength 0.6mm
\begin{picture}(120,40)

\put(-10,10){\arc{10}{3.14}{0}}
\put(-10,19){\makebox(0,0){{\footnotesize $b$}}}
\put(-10,15){\makebox(0,0){{\footnotesize $\bullet$}}}
\put(-5,5){\makebox(0,0){{\footnotesize $\iota$}}}
\put(-15,6){\makebox(0,0){{\footnotesize $\bar{\iota}$}}}

\put(10,10){\makebox(0,0){{\footnotesize $\in H(\Theta)$}}}

\put(60,20){\arc{8}{3.14}{0}} \put(68,20){\arc{8}{0}{3.14}}
\path(72,20)(75,30) \put(60,24){\makebox(0,0){{\footnotesize
$\bullet$}}} \path(56,20)(57,0) \path(68,0)(70,8)(79,25)(80,30)

\put(85,15){\makebox(0,0){{\footnotesize $=$}}}

\put(106,20){\arc{8}{0}{3.14}} \put(114,20){\arc{8}{3.14}{0}}
\put(114,24){\makebox(0,0){{\footnotesize $\bullet$}}}
\path(118,20)(115,0) \path(102,20)(99,30)
\path(104,0)(103,8)(101,11)(95,25)(94,30)

\put(57,-4){\makebox(0,0){{\footnotesize $\bar{\iota}$}}}
\put(68,-4){\makebox(0,0){{\footnotesize $\iota$}}}
\put(104,-4){\makebox(0,0){{\footnotesize $\bar{\iota}$}}}
\put(115,-4){\makebox(0,0){{\footnotesize $\iota$}}}
\put(75,34){\makebox(0,0){{\footnotesize $\bar{\iota}$}}}
\put(80,35){\makebox(0,0){{\footnotesize $\iota$}}}
\put(94,34){\makebox(0,0){{\footnotesize $\bar{\iota}$}}}
\put(99,34){\makebox(0,0){{\footnotesize $\iota$}}}

\put(60,28){\makebox(0,0){{\footnotesize $b$}}}
\put(114,28){\makebox(0,0){{\footnotesize $b$}}}

\end{picture}
\end{center}
\caption{The algebra ${\cal Z}(\Theta)$} \label{zatop}
\end{figure}
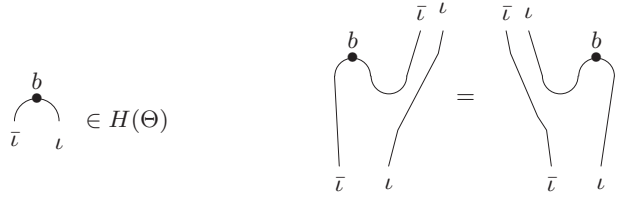

\begin{figure}[htb]
\begin{center}
\unitlength 0.6mm
\begin{picture}(90,30)

\put(-10,10){\makebox(0,0){{\footnotesize ${\cal
L}_{\bar{\iota}}:$}}}

\put(10,10){\arc{10}{3.14}{0}}
\put(10,19){\makebox(0,0){{\footnotesize $b$}}}
\put(10,15){\makebox(0,0){{\footnotesize $\bullet$}}}
\put(15,5){\makebox(0,0){{\footnotesize $\iota$}}}
\put(5,6){\makebox(0,0){{\footnotesize $\bar{\iota}$}}}

\put(30,10){\makebox(0,0){{\footnotesize $\longrightarrow$}}}

\put(50,10){\arc{10}{3.14}{0}} \put(60,10){\arc{10}{0}{3.14}}
\path(45,10)(45,0) \path(65,10)(65,20)

\put(45,-4){\makebox(0,0){{\footnotesize $\iota$}}}
\put(65,24){\makebox(0,0){{\footnotesize $\iota$}}}
\put(62,0){\makebox(0,0){{\footnotesize $b^\ast$}}}
\put(60,5){\makebox(0,0){{\footnotesize $\bullet$}}}

\end{picture}
\end{center}
\caption{left Frobenius reciprocity ${\cal L}_{\bar{\iota}}$:
Hom$(\hbox{id},\theta)\to\hbox{Hom}(\iota,\iota)$}
\label{LFR}
\end{figure}
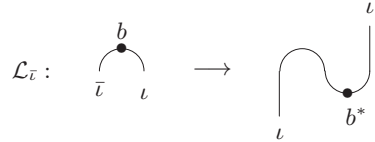

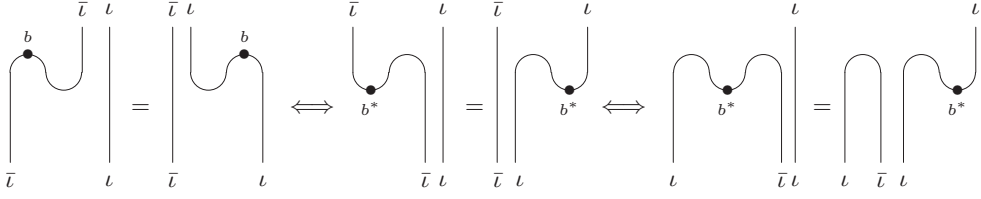
\begin{figure}[htb]
\begin{center}
\unitlength 0.6mm
\begin{picture}(100,45)

\put(-55,20){\arc{8}{3.14}{0}} \put(-47,20){\arc{8}{0}{3.14}}
\path(-37,0)(-37,30)
\path(-59,20)(-59,0) \path(-43,20)(-43,30)
\put(-55,24){\makebox(0,0){{\footnotesize $\bullet$}}}
\put(-55,28){\makebox(0,0){{\tiny $b$}}}

\put(-30,12){\makebox(0,0){{\footnotesize $=$}}}

\path(-23,0)(-23,30)
\put(-15,20){\arc{8}{0}{3.14}}
\put(-7,20){\arc{8}{3.14}{0}}
\path(-19,20)(-19,30) \path(-3,20)(-3,0)
\put(-7,24){\makebox(0,0){{\footnotesize $\bullet$}}}
\put(-7,28){\makebox(0,0){{\tiny $b$}}}

\put(-59,-4){\makebox(0,0){{\footnotesize $\bar{\iota}$}}}
\put(-37,-4){\makebox(0,0){{\footnotesize $\iota$}}}
\put(-23,-4){\makebox(0,0){{\footnotesize $\bar{\iota}$}}}
\put(-3,-4){\makebox(0,0){{\footnotesize $\iota$}}}

\put(-43,34){\makebox(0,0){{\footnotesize $\bar{\iota}$}}}
\put(-37,34){\makebox(0,0){{\footnotesize $\iota$}}}
\put(-23,34){\makebox(0,0){{\footnotesize $\bar{\iota}$}}}
\put(-19,34){\makebox(0,0){{\footnotesize $\iota$}}}

\put(8,12){\makebox(0,0){{\footnotesize $ \Longleftrightarrow $}}}

\put(21,20){\arc{8}{0}{3.14}}
\put(29,20){\arc{8}{3.14}{0}}
\path(37,0)(37,30)
\path(17,20)(17,30) \path(33,20)(33,0)
\put(21,16){\makebox(0,0){{\footnotesize $\bullet$}}}
\put(21,12){\makebox(0,0){{\tiny $b^\ast$}}}

\put(44,12){\makebox(0,0){{\footnotesize $=$}}}

\put(57,20){\arc{8}{3.14}{0}} \put(65,20){\arc{8}{0}{3.14}}
\path(49,0)(49,30)
\path(53,20)(53,0) \path(69,20)(69,30)
\put(65,16){\makebox(0,0){{\footnotesize $\bullet$}}}
\put(65,12){\makebox(0,0){{\tiny $b^\ast$}}}

\put(33,-4){\makebox(0,0){{\footnotesize $\bar{\iota}$}}}
\put(37,-4){\makebox(0,0){{\footnotesize $\iota$}}}
\put(49,-4){\makebox(0,0){{\footnotesize $\bar{\iota}$}}}
\put(54,-4){\makebox(0,0){{\footnotesize $\iota$}}}

\put(17,34){\makebox(0,0){{\footnotesize $\bar{\iota}$}}}
\put(37,34){\makebox(0,0){{\footnotesize $\iota$}}}
\put(49,34){\makebox(0,0){{\footnotesize $\bar{\iota}$}}}
\put(69,34){\makebox(0,0){{\footnotesize $\iota$}}}

\put(77,12){\makebox(0,0){{\footnotesize $\Longleftrightarrow $}}}

\put(92,20){\arc{8}{3.14}{0}} 
\put(100,20){\arc{8}{0}{3.14}} \put(108,20){\arc{8}{3.14}{0}}

\path(88,20)(88,0) \path(112,20)(112,0)
\path(115,0)(115,30)
\put(100,16){\makebox(0,0){{\footnotesize $\bullet$}}}
\put(100,12){\makebox(0,0){{\tiny $b^\ast$}}}

\put(121,12){\makebox(0,0){{\footnotesize $=$}}}

\put(130,20){\arc{8}{3.14}{0}}
\put(143,20){\arc{8}{3.14}{0}} \put(151,20){\arc{8}{0}{3.14}}
\path(126,20)(126,0)\path(134,20)(134,0)
\path(139,20)(139,0)\path(155,20)(155,30)
\put(151,16){\makebox(0,0){{\footnotesize $\bullet$}}}
\put(151,12){\makebox(0,0){{\tiny $b^\ast$}}}

\put(88,-4){\makebox(0,0){{\footnotesize $\iota$}}}
\put(112,-4){\makebox(0,0){{\footnotesize $\bar{\iota}$}}}
\put(115,-4){\makebox(0,0){{\footnotesize $\iota$}}}
\put(126,-4){\makebox(0,0){{\footnotesize $\iota$}}}
\put(134,-4){\makebox(0,0){{\footnotesize $\bar{\iota}$}}}
\put(139,-4){\makebox(0,0){{\footnotesize $\iota$}}}

\put(115,34){\makebox(0,0){{\footnotesize $\iota$}}}
\put(155,34){\makebox(0,0){{\footnotesize $\iota$}}}

\end{picture}
\end{center}
\caption{Proving $b\in{\cal Z}(\Theta) \Longleftrightarrow {\cal
L}_{\bar{\iota}}(b)w_1= w_1{\cal L}_{\bar{\iota}}(b)$ }
\label{diagprove}
\end{figure}

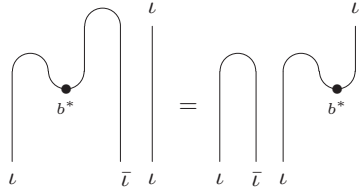
\begin{figure}[htb]
\begin{center}
\unitlength 0.6mm
\begin{picture}(120,35)

\put(30,20){\arc{8}{3.14}{0}} \put(38,20){\arc{8}{0}{3.14}}
\path(42,20)(42,30) \put(38,16){\makebox(0,0){{\footnotesize
$\bullet$}}} \path(26,20)(26,0)
\put(38,12){\makebox(0,0){{\tiny $b^\ast$}}}
\put(46,30){\arc{8}{3.14}{0}} \path(50,30)(50,0)

\path(57,30)(57,0)

\put(65,12){\makebox(0,0){{\footnotesize $=$}}}

\put(76,20){\arc{8}{3.14}{0}} \path(80,20)(80,0)
\path(72,20)(72,0)

\put(90,20){\arc{8}{3.14}{0}} \put(98,20){\arc{8}{0}{3.14}}
\path(102,20)(102,30) \put(98,16){\makebox(0,0){{\footnotesize
$\bullet$}}} \path(86,20)(86,0)
\put(98,12){\makebox(0,0){{\tiny $b^\ast$}}}

\put(102,34){\makebox(0,0){{\footnotesize $\iota$}}}

\put(57,34){\makebox(0,0){{\footnotesize ${\iota}$}}}

\put(26,-4){\makebox(0,0){{\footnotesize $\iota$}}}
\put(51,-4){\makebox(0,0){{\footnotesize $\bar{\iota}$}}}
\put(57,-4){\makebox(0,0){{\footnotesize ${\iota}$}}}
\put(72,-4){\makebox(0,0){{\footnotesize $\iota$}}}
\put(80,-4){\makebox(0,0){{\footnotesize $\bar{\iota}$}}}
\put(86,-4){\makebox(0,0){{\footnotesize $\iota$}}}

\end{picture}
\end{center}
\caption{Diagrammatic representation for $ {\cal
L}_{\bar{\iota}}(b)w_1= w_1{\cal L}_{\bar{\iota}}(b)$ }
\label{lbw1}
\end{figure}

\begin{lemma}
Let $\Theta=(\theta,v,w)$ be a Q-system with $N\subset M$ the
associated inclusion. Then $H(\Theta)\simeq N^\prime\cap M$ and
moreover ${\cal Z}(\Theta)\simeq M^\prime\cap M$.
\lablth{frobrecip}
\end{lemma}
\bproof As in e.g.\ \cite{BEK1} the left Frobenius reciprocity
map ${\cal L}_{\bar{\iota}}:$
Hom$(\hbox{id},\theta)\to\hbox{Hom}(\iota,\iota)$ is defined as
${\cal L}_{\bar{\iota}}(b)$=$d_\iota \iota(b^\ast)w_1$ with
$w_1\in\hbox{Hom}(\hbox{id},\iota\bar{\iota})$, see Fig.\
\ref{LFR}. Note that Hom$(\iota,\iota)=N^\prime\cap M$. Since
$M=Nw_1$ pointwise \cite{longo2} and ${\cal L}_{\bar{\iota}}(b)$
commutes with $N$, we need to prove that $b\in{\cal Z}(\theta)$
if and only if ${\cal L}_{\bar{\iota}}(b)w_1=w_1{\cal
L}_{\bar{\iota}}(b)$. The later equation is displayed in
\fig{lbw1}. The algebraic verification of
\fig{diagprove} is as follows.
Let $b\in H(\Theta)$ such that ${\cal L}_{\bar{\iota}}(b)w_1=w_1
{\cal L}_{\bar{\iota}}(b)w_1$, i.e.\  $\iota(b^\ast)w_1
w_1=w_1\iota(b^\ast)w_1$, thus applying $\bar{\iota}$ on both
sides, $\bar{\iota}\iota(b^\ast)w w=w\bar{\iota}\iota(b^\ast)w$.
Multiplication by $v^\ast$ on the LHS yields
$v^\ast\theta(b^\ast)w^2=v^\ast w\theta(b^\ast)$. But since
$\Theta$ is a Q-system, $v^\ast w=1/d$ where $d=d(\iota)$ and
$v^\ast\theta(b^\ast)=b^\ast v^\ast$ because $v\in\hbox{Hom}({\rm
id,\theta})$. Therefore we get $b^\ast v^\ast w^2=
\theta(b^\ast)w/d$, so again by definition of Q-system and then
applying the $\ast$-operation on both sides we get $w^\ast
b=w^\ast \theta(b)$. Therefore $b\in{\cal Z}(\Theta)$.
Let us now assume that $b\in{\cal Z}(\Theta)$, thus applying
$\ast$ on both sides and then $\iota$, we get
$\iota(b^\ast)\iota(w)=\iota \theta(b^\ast)\iota(w)$.
Multiplication by $w_1$ on the LHS of both sides yields
$\iota(b^\ast)\iota(w)w_1=\iota \theta(b^\ast)\iota(w)w_1$.  We
have
$\iota(b^\ast)\iota(w)w_1=\iota(b^\ast)\gamma(w_1)w_1=\iota(b^\ast)w_1w_1$
because $w_1\in\hbox{Hom}({\rm id},\gamma)$ and $\iota
\theta(b^\ast)\iota(w)w_1=\gamma(\iota(b^\ast)w_1)w_1=w_1\iota(b^\ast)w_1$
using again that $w_1\in\hbox{Hom}({\rm id},\gamma)$. Hence we
conclude that $\iota(b^\ast)w_1w_1=w_1\iota(b^\ast)w_1$, thus
${\cal L}_{\bar{\iota}}(b)w_1=w_1 {\cal L}_{\bar{\iota}}(b)w_1$.
\eproof

Let $\Theta$=$(\theta,v,w)$ be a Q-system, $N\subset M$ the
associated braided inclusion and $p$ a projection in the finite
dimensional algebra $M^\prime\cap
M\subset\hbox{Hom}(\iota,\iota)$. Then we consider the cut-down
Q-system $\Theta_p$ where $\theta_p=\bar{\iota}_p\iota_p$ with
$\iota_p(x)=p\iota(x)$. Note that $\theta_p=\bar{\iota} p
\iota\prec \bar{\iota}\iota\in\Sigma(\NXN)$. Denote by $Z_p$ the
modular invariant produced by $\Theta_p$. Then if we take a
(finite) partition of unity $\sum_p p={\bf 1}$ of the algebra
$M^\prime\cap M$ by minimal projections $\{p\}$,
$\Theta=\bigoplus_p\Theta_p$.

\begin{lemma}
Let $N\subset M$ be a braided inclusion and $\{p\}$ be a
partition of unity by minimal projections of $M^\prime\cap M$.
Then $Z_{N\subset M}$ is a sum of normalised sufferable modular invariants
$Z_{N\subset M}=\sum_p Z_p$, so that $Z_{0,0}$ is the dimension of
$M^\prime\cap M$. \lablth{sumsmod}
\end{lemma}
\bproof Since every projection $p$ is minimal, $p(M^\prime\cap
M)p=\bbC p$ so $M_p$ is a factor and then $[Z_p]_{00}=1$.
Observe by \cite[Sect. 5]{E1} that $\a$-induction from $N$ to
$M$ is trivial on the relative commutant $N^\prime\cap M$, thus
trivial on $M^\prime\cap M$, so that $\a$-induction is the sum of
the cut-down $\a$-inductions. Consequently, Z is additive
on Q-systems. For a formal proof see the description
in the remark before Theorem 5.7 in \cite{BEK1} of $Z$
in terms of $N$-$N$ and $N$-$M$ data, from which
additivity is clear, or use \cite[Proposition 5.3]{frs}.\eproof

\begin{theorem}
Let  $N\subset M_a$ and $N\subset M_b$ be braided inclusions,
whose modular invariants we denote by $Z_a$ and $Z_b$
respectively.
Then the product $Z_aZ_b^{t}$ is sufferable and decomposes as an
integral sum of normalised sufferable modular invariants.
\lablth{nice}
\end{theorem}
\bproof
Suppose $\Theta_a$ and $\Theta_b$ are Q-systems for
the inclusions $N\subset M_a$ and $N\subset M_b$ respectively.
Let $\Theta_{ab^{\rm opp}}$ denote the Q-system
$\Theta_a\otimes\Theta_b^{{\rm opp}}$
with $N\subset M_{ab^{opp}}$ the associated inclusion. Then
by Lemma \ref{categ} and \cite[Proposition 5.3]{frs}, $Z_{N\subset
M_{ab^{{\rm opp}}}} = Z_aZ_b^t$. On the other hand, let us take a partition of
unity of minimal projections $\{p\}$ of $M_{ab^{opp}}^\prime\cap
M_{ab^{{\rm opp}}}$. Then $Z_{N\subset M_{ab^{\rm opp}}} =\sum Z_p$ by Lemma
\ref{sumsmod} and the result is proven. \eproof

A special case of this result when $a=b$, and  $(Z_aZ_a^{t})_{00} = 2^Ë$
was derived in \cite{E1}.

\begin{corollary}
The set of sufferable modular invariants yields a fusion rule algebra.
\end{corollary}

\subsection{Q-systems from $6j$-symbols}

Here we present in co-ordinate form  the Q-system relations on the
isometries $v$ and $w$. We consider only the case where the
Verlinde fusion rules are multiplicty free and where $\theta$ is
also multiplicity free i.e.\ $\lan\theta, \la \ran=0,1$.

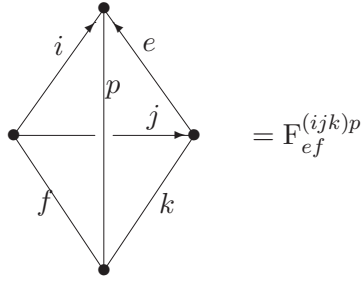
\begin{figure}[htb]
\begin{center}
\unitlength 0.6mm
\begin{picture}(70,70)
\thinlines \path(40,30)(20,58)(0,30)(20,0)(40,30)
\path(20,58)(20,0) \path(0,30)(18,30) \path(22,30)(40,30)

\put(40,30){\makebox(0,0){$\bullet$}}
\put(20,58){\makebox(0,0){$\bullet$}}
\put(0,30){\makebox(0,0){$\bullet$}}
\put(20,0){\makebox(0,0){$\bullet$}}

\put(30,50){\makebox(0,0){$e$}} \put(31,34){\makebox(0,0){$j$}}
\put(10,50){\makebox(0,0){$i$}} \put(7,15){\makebox(0,0){$f$}}
\put(34,15){\makebox(0,0){$k$}} \put(22,40){\makebox(0,0){$p$}}
\put(38,30){\vector(1,0){0}} \put(22,55){\vector(-1,1){0}}
\put(18,55){\vector(1,1){0}}
\put(65,30){\makebox(0,0){$=\hbox{F}_{ef}^{(ijk)p}$}}
\end{picture}
\end{center}
\caption{$6j$-symbols} \label{symbols}
\end{figure}

Let $s_i\in\hbox{Hom}(\la_i,\theta)$ be Cuntz generators, i.e.\
$s_i^\ast s_j=\delta_{i,j}{\bf 1}$, $\sum_i s_i s_i^\ast={\bf 1}$,
such that $\theta(a)=\sum s_i \la_i (a)s_i^\ast$. Similarly fix an
orthonormal basis $\{v_{ij}^k\}$ of Hom$(\la_k,\la_i\la_j)$. Then
any $w\in\hbox{Hom}(\theta,\theta^2)$ can be written as $w$=$\sum
w_{ij}^k s_i\la_i(s_j)v_{ij}^ks_k^\ast$ with $w_{ij}^k$ complex
numbers. Since we assume dim Hom$(id,\theta)=1$, we can take
$v=s_0$. We  note by \cite{KI}, that the Q-relations
$w^\ast\theta(w)=ww^\ast$ and $w^2=\theta(w)w$ are equivalent as
long as we assume the other relation $dv^\ast
w=dw^\ast\theta(v)={\rm id}$. So $w$ has to be an isometry $w^\ast
w={\rm id}$ and $dv^\ast w=w^\ast\theta(v)={\rm id}$ and
$w^2=\theta(w)w$. These constraints are respectively
\begin{eqnarray}
\sum_{i,j}|w_{ij}^k|^2=1,
w_{0k}^k=w_{k0}^k=1/d,\  \hbox{and}\
w_{ij}^lw_{lk}^p=\sum_{e\prec\theta}w_{jk}^ew_{ie}^p
{\rm F}_{el}^{(ijk)p}
\label{mult}
\end{eqnarray}
for all $i,j,l,k,p\prec\theta$
and where the numbers $F_{ef}^{(ijk)p}\in\bbC$ are the associated
$6j$-symbols \cite{EK}, see \fig{symbols}.
We  calculate
\begin{eqnarray*}
w^\ast w&=&\sum\Big(w_{rp}^t s_r\la_r(s_p) v_{rp}^t s_t^\ast\Big)^\ast
\Big(w_{ij}^k s_i \la_i(s_j) v_{ij}^k s_k^\ast \Big)\\
&=&\sum \bar{w}_{rp}^t w_{ij}^k s_t
v_{rp}^{t\ast}\la_r(s_p^\ast)s_r^\ast\cdot s_i\la_i(s_j) v_{ij}^k
s_k^\ast=\sum |w_{ij}^t|^2 s_ts_t^\ast,
\end{eqnarray*}
hence $\sum_{i,j}|w_{ij}^t|^2=1$ for all $t$. Also
$$d v^\ast w=d s_0^\ast\sum w_{ij}^k s_i\la_i(s_j)v_{ij}^ks_k^\ast=d \sum
w_{0j}^k s_jv_{0j}^ks_k^\ast=\sum w_{0k}^k s_ks_k^\ast,$$
hence $w_{0k}^k=1/d$ for all $k$ and similarly from $d
w^\ast\theta(v)={\rm id}$ we obtain $w_{k0}^k=1/d$ for all $k$.
We now have
\begin{eqnarray*}
ww&=&\sum w_{ij}^k w_{pq}^r s_i\la_i(s_j) v_{ij}^k s_k^\ast\cdot
s_p\la_p(s_q) v_{pq}^rs_r^\ast\\
&=&\sum w_{ij}^k w_{kq}^r s_i\la_i(s_j) \big(v_{ij}^k \la_k(s_q)\big)
v_{kq}^r s_r^\ast\\
&=& \sum w_{ij}^k w_{kq}^r s_i\la_i(s_j) \la_i\la_j(s_q)v_{ij}^k v_{kq}^r
s_r^\ast
\end{eqnarray*}
whereas
\begin{eqnarray*}
\theta(w)w&=& \sum w_{ij}^k w_{pq}^r\theta\Big(s_i\la_i(s_j) v_{ij}^k
s_k^\ast \Big)\cdot s_p\la_p(s_q) v_{pq}^rs_r^\ast\\
&=&\sum w_{ij}^k w_{pq}^r
s_p\la_p(s_i)\la_p\la_i(s_j)\la_p(v_{ij}^k)\la_p(s_k^\ast)\la_p(s_q)
v_{pq}^r s_r^\ast\\
&=& \sum w_{ij}^k w_{pk}^r s_p\la_p(s_i)\la_p\la_i(s_j)\la_p(v_{ij}^k)
v_{pk}^r s_r^\ast\\
&=&\sum w_{ij}^k w_{pk}^r {\rm
F}_{kf}^{(pij)r}s_p\la_p(s_i)\la_p\la_i(s_j) v_{pi}^f v_{fj}^r s_r^\ast
\end{eqnarray*}
yielding the final relation of Eqs. (\ref{mult}).

In the sequel we will use the following normalisations
\cite[\eerf{2.38}]{frs}
\begin{eqnarray}
{\rm F}_{pj}^{(0jk)p}={\rm F}_{ki}^{(i0k)p}=\quad {\rm
F}_{jp}^{(ij0)p}=1 \label{norm}
\end{eqnarray}
whenever they are allowed to be non-zero by the fusion rules.
It is usually tremendous work to prove that a given
$\theta\in\Sigma(\NXN)$ is a dual canonical endomorphism
(requiring the solution of
non-linear equations). Here we construct Q-systems that will be
of use in the sequel. The Frobenius-Schur indicator
FS$_\la=\omega^{-1}\sum_{\mu,\nu}N_{\mu,\nu}^\la d_\mu d_\nu
\omega_\mu^2/\omega_\nu^2$ vanishes unless $\la$ is self-conjugate
and is $+1$ or $-1$ depending on whether $\la$ is real or
pseudo-real, respectively (see e.g. \cite{frs}). For
a simple current $\sigma$ (i.e.\ $d_\sigma=1$) one has
FS$_\sigma=\omega^2_\sigma$. There is a well known obstruction
for finding a periodic inner perturbation of a given automorphism
$\sigma$ which is
a priori only periodic modulo the inner automorphisms $[\sigma^n]=[\hbox{id}]$,
so that $(Ad(v)\sigma)^n= \hbox{id}$ \cite{alain1}.
We consider this obstruction in the context of Q-systems.

\begin{lemma} Let $\NXN$ be a nondegenerate braided system and
suppose that there is a simple current $\sigma\in\NXN$ such that
$[\sigma^2]=[\hbox{id}]$. Then $\theta=\hbox{id}\oplus\sigma$ is
a dual canonical endomorphism of a braided subfactor if and only
if the Frobenius-Schur indicator ${\rm FS}_\sigma=1$.
\lablth{newfrob2}
\end{lemma}
\bproof The last Eq.\ of (\ref{mult})
reduces to
$$w_{ij}^lw_{lk}^p=w_{jk}^0w_{i0}^p {\rm F}_{0l}^{(ijk)p}+
w_{jk}^1w_{i1}^p \rm{F}_{1l}^{(ijk)p}$$ where we put $1=\sigma$,
$0=\hbox{id}$ and $i,j,k,l=0,1$. We
 know the 6$j$-symbols by e.g. \cite[Eqs. (2.38)
and (2.47)]{frs}:
$${\rm F}_{p,j}^{(0jk)p}=
{\rm F}_{ki}^{(i0k)p}={\rm F}_{jp}^{(ij0)p}=1,\quad
\hbox{{FS}}_\sigma=\rm{F}_{00}^{(111)1}.$$  When we
insert these 6$j$-symbols into the above equation we see that we
get consistent constants $w_{ij}^k$ if and only if
FS$_\sigma=1$. We then get the following solution:
$$ w_{00}^0 = w_{11}^0 = w_{01}^1=w_{10}^1=1/\sqrt{2}\quad
\hbox{and}\quad w_{11}^1 = w_{00}^1 = w_{01}^0 = w_{10}^0=0,$$
which yield a solution for all of  Eqs.\ (\ref{mult}).
\eproof
We remark that FS$_\sigma=1$ is Rehren's
condition \cite{rehren} $\omega^2_\sigma=1$ in the net setting.
Let $p$ be a prime number and $\NXN=\{\la_{(m,n)}\}$ be the
system of endomorphisms of $N=M_0\rtimes\Delta(\bbZ_p)$ obtained
from the quantum double of $M_0\subset M_0\rtimes\bbZ_p$ whose
sectors obey the $\bbZ_p\times\bbZ_p$ fusion rules and where
$M_0$ is a type III factor. We also fix a non-degenerate braiding
on $\NXN$. Now we look for $\theta\in\Sigma(\NXN)$ such that
$\theta$ is a dual canonical endomorphism.
\begin{lemma}
(i)\ Let $p\not=2$ be a prime number and $H$ a subgroup of
$\bbZ_p\times\bbZ_p$ isomorphic to $\bbZ_p$ and consider
$\theta_H=\bigoplus_{(m,n)\in H}\la_{(m,n)}$. Then $\theta_H$ is
a dual canonical endomorphism of some braided subfactor $N\subset
M$.\\
(ii) The endomorphism $\theta_{\bbZ_p\times\bbZ_p}$ is a dual
canonical endomorphism. \lablth{newfrob}
\end{lemma}
\bproof The 6j-symbols are, see e.g. \cite[\eerf{7.340}]{fk},
$${\rm F}_{j+k,i+j}^{(ijk)i+j+k}=1,\quad i,j,k\in
\bbZ_p\times\bbZ_p$$ and the others vanish. Note also this
normalization is consistent with the one used in \erf{norm}. Next
we consider the subfactor $N=M_0\rtimes \Delta(\bbZ_p)\subset
M_0\rtimes (\bbZ_p\times\bbZ_p)$. Then the $\theta_{{\scriptsize
\Delta}(\bbZ_p)}$=$(0,0)\oplus(1,1)\oplus\cdots\oplus(p,p)$ is its
dual canonical endomorphism by \cite{KY1} and therefore
$\theta_{{\scriptsize \Delta}(\bbZ_p)}$ is a dual canonical
endomorphism. The 6j-symbols for $\Delta(\bbZ_p)$ are obtained by
restriction of those for $\bbZ_p\times\bbZ_p$, i.e.
$${\rm F}_{j+k,i+j}^{(ijk)i+j+k}=1,\quad i,j,k\in
\Delta(\bbZ_p)$$ and ${\rm F}_{m,n}^{(ijk)h}=0$ if $h\not=i+j+k$
or $m\not=j+k$ or $n\not=i+j$ with $i,j,k,h,m,n$ in
$\Delta(\bbZ_p)$. Hence if we take another copy $H$ of $\bbZ_p$ in
$\bbZ_p\times\bbZ_p$ we get exactly the same formula for its
6j-symbols, i.e.
$${\rm F}_{j+k,i+j}^{(ijk)i+j+k}=1,\quad i,j,k\in H$$
and ${\rm F}_{m,n}^{(ijk)h}=0$ if $h\not=i+j+k$ or $m\not=j+k$ or
$n\not=i+j$ with $i,j,k,h,m,n$ in $H$ (because of the fusion
rules). Thus solving \erf{mult} for $\theta_H$ is equivalent to
solving the same equation for $\theta_{{\scriptsize
\Delta}(\bbZ_p)}$. Since we already know that such a solution
exists for $\Delta(\bbZ_p)$ we conclude that $\theta_H$ is also a
dual canonical endomorphism of some braided subfactor $N\subset
M$.

Let us now prove (ii). From part (i), $\theta_{\bbZ_p\times 0}$
and $\theta_{0\times\bbZ_p}$ are dual canonical endomorphisms, so
is their product $\theta_{\bbZ_p\times\bbZ_p}$. \eproof

\section{The subfactor $M_0\rtimes\Delta(G)\subset M_0
\rtimes(G\times G)$}
\label{kosakietal}

We identify the intermediate subfactors of the asymptotic
subfactor of a finite group $G$ that arise from subgroups of $G$.

\begin{figure}[htb]
\begin{center}
\unitlength 0.6mm
\begin{picture}(65,50)
\thinlines

\put(10,34){\makebox(0,0){${\widehat G}$}}
\put(10,-5){\makebox(0,0){${\widehat H}$}}
\put(8,2){\vector(0,1){25}}
\put(12,28){\vector(0,-1){25}}

{\scriptsize
\put(17,17){\makebox(0,0){rest}}
\put(3,17){\makebox(0,0){ind}}
\put(60,34){\makebox(0,0){$\matrix{\leftarrow\ M-M \hbox{ system},\cr
 \hbox{q\ dim=dim}(\pi), \pi\in{\widehat G}}$ }}
\put(60,-5){\makebox(0,0){$\matrix{\leftarrow M-N \hbox{ system},\cr
{}\quad \hbox{q\ dim=dim}(\pi)\sqrt{|G/H|}, \pi\in{\widehat H}}$}}
}
\end{picture}
\end{center}
\caption{Dual principal graph of $M_0\rtimes H\subset M_0\rtimes G$}
\label{dualgraph}
\end{figure}
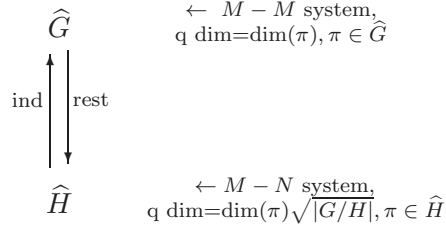

\begin{proposition}
There is a bijection between normal subgroups $N\triangleleft G$
and subgroups $H$ of $G\times G$ containing $\Delta(G)$.
\lablth{inter-sub}
\end{proposition}
\bproof Let $H$ be a subgroup of $G\times G$ such that
$\Delta\subset H,$ where $\Delta$ denotes the diagonal subgroup of $G\times G$.
If $(h_1,h_2)\in H$, $(h_1^{-1}, h_1^{-1})\in \Delta\subset H$ implies
$(e,h_2h_1^{-1})\in H$. Then $H^\prime:=\{h\in G: (e,h)\in H\}$ is a
subgroup of $G$. On the other hand, for $s\in G,\ h\in H^\prime$,
$(e,shs^{-1})=(s,s)(e,h)(s^{-1},s^{-1})\ \in H,$
therefore $shs^{-1}\in H^\prime$ for all $s$ and so $H^\prime$ is a normal
subgroup of $G$ (and $\{e\}\times H^\prime\subset H$). Conversely, take a
normal subgroup $H^\prime$ of $G$, and define 
$H:=\Delta\cdot H^\prime\times H^\prime=
\{(gH^\prime,gH^\prime): g\in G\}$.
Then $\Delta\subset H\subset G\times G$ and clearly the above
$H\leftrightarrow H^\prime$ construction is a bijection.
\eproof
Then from each normal $N$ subgroup of $G$ we get an intermediate subfactor
$M_0\rtimes\Delta(G)\subset M_0\rtimes N \subset M_0\rtimes(G\times G)$.
Clearly $\Delta(G)$ is normal in 
$G\times G$ if and only if $G$ is an abelian group.
In the more general case where we fix a subgroup $H$ of $G$, the
principal and dual graphs of $N=M_0\rtimes H\subset M_0\rtimes G=M$
are computed in \cite{KY1} through the Mackey induction/restriction
machinery. The dual principal graph is as in \fig{dualgraph} with
$\hbox{res}: \pi\to\pi\mid_H,\quad \hbox{ind}:\sigma\to
\hbox{ind}_H^G\sigma$ being restriction and induction of group
representations
respectively.

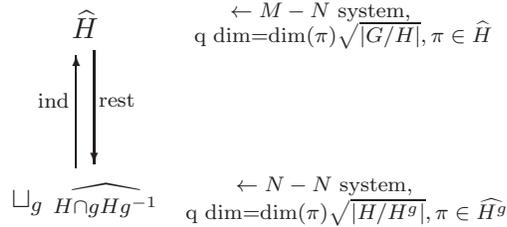
\begin{figure}[htb]
\begin{center}
\unitlength 0.6mm
\begin{picture}(65,45)
\thinlines

\put(10,34){\makebox(0,0){${\widehat H}$}}
\put(10,-5){\makebox(0,0){$\sqcup_g\ \widehat{{}_{H\cap gHg^{-1}}}$}}
\put(8,2){\vector(0,1){25}}
\put(12,28){\vector(0,-1){25}}

{\scriptsize
\put(17,17){\makebox(0,0){rest}}
\put(3,17){\makebox(0,0){ind}}
\put(63,34){\makebox(0,0){$\matrix{\leftarrow M-N \hbox{ system},\cr
{}\qquad \hbox{q\ dim=dim}(\pi)\sqrt{|G/H|}, \pi\in{\widehat H}}$ }}
\put(63,-5){\makebox(0,0){$\matrix{\leftarrow N-N \hbox{ system},\cr
{}\qquad \hbox{q\ dim=dim}(\pi)\sqrt{|H/ H^g|}, \pi\in{\widehat {H^g}}}$}}
}

\end{picture}
\end{center}
\caption{Principal graph of $M_0\rtimes H\subset M_0\rtimes G$}
\label{pdualgraph}
\end{figure}

Strictly speaking, the dual principal graph is the connected
component of the above restriction/induction procedure containing
the trivial representation of $G$. The principal graph is more
subtle. First we take representatives of the double cosets
$H\setminus G/ H$. Then for each double coset $HgH$ we fix a
representative $g$ and consider the stabilizer subgroup
$H^g:=H\cap gHg^{-1}$. Then the set of vertices of the principal
graph of $N\subset M$ is labelled by pairs $(g,\sigma)$ where $g$
runs over a set of representatives of the double cosets and
$\sigma$ are irreducible representations of $H^g$ \cite{KY}.
Finally the principal graph is the connected component containing
the trivial representation of $H$ of the graph, see
\fig{pdualgraph}. By \cite[Proposition 31]{KY1}, the vertex set of
the principal graph of $N\subset M$ is given by $\widehat{G/{\cal
N}}\sqcup\widehat{H/{\cal N}}$ where ${\cal N}$ is the largest
normal subgroup of $G$ containing $H$.
Furthermore, $\widehat{G/{\cal N}}$ labels the $N$-$N$
irreducible sectors whereas $\widehat{H/{\cal N}}$ labels the
$N$-$M$ ones. We note that the principal graph of $M_0\rtimes
H\subset M_0\rtimes G$ is the dual principal graph of the inclusion
$M_0^G\subset M_0^H$ (see e.g. \cite{J1} or \cite{EK}). The more
difficult fusion of the $N$-$N$ sectors is also carried out by
Kosaki et al. \cite{KY1,KY}.

When we are in the situation $\Delta(G)\subset G\times G$,
the center ${\cal Z}_G $ of $G$ is trivial if and only if ${\cal
N}$ is trivial. Therefore, the above Mackey induction/restriction
graphs are connected if and only if ${\cal Z}_G=\{e\}$. This
happens, for instance, when $G$ is $S_3$ the symmetric group on 3
letters. Moreover, since
$$\Delta(G)(g,h)\Delta(G)=\Delta(G)(e,g^{-1}h)\Delta(G)$$
for any $g,h \in G$, we easily see that every double coset
$\Delta(G)(g,h)\Delta(G)$ gives rise to a conjugacy
class $C_{g^{-1}h}$ of $\Delta(G)$. Since the stabilizer of
$(g,h)$ equals the centraliser of $g^{-1}h$, i.e.
$$(g,h)\Delta(G)(g^{-1},h^{-1})\cap\Delta(G)=\{x\in G:
xg^{-1}h=g^{-1}hx \},$$
we have the group theoretical relation between the framework of Coste
et al. \cite{Gan} and the one in \cite{KY1} which we will use in the sequel.

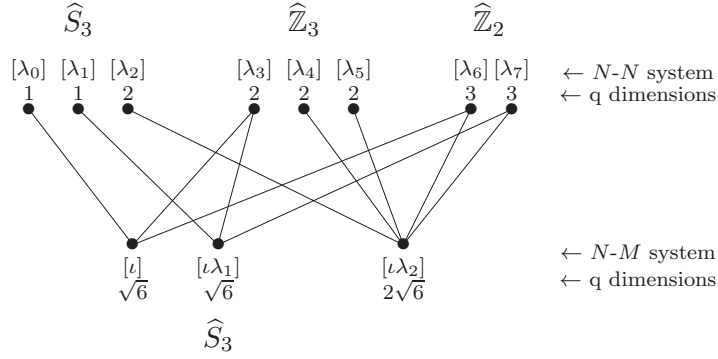
\begin{figure}[htb]
\begin{center}
\unitlength 0.6mm
\begin{picture}(80,85)
\thinlines \put(-7,20){\makebox(0,0){$\bullet$}}
\put(12,20){\makebox(0,0){$\bullet$}}
\put(53,20){\makebox(0,0){$\bullet$}}

\put(-30,50){\makebox(0,0){$\bullet$}}
\put(-19,50){\makebox(0,0){$\bullet$}}
\put(-8,50){\makebox(0,0){$\bullet$}}
\put(20,50){\makebox(0,0){$\bullet$}}
\put(31,50){\makebox(0,0){$\bullet$}}
\put(42,50){\makebox(0,0){$\bullet$}}
\put(68,50){\makebox(0,0){$\bullet$}}
\put(77,50){\makebox(0,0){$\bullet$}}

\path(-30,50)(-7,20)(20,50)(12,20)(77,50)(53,20)(-8,50)

\path(-7,20)(68,50) \path(12,20)(-19,50)  \path(53,20)(68,50)
\path(53,20)(42,50) \path(53,20)(31,50)

\put(-19,70){\makebox(0,0){${\widehat S_3}$}}
\put(31,70){\makebox(0,0){${\widehat \bbZ_3}$}}
\put(72,70){\makebox(0,0){${\widehat \bbZ_2}$}}
\put(12,0){\makebox(0,0){${\widehat S_3}$}}

{\scriptsize \put(-7,12){\makebox(0,0){$ \matrix {[\iota]\cr
\sqrt{6}}$}}
\put(12,12){\makebox(0,0){$\matrix{[\iota\la_1]\cr\sqrt{6}}$}}
\put(53,12){\makebox(0,0){$\matrix{[\iota\la_2]\cr 2\sqrt{6}}$}}

\put(-30,56){\makebox(0,0){$ \matrix {[\la_0]\cr 1}$}}
\put(-19,56){\makebox(0,0){$ \matrix {[\la_1]\cr 1}$}}
\put(-8,56){\makebox(0,0){$ \matrix {[\la_2]\cr 2}$}}
\put(20,56){\makebox(0,0){$ \matrix {[\la_3] \cr 2}$}}
\put(31,56){\makebox(0,0){$ \matrix {[\la_4]\cr 2}$}}
\put(42,56){\makebox(0,0){$ \matrix {[\la_5]\cr 2}$}}
\put(68,56){\makebox(0,0){$ \matrix {[\la_6] \cr 3}$}}
\put(77,56){\makebox(0,0){$ \matrix {[\la_7]\cr 3}$}}

\put(105,58){\makebox(0,0){$\leftarrow$ $N$-$N$ system}}
\put(105,18){\makebox(0,0){$\leftarrow$ $N$-$M$ system}}
\put(105,53){\makebox(0,0){$\leftarrow$ q\ dimensions}}
\put(105,12){\makebox(0,0){$\leftarrow$ q\ dimensions}} }
\end{picture}
\end{center}
\caption{Principal graph of $M_0\rtimes \Delta(S_3)\subset
M_0\rtimes(S_3\times S_3)$} \label{prgraph}
\end{figure}

We now illustrate the above algorithm in some concrete examples.
In the first example we display the principal graph of
$N=M_0\rtimes\Delta(S_3)\subset M_0\rtimes(S_3\times S_3)=M$
in \fig{prgraph}. Here we have three
$S_3\equiv\Delta(S_3)$ double cosets:\\

$\bullet$ $S_3(e,e)S_3$ with 6 elements,

$\bullet$ $S_3((132),e)S_3$ with 12 elements and

$\bullet$ $S_3((132),(12))S_3$ with 18 elements. \\

We further have
$$S_3\cap(e,e)S_3{(e,e)}^{-1}\simeq S_3,\quad
S_3\cap((132),e)S_3{(132,e)}^{-1}\simeq\bbZ_3$$
$$S_3\cap((132),(12))S_3{((132),(12))}^{-1}\simeq\bbZ_2.$$
The three first vertices $[\la_0],[\la_1], [\la_2]$ of
\fig{prgraph} label the trivial $\bf 1$, sign $\epsilon$ and 2
dimensional $\sigma$ irreducible representations of the group
$S_3$, respectively. The vertices $[\la_3], [\la_4],[\la_5]$
label the trivial and the other two irreducible representations
of $\bbZ_3$, respectively. Finally, the last two vertices
$[\la_6], [\la_7]$ label the two irreducible inequivalent
representations of $\bbZ_2$. In the sequel, it will become clear
why we choose this notation for the vertices of the principal
graph of the above subfactor $N\subset M$.

\begin{figure}[htb]
\begin{center}
\unitlength 0.6mm
\begin{picture}(80,90)
\thinlines

\put(0,20){\makebox(0,0){$\bullet$}}
\put(20,20){\makebox(0,0){$\bullet$}}
\put(40,20){\makebox(0,0){$\bullet$}}

\put(-35,50){\makebox(0,0){$\bullet$}}
\put(-20,50){\makebox(0,0){$\bullet$}}
\put(-5,50){\makebox(0,0){$\bullet$}}
\put(10,50){\makebox(0,0){$\bullet$}}

\put(25,50){\makebox(0,0){$\bullet$}}
\put(40,50){\makebox(0,0){$\bullet$}}
\put(55,50){\makebox(0,0){$\bullet$}}
\put(70,50){\makebox(0,0){$\bullet$}}
\put(85,50){\makebox(0,0){$\bullet$}}

\path(-35,50)(0,20)(25,50)
\path(-20,50)(20,20)(10,50)
\path(-5,50)(40,20)(85,50)(0,20)
\path(20,20)(85,50)
\path(40,50)(40,20)(55,50)
\path(40,20)(70,50)

\put(20,-5){\makebox(0,0){$\widehat{S_3}$}}
\put(18,75){\makebox(0,0){$\widehat{S_3}\times\widehat{S_3}$}}

{\scriptsize \put(0,10){\makebox(0,0){$\matrix{[\bar{\iota}]\cr
\sqrt{6}}$}}
\put(20,10){\makebox(0,0){$\matrix{[\la_1\bar{\iota}]\cr
\sqrt{6}}$}}
\put(40,10){\makebox(0,0){$\matrix{[\la_2\bar{\iota}]\cr
2\sqrt{6}}$}}

\put(-35,59){\makebox(0,0){$\matrix{[{\bf 1},{\bf 1}]\cr 1}$}}
\put(-20,59){\makebox(0,0){$\matrix{[{\bf 1},\epsilon]\cr 1}$}}
\put(-5,59){\makebox(0,0){$\matrix{[{\bf 1},\sigma]\cr 2}$}}
\put(10,59){\makebox(0,0){$\matrix{[\epsilon,{\bf 1}]\cr 1}$}}
\put(25,59){\makebox(0,0){$\matrix{[\epsilon,\epsilon]\cr 1}$}}
\put(40,59){\makebox(0,0){$\matrix{[\epsilon,\sigma]\cr 2}$}}
\put(55,59){\makebox(0,0){$\matrix{[\sigma,{\bf 1}]\cr 2}$}}
\put(70,59){\makebox(0,0){$\matrix{[\sigma,\epsilon]\cr 2}$}}
\put(85,59){\makebox(0,0){$\matrix{[\sigma,\sigma]\cr 2}$}}
\put(110,60){\makebox(0,0){$\leftarrow$ $M$-$M$ system}}
\put(105,10){\makebox(0,0){$\leftarrow$ q\ dimensions}}
\put(110,56){\makebox(0,0){$\leftarrow$ q\ dimensions}}
\put(105,15){\makebox(0,0){$\leftarrow$ $M$-$N$ system}}
}

\end{picture}
\end{center}
\caption{Dual principal graph of $M_0\rtimes 
S_3\subset M_0\rtimes(S_3\times S_3)$}
\label{dprgraph}
\end{figure}
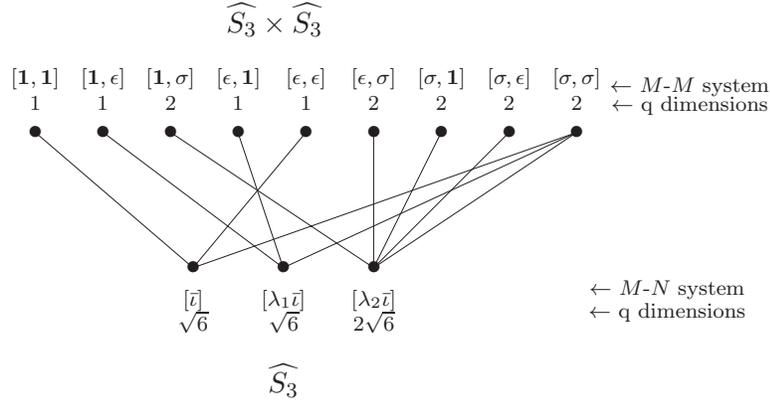

The dual principal graph of the subfactor $M_0\rtimes S_3\subset
M_0\rtimes S_3\times S_3$ is derived from the irreducible
representations of $S_3$ whose (commutative) fusion is as follows:
$\epsilon\epsilon= {\bf 1}=\epsilon{\bf 1},\quad
\sigma\sigma={\bf 1}+\epsilon+\sigma$. In this way the Mackey
induction/restriction between $\widehat{S_3}\times\widehat{S_3}$
and $\widehat{S_3}$ gives rise to the graph displayed in
\fig{dprgraph} which is in turn the dual principal graph of the
above subfactor by \cite{KY1}.

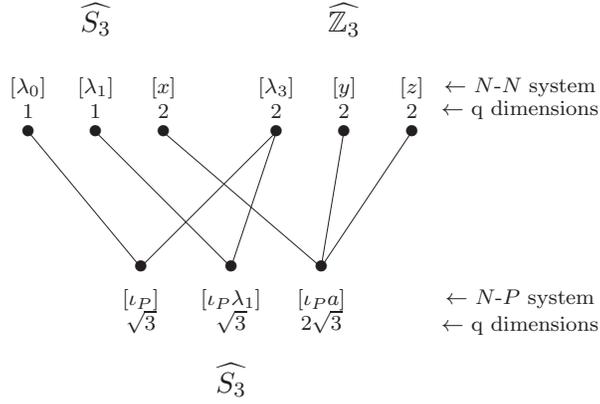
\begin{figure}[htb]
\begin{center}
\unitlength 0.6mm
\begin{picture}(70,80)
\thinlines

\put(10,20){\makebox(0,0){$\bullet$}}
\put(30,20){\makebox(0,0){$\bullet$}}
\put(50,20){\makebox(0,0){$\bullet$}}

\put(-15,50){\makebox(0,0){$\bullet$}}
\put(0,50){\makebox(0,0){$\bullet$}}
\put(15,50){\makebox(0,0){$\bullet$}}
\put(40,50){\makebox(0,0){$\bullet$}}
\put(55,50){\makebox(0,0){$\bullet$}}
\put(70,50){\makebox(0,0){$\bullet$}}

\put(30,-5){\makebox(0,0){$\widehat{S_3}$}}
\put(0,75){\makebox(0,0){$\widehat{S_3}$}}
\put(55,75){\makebox(0,0){$\widehat{\bbZ_3}$}}

\path(-15,50)(10,20)(40,50)(30,20)(0,50)
\path(15,50)(50,20)(55,50)
\path(50,20)(70,50)

{\scriptsize \put(10,10){\makebox(0,0){$ \matrix {[\iota_P]\cr
\sqrt{3}}$}} \put(30,10){\makebox(0,0){$ \matrix
{[\iota_P\la_1]\cr \sqrt{3}}$}} \put(50,10){\makebox(0,0){$
\matrix {[\iota_P a]\cr 2\sqrt{3}}$}}

\put(-15,57){\makebox(0,0){$ \matrix {[\la_0]\cr 1}$}}
\put(0,57){\makebox(0,0){$ \matrix {[\la_1]\cr 1}$}}
\put(15,57){\makebox(0,0){$ \matrix {[x]\cr 2}$}}
\put(40,57){\makebox(0,0){$ \matrix {[\la_3]\cr 2}$}}
\put(55,57){\makebox(0,0){$ \matrix {[y]\cr 2}$}}
\put(70,57){\makebox(0,0){$ \matrix {[z]\cr 2}$}}

\put(94,60){\makebox(0,0){$\leftarrow$ $N$-$N$ system}}
\put(94,13){\makebox(0,0){$\leftarrow$ $N$-$P$ system}}
\put(94,55){\makebox(0,0){$\leftarrow$ q\ dimensions}}
\put(94,7){\makebox(0,0){$\leftarrow$ q\ dimensions}}
}

\end{picture}
\end{center}
\caption{Mackey induction/restriction for 
$S_3\subset S_3\cdot(\bbZ_3\times \bbZ_3)$}
\label{prgraph2}
\end{figure}

The only non-trivial intermediate subfactor $M_0\rtimes S_3\subset
P\subset M_0\rtimes (S_3\times S_3)$ arises from the normal
subgroup $\bbZ_3\equiv\{e,(123),(132)\}$ of $S_3$ by Proposition
\ref{inter-sub}, i.e., $P=M_0\rtimes K$ where
$K:=S_3\cdot(\bbZ_3\times \bbZ_3)\simeq
(\bbZ_3\times\bbZ_3)\rtimes\bbZ_2$. The group $H$ is a
non-commutative subgroup of $S_3\times S_3$ of order 18 (equal to
pairs in $S_3\times S_3$ with the same parity). The subfactor
$ M_0\rtimes S_3\subset M_0\rtimes K$ thus has Jones index equals 3.
Next we compute its principal graph. We obtain two double cosets
for $S_3\subset S_3\cdot(\bbZ_3\times \bbZ_3)$; namely
$S_3(e,e)S_3\simeq S_3$ and $S_3((12),(23))S_3$ with 12 elements.
We then obtain
$$S_3\cap(e,e)S_3{(e,e)}^{-1}\simeq S_3, \quad
S_3\cap((12),(23))S_3{((12),(23))}^{-1}\simeq \bbZ_3.$$ In this
way, we display the induction/restriction graph in \fig{prgraph2}
where the bottom vertices and the first 3 vertices on the top
label the irreducible representations of $S_3$ and the other 3
vertices label the irreducible representations of $\bbZ_3$. We
explain the notation of the labels. The $N$-$N$ labels of the
connected component containing the vertex $[\iota_P]$ in
\fig{prgraph2} is a subset of those $N$-$N$ vertices in
\fig{prgraph}, the others (i.e.\ $x,y,z$) are unrelated to the
latter ones \cite{KY1},  (cf.\ the proof of Proposition
\ref{nicetheta})

Note that the principal graph of
$N= M_0\rtimes S_3\subset M_0\rtimes (S_3\cdot(\bbZ_3\times\bbZ_3))=P$ is
the connected component containing the trivial representation of $S_3$.
Therefore it is the Dynkin diagram $A_5$. This can alternatively be
seen by noting that $S_3$ is non-normal (index 3) subgroup 
of $S_3\cdot(\bbZ_3\times\bbZ_3)$.

\section{Modular data and modular invariants for cyclic groups}
\label{cyclicprime}

The modular data for finite abelian groups have been carried out explicitly
in \cite[Section 3.1]{Gan}. In particular for cyclic groups $\bbZ_d$,
the primary fields are parametrized by pairs $(m,n)$ for $m,n\in\bbZ_d$
whose conjugate is $(-m,-n)$, and the $S$ and $T$ matrices:

\begin{eqnarray}
S_{(m,n),(m^\prime,n^\prime)}&=&
{1\over d}\hbox{exp}[-2\pi\sqrt{-1}(nm^\prime+mn^\prime)/d],\nonumber \\
T_{(m,n),(m,n)}&=&\hbox{exp}[(2\pi\sqrt{-1} nm/d)].
\label{stgroup}
\end{eqnarray}

When $d=p$ is a prime number, the list of all modular invariants
for the cyclic groups $\bbZ_p$ have been computed in \cite{Gan}.
For any prime $p$, there are precisely four non-permutation
modular invariants: ${\cal Z}_4=xx^\ast, {\cal Z}_7=xy^\ast,$
${\cal Z}_8=yx^\ast,$ and ${\cal Z}_5=yy^\ast$ where
$x=\sum_{i=0}^{p-1}\ch_{i0}$ and $y=\sum_{j=0}^{p-1}\ch_{0j}.$
For any prime $p>2$, there will be precisely four permutation
modular invariants: ${\cal Z}_1= \sum_{i,j=0}^{p-1}
\ch_{ij}\ch_{ij}^\ast$, ${\cal
Z}_2=\sum_{i,j=0}^{p-1}\ch_{ij}\ch_{ji}^\ast$, ${\cal Z}_3=
\sum_{i,j=0}^{p-1}\ch_{ij}\ch_{-j,-i}^\ast$ and ${\cal
Z}_6=\sum_{i,j=0}^{p-1}\ch_{ij}\ch_{-i,-j}^\ast$ the charge
conjugation. When $p=2$, ${\cal Z}_6={\cal Z}_1$ and ${\cal
Z}_3={\cal Z}_2$. This is the complete list of modular invariants
for any prime $p$. The modular data for $\bbZ_3$ are
{\footnotesize
$$
S={1\over 3}\pmatrix{1&1&1&1&1&1&1&1&1\cr
1&1&1&\overline{\zeta}&\overline{\zeta}&
\overline{\zeta}
&\zeta&\zeta&\zeta  \cr
 1&1&1&\zeta&\zeta&
\zeta&\overline{\zeta}&
\overline{\zeta}&\overline{\zeta} \cr
1&\overline{\zeta}&\zeta&
1&\overline{\zeta}&\zeta&1&\overline{\zeta}&
\zeta \cr
 1&\overline{\zeta} &\zeta&
\overline{\zeta}&\zeta&
1&\zeta&1&\overline{\zeta} \cr
1&\overline{\zeta}&\zeta&\zeta&1&
\overline{\zeta}&\overline{\zeta}&
\zeta&1 \cr
1&\zeta&\overline{\zeta}&1&\zeta&\overline{\zeta}&1&\zeta&
\overline{\zeta} \cr
1&\zeta&\overline{\zeta}&
\overline{\zeta}&1&\zeta&
\zeta&\overline{\zeta}&1 \cr
1&\zeta&\overline{\zeta}&\zeta&
\overline{\zeta}&1&\overline{\zeta}&1&\zeta },\quad
T=\hbox{diag}(1,1,1,1, \overline{\zeta},
\zeta,1,\zeta,\overline{\zeta}),$$
}
where $\zeta=\exp[2\pi\sqrt{-1}/3]$ and $\overline{\zeta}
=\exp[-2\pi\sqrt{-1}/3].$

\subsection{Subfactors for the quantum $\bbZ_3$ modular invariants}

The $\bbZ_2$ case has been carried out in \cite[Page 25]{BE5},
using a different notation. After relabelling the primary fields
they match in the following manner: ${\cal Z}_2=W$, ${\cal
Z}_4=X_c$, ${\cal Z}_5=X_s$, ${\cal Z}_7=Q^t$ and ${\cal Z}_8=Q$.
In Table \ref{fusionsuff2} we provide all products of modular
invariants of the quantum $\bbZ_2$ double model.
\begin{table}
{\scriptsize
$$\begin{array}{|c|c  c c c c c|}  \hline\hline & 1 & W &
X_s&X_c& Q^t & Q \\ \hline\hline
1& 1&W&X_s&X_c& Q& Q^t \\
W&W&1&Q^t&Q&X_c&X_s\\
X_s&X_s&Q&2X_s&Q&2Q&X_s\\
X_c&X_c&Q^t&Q^t&2X_c&X_c& 2Q^t\\
Q^t&Q^t&X_c&2Q^t&X_c&2X_c&Q^t\\
Q&Q&X_s& X_s &2Q&Q&2X_s\\
\hline
\end{array}$$
} \caption{Fusion ${\cal Z}_a{\cal Z}_b^t$ of $\bbZ_2$ modular
invariants} \label{fusionsuff2}
\end{table}
The heterotic modular invariants ${\cal Z}_7$ and ${\cal Z}_8$
arise from the subfactor $N\subset N\rtimes(\bbZ_2\times\bbZ_2)$
producing another two modular invariants  (namely ${\cal Z}_4$
and ${\cal Z}_5$) through intermediate subfactors by the
application of the type I parent theorem. The last non-trivial
modular invariant arise from an inclusion $N\subset
N\rtimes_{\la_v}\Delta(\bbZ_2)$. We note that here
$N=\pi_0({\mathit{SO}}(16l))^{\prime\prime}$ with $\pi_0$ the level
$1$ vacumm representation \cite{W}.

Let us concentrate on the quantum double of $\bbZ_3$. For each
primary $(m,n)$ of $\bbZ_3:=\{0,1,2\}$ modular data
let $\lambda_{mn}$ denote a corresponding endomorphism of a type
III factor $N$. Let $\NXN$ be the system arising from these
endomorphisms which as sectors obey the $\bbZ_3\times\bbZ_3$
fusion rules \cite{longo2}. The following are the subgroups of
$\bbZ_3\times\bbZ_3$:

$\bullet\quad H_1=\{(0,0)\},$

$\bullet\quad H_2=\{(0,0),(1,2),(2,1)\},$

$\bullet\quad H_3=\{(0,0),(1,1),(2,2)\}=\Delta(\bbZ_3),$

$\bullet\quad H_4=\bbZ_3\times\{0\},$

$\bullet\quad H_5=\{0\}\times\bbZ_3,$

$\bullet\quad H_6=\bbZ_3\times\bbZ_3.$

 Let $N=M_0\rtimes \bbZ_3$
for a fixed type III factor $M_0$ in which we realise the braided
system $\NXN$. By Lemma \ref{newfrob}, every subgroup $H$ of
$\bbZ_3\times\bbZ_3$ gives rise to a dual canonical endomorphism
$\theta_H=\oplus_{(m,n)\in H}\la_{(m,m)}$.

\begin{proposition}
In the case $p=3$, the subfactors $N\subset N\rtimes H_1$,
$N\subset N\rtimes H_2$ and $N\subset N\rtimes H_3$ produce
the permutation modular invariants
${\cal Z}_1$, ${\cal Z}_2$ and ${\cal Z}_3$, respectively.
\lablth{h1h2}
\end{proposition}
\bproof
The dual canonical sector of $N\subset N\rtimes H_2$
is $[\theta_2]=[\lambda_{00}]\oplus[\la_{12}]\oplus[\la_{21}].$
In this way, we obtain $\lan\theta_2\la_{mn}, \lambda_{mn}\ran=1$
for all primary fields of the model. 
On the other hand, $\lan\theta_2\lambda_{0i},\lambda_{0j}\ran=
\delta_{i,j}$ where $\delta$ denotes the Kronecker function.
So we have so far three irreducible $N$-$M$ sectors
$[\iota\lambda_{00}], [\iota\lambda_{01}],[\iota\lambda_{02}]$.
However, since
$\lan\theta_2\lambda_{0i},
\lambda_{kj}\ran=\delta_{i-k,j}$
we have no further irreducible $N$-$M$ sectors.
Then $\Tr(Z_{N\subset M})=3$ from where we conclude that
$Z_{N\subset M}=
{\cal Z}_4, {\cal Z}_5, {\cal Z}_2\ \hbox{or}\ {\cal Z}_3$.
As $\lan\theta_2\lambda_{00},\lambda_{0i}\ran=0$ for
$i=1,2$, $Z_{N\subset M}$ is neither ${\cal Z}_4$ nor ${\cal Z}_5$
because of the inequality
$\lan\alpha_\la^\pm,\alpha_\mu^\pm\ran\leq\lan\theta
\la,\mu\ran$ as in \cite[\eerf{37}]{BEK2}.
Now let us appeal to \cite[Proposition 3.3]{E1}. On one hand,
$$\bigoplus_{i,j,k,r}
{{\cal Z}_3}_{(ij,kr)}[\la_{ij}\bar{\la}_{kr}]=
3[\la_{00}]\oplus 3[\la_{11}]\oplus 3[\la_{22}]$$
and on the other,
$$\bigoplus_{a\in\MXN}[a\bar{a}]=[\bar{\la}_{00}\bar{\iota}\iota\la_{00}]
\oplus[\bar{\la}_{12}\bar{\iota}\iota\la_{12}]\oplus[\bar{\la}_{21}
\bar{\iota}\iota\la_{21}]=3[\la_{00}]\oplus 3[\la_{12}]\oplus 3[\la_{21}]$$
which is not the same sector.
So $Z_{N\subset M}$ cannot be ${\cal Z}_3$ by \cite[Proposition
3.3]{E1},
and must be ${\cal Z}_2$.

The dual canonical sector of $N\subset N\rtimes H_3$ is
$[\theta_3]=[\la_{00}]\oplus[\la_{11}]
\oplus[\lambda_{22}].$
Since for $\theta_3$ we have
$\lan\theta_3\lambda_{0i},\lambda_{kj}
\ran=\delta_{i+k,j}$
we get exactly 3 irreducible $N$-$M$ sectors $[\iota\la_{00}],[\iota\la_{01}],
[\iota\lambda_{02}]$, so
$\hbox{Tr}(Z_{N\subset M})=3$ implying that
$Z_{N\subset M}={\cal Z}_4, {\cal Z}_5,
{\cal Z}_2\ \hbox{or}\ {\cal Z}_3.$
Elimination of both ${\cal Z}_4$ and ${\cal Z}_5$ is similar to the
above case.

We also have
$$\bigoplus_{a\in\MXN}[a\bar{a}]=[\bar{\la}_{00}\bar{\iota}\iota\la_{00}]
\oplus[\bar{\la}_{11}\bar{\iota}\iota\la_{11}]\oplus[\bar{\la}_{22}
\bar{\iota}\iota\la_{22}]=3[\la_{00}]\oplus 3[\la_{11}]\oplus
3[\la_{22}]$$ and on the other hand,
$$\bigoplus_{i,j,k,r}
{{\cal Z}_2}_{(ij,kr)}[\la_{ij}\bar{\la}_{kr}]=
3[\la_{00}]\oplus 3[\la_{12}]\oplus 3[\la_{21}].$$
Therefore again by \cite[Proposition 3.3]{E1},
$Z_{N\subset N\rtimes H_3}={\cal Z}_3$.
\eproof

\begin{proposition}
The quantum double modular invariants ${\cal Z}_4$ and ${\cal
Z}_5$ are realised by the subfactors $N\subset N\rtimes
(\bbZ_3\times\{0\})$ and $N\subset N\rtimes(\{0\}\times\bbZ_3)$.
\end{proposition}

\bproof For the subfactor $N\subset N\rtimes(\{0\}\times\bbZ_3)$,
$[\theta_5]=[\la_{00}]\oplus[\la_{01}]\oplus[\la_{02}]$
is its
dual canonical sector. Computing we find that
$\lan\theta_5\la_{ij},\la_{ij}\ran=1, \hbox{ for all}\ i, j,$
$\lan\theta_5\la_{00},\la_{01}\ran=$
$\lan\theta\la_{00},\la_{02}\ran=$
$\lan\theta_5\la_{10},\la_{11}\ran=$
$\lan\theta_5\la_{10},\la_{12}\ran=$
$\lan\theta_5\la_{20},\la_{21}\ran$
$=\lan\theta\la_{20},\la_{22}\ran=1$.
 The other intertwiner spaces are zero.
Hence $[\iota\la_{00}]=[\iota\la_{01}]=[\iota\la_{02}],
[\iota\la_{10}]=[\iota\la_{11}]=[\iota\la_{12}],
[\iota\la_{20}]=[\iota\la_{21}]=[\iota\la_{22}].$
Then $\Tr(Z_{N\subset M})=3$, implying that
$Z_{N\subset M}={\cal Z}_4, {\cal Z}_5, {\cal Z}_2\ \hbox{or}\ {\cal Z}_3.$
Since $\lan\theta\la_{00},\la_{10}\ran=0$, $Z_{N\subset M}$ cannot be
${\cal Z}_4$.
Now we employ \cite[Proposition 3.3]{E1}.
First, we use $\theta_5$ and obtain:
$$[\epsilon]:=\bigoplus_{a\in\MXN}[a\bar{a}]=
[\bar{\la}_{00}\bar{\iota}\iota\la_{00}]
\oplus[\bar{\la}_{10}\bar{\iota}\iota\la_{10}]\oplus[\bar{\la}_{20}
\bar{\iota}\iota\la_{20}]=3[\la_{00}]\oplus 3[\la_{01}]\oplus 3[\la_{02}].$$
For $Z={\cal Z}_2$ or ${\cal Z}_3$,
$$\bigoplus_{i,j,k,r}Z_{(ij,kr)}[\la_{ij}\bar{\la}_{kr}]$$
has been computed in the proof of Proposition \ref{h1h2}. It is in
either case different from the above computation of $[\epsilon]$.
So $Z_{N\subset M}$ cannot be the permutation modular invariants
${\cal Z}_2$ and ${\cal Z}_3$. Therefore
$Z_{N\subset N\rtimes(\{0\}\times\bbZ_3)}={\cal Z}_5.$
In a similar way, we prove that
$N\subset N\rtimes(\bbZ_3\times\{0\})$
produces ${\cal Z}_4$.
\eproof

\begin{proposition}
The modular invariants ${\cal Z}_6,{\cal Z}_7, {\cal Z}_8$ are
sufferable.
\end{proposition}
\bproof
We already know that ${\cal Z}_2$ and ${\cal Z}_3$ are sufferable,
then so is ${\cal Z}_3{\cal Z}_2={\cal Z}_6$ by Section \ref{frobsub}.
The same with ${\cal Z}_7={\cal Z}_2{\cal Z}_5$.
\eproof

\begin{corollary} For $i=1,2,3,4,5$, the symmetric quantum
$\bbZ_3$ double modular invariant ${\cal Z}_i$ is produced by the
subfactor $N\subset N\rtimes H_i$. Moreover $H_6$ produces ${\cal
Z}_6$, ${\cal Z}_7$ and ${\cal Z}_8$.
\end{corollary}

The global index of $\NXN$ is $\omega=9$. In the
$H_1,H_2,H_3,H_6$ cases, we have $\MXMo=\MXMpm=\MXM$ and as
sectors they are $\bbZ_3\times\bbZ_3$. In the cases
$H_4,H_5,H_7,H_8$, $\MXMo=\{(0,0)\}$, the sectors of $\MXMpm$ are
isomorphic to $\bbZ_3$ and those from $\MXM$ isomorphic to
$\bbZ_3\times\bbZ_3$. The indecomposable nimreps can be derived
by \cite[Proposition 4]{Gan2}. We also have Table
\ref{fusionsuff0} for the (unique) fusion among the (sufferable)
quantum $\bbZ_3$ modular invariants. Note that 
$H^2(\bbZ_3\times\bbZ_3,S^1)$ =$\bbZ_3$
is responsible for $H_6$ producing
three different modular invariants.
\unitlength 1.4mm \thinlines
\begin{table}[htb]
\begin{center}
{\scriptsize
\begin{tabular}{|c|c c c c c c c c|}  \hline\hline
& ${\cal Z}_1$ & ${\cal Z}_2$ & ${\cal Z}_3$ & ${\cal Z}_{4}$&
${\cal Z}_5$ & ${\cal Z}_6$ &${\cal Z}_7$&${\cal Z}_8$\\ \hline\hline
${\cal Z}_1$& ${\cal Z}_1$ &${\cal Z}_2$ &${\cal Z}_3$ &
${\cal Z}_{4}$ & ${\cal Z}_5$ & ${\cal Z}_6$ & ${\cal Z}_8$&${\cal Z}_7$\\
${\cal Z}_2$& ${\cal Z}_2$ &${\cal Z}_1$ &
${\cal Z}_6$ &${\cal Z}_{8}$ &${\cal Z}_7$
&${\cal Z}_3$&${\cal Z}_4$ &${\cal Z}_5$\\
${\cal Z}_3$&${\cal Z}_3$ &${\cal Z}_6$ &${\cal Z}_1$ 
&${\cal Z}_{8}$ &${\cal Z}_7$ &${\cal Z}_2$&${\cal Z}_4$
&${\cal Z}_5$\\
${\cal Z}_4$&${\cal Z}_4$&${\cal Z}_7$ &${\cal Z}_7$&
$3{\cal Z}_4$ &${\cal Z}_7$ & ${\cal Z}_4$ & ${\cal Z}_4$&
$3{\cal Z}_7$\\
${\cal Z}_5$& ${\cal Z}_5$&${\cal Z}_8$ &${\cal Z}_8$ 
&${\cal Z}_8$ &$3{\cal Z}_5$
&${\cal Z}_5$& 3${\cal Z}_8$&${\cal Z}_5$\\
${\cal Z}_6$& ${\cal Z}_6$&${\cal Z}_3$ &${\cal Z}_2$ & ${\cal Z}_4$ & 
${\cal Z}_5$
&${\cal Z}_1$ &${\cal Z}_8$&${\cal Z}_7$\\
${\cal Z}_7$&${\cal Z}_7$&${\cal Z}_4$&${\cal Z}_4$&
${\cal Z}_4$&3${\cal Z}_7$&${\cal Z}_7$&$3{\cal Z}_4$&${\cal Z}_7$\\
${\cal Z}_8$&${\cal Z}_8$&${\cal Z}_5$&
${\cal Z}_5$&$3{\cal Z}_8$&${\cal Z}_5$&${\cal Z}_8$&${\cal Z}_8$&
3${\cal Z}_5$\\
\hline
\end{tabular}
}
\end{center}
\caption{Fusion ${\cal Z}_a{\cal Z}_b^t$ of $\bbZ_3$ level 0 modular
invariants} \label{fusionsuff0}
\end{table}

\subsection{Quantum $\bbZ_3$ double level 1 and 2 modular invariants}

We have that $H^3(\bbZ_3, S^{1})=\bbZ_3=\{0,1,2\}$ by e.g.
\cite{Gan}. For $k=1$ (the level 1 case), all the primary fields
are simple currents again and the $S$ and $T$ matrices are as
follows: {\footnotesize
$$
S={1\over 3}\pmatrix{1&1&1&1&1&1&1&1&1\cr
1&1&1&\bar{z}&\bar{z}&\bar{z}&z&z&z\cr
1&1&1&z&z&z&\bar{z}&\bar{z}&\bar{z}\cr
1&\bar{z}&z&\bar{w}&v&u&\bar{v}&w&\bar{u}\cr
1&\bar{z}&z&v&u&\bar{w}&\bar{u}&\bar{v}&w\cr
1&\bar{z}&z&u&\bar{w}&v&w&\bar{u}&\bar{v}\cr
1&z&\bar{z}&\bar{v}&\bar{u}&w&u&v&\bar{w}\cr
1&z&\bar{z}&w&\bar{v}&\bar{u}&v&\bar{w}&u\cr
1&z&\bar{z}&\bar{u}&w&\bar{v}&\bar{w}&u&v\cr},\quad
T=\hbox{diag}(1,1,1,u,v,\bar{w},v,u,\bar{w}),$$
}
where $u=\exp[2\pi i/9]$, $v=\exp[8\pi i/9]$, $w=\exp[4\pi i/9]$,
$z=\exp[-2\pi i/3]$.

These modular data arise from the affine ${\mathit{SU}}(9)_1$.
The modular invariants are:\\ $Z_1=\sum_{i,j}\ch_{ij}\ch_{ij}^\ast$,
$Z_2=\ch_{00}\ch_{00}^\ast+(\ch_{01}\ch_{02}^\ast+\ch_{10}\ch_{21}^\ast+
\ch_{11}\ch_{20}^\ast+\ch_{12}\ch_{22}^\ast + c.c.)$ and
$Z_3=\sum_{i,j}\ch_{i0}\ch_{j0}^\ast$.

Since $T_{10,10}$ is a 9th root of unity, by \cite[Lemma
4.4]{rehren} we can choose representatives in each simple current
sector such that $\{\la_{ij}\}\simeq\bbZ_9$. Moreover the sector
$[\theta]=\bigoplus_{ij}[\la_{ij}]$ is a dual canonical sector of
the subfactor $N\subset N\rtimes\bbZ_9$. This produces in turn
the above $Z_2$ permutation invariant. Also as
$\{\la_{00},\la_{01},\la_{02}\}$ is a subgroup of $\bbZ_9$,
$[\theta]=[\la_{00}]\oplus[\la_{01}]\oplus[\la_{02}]$ is a dual
canonical endomorphism of $N\subset N\rtimes(\{0\}\times\bbZ_3)$
which produces the above $Z_3$ modular invariant. In this way, we
conclude that all the quantum $\bbZ_3$ double level 1 modular
invariants are sufferable and their fusion rules are as in Table
\ref{fusionsuffsu91}.
\begin{table}
{\scriptsize
$$\begin{array}{|c|c c c|}  \hline\hline
& Z_1 & Z_2 & Z_3  \\ \hline\hline
Z_1& Z_1&Z_2&Z_3\\
Z_2&Z_2&Z_1&Z_3\\
Z_3&Z_3&Z_3&3Z_3\\
\hline
\end{array}$$
} \caption{Fusion $Z_a Z_b^t$ of $\bbZ_3$ levels 1 or 2 modular invariants}
\label{fusionsuffsu91}
\end{table}
The level $k=2$ modular data is obtained by the complex
conjugation of those for the level $k=1$ case.

\section{Subfactors and nimreps for the quantum $S_3$ double modular invariants}
\label{s3}

\subsection{The modular data and the Verlinde fusion matrices}

The symmetric group $S_3=\bbZ_3\rtimes\bbZ_2$ case has eight
primary fields:
$(e,\psi_i)$ for $i=0,1,2$ where $\psi_0,$ $\psi_1,$ $\psi_2$
are the characters associated to the trivial, sign and two
dimensional representation (resp.) of $S_3$; $((123),\varphi_k)$,
$k=0,1,2$ for the 3 characters of $\bbZ_3$; and $((12), \varphi_k^\prime)$,
for the 2 characters of $\bbZ_2$. We label these primary fields
$\Phi$ by $0, 1,\dots,7$ as in \cite{Gan}. The modular data are:
{\footnotesize
\begin{eqnarray}
S&=&{1\over 6}\pmatrix{1&1&2&2&2&2&3&3\cr
1&1&2&2&2&2&-3&-3\cr 2&2&4&-2&-2&-2&0&0\cr
2&2&-2&4&-2&-2&0&0\cr 2&2&-2&-2&-2&4&0&0\cr
2&2&-2&-2&4&-2&0&0\cr 3&-3&0&0&0&0&3&-3\cr
3&-3&0&0&0&0&-3&3\cr}\nonumber \\
T&=&\hbox{diag}(1,1,1,1,\hbox{exp}[2\pi\sqrt{-1}/3],
\hbox{exp}[4\pi\sqrt{-1}/3],1,-1). \nonumber
\end{eqnarray}
}
The primary fields 0 and 1 are the only simple currents of our
model.
Next we determine explicitly all Verlinde fusion matrices
\cite{Ve}, $N_i,$ $i=0,1,\dots,7$ since they will be employed in
the calculation of the dimensions of intertwiner spaces later.
They are computed from the Verlinde formula \cite{Ve} and can
also be derived from \cite[Page 279]{KY}, which as quadratic forms are:
\begin{eqnarray*}
N_0&=&|\ch_0|^2+|\ch_1|^2+|\ch_2|^2+|\ch_3|^2+|\ch_4|^2+|\ch_5|^2+|\ch_6|^2+
|\ch_7|^2,\\
N_1&=&\ch_0\ch_1^\ast+\ch_1\ch_0^\ast+|\ch_2|^2+|\ch_3|^2+|\ch_4|^2+
|\ch_5|^2+\ch_6\ch_7^\ast+\ch_7\ch_6^\ast,\\
N_2&=&\ch_0\ch_2^\ast+\ch_2\ch_0^\ast+\ch_1\ch_2^\ast+
\ch_2\ch_1^\ast+|\ch_2|^2+
\ch_3(\ch_4^\ast+\ch_5^\ast)+(\ch_4+\ch_5)\ch_3^\ast+\ch_4\ch_5^\ast+\\
&&\ch_5\ch_4^\ast+|\ch_6+\ch_7|^2,\\
N_3&=&\ch_0\ch_3^\ast+\ch_3\ch_0^\ast+\ch_1\ch_3^\ast+
\ch_3\ch_1^\ast +|\ch_3|^2+\ch_2(\ch_4^\ast+\ch_5^\ast)+
(\ch_4+\ch_5)\ch_2^\ast+\ch_4\ch_5^\ast+\\ &&
\ch_5\ch_4^\ast+|\ch_6+\ch_7|^2,\\
N_4&=&\ch_0\ch_4^\ast+\ch_4\ch_0^\ast +\ch_1\ch_4^\ast
+\ch_4\ch_1^\ast+|\ch_4|^2+
 \ch_5(\ch_2+\ch_3)^\ast+ (\ch_2+\ch_3)\ch_5^\ast +\ch_3\ch_2^\ast+\\
 &&\ch_2\ch_3^\ast+|\ch_6+\ch_7|^2 ,\\
N_5&=&\ch_0\ch_5^\ast+\ch_5\ch_0^\ast+\ch_1\ch_5^\ast+\ch_5\ch_1^\ast+
|\ch_5|^2+ \ch_4(\ch_2^\ast+\ch_3^\ast)+ (\ch_2+\ch_3)\ch_4^\ast
+\ch_2\ch_3^\ast+\\
&&\ch_3\ch_2^\ast+|\ch_5|^2+|\ch_6+\ch_7|^2,\\
N_6&=&\ch_0\ch_6^\ast+\ch_6\ch_0^\ast+\ch_1\ch_7^\ast+\ch_7\ch_1^\ast+
(\ch_2+\ch_3+\ch_4+\ch_5)(\ch_6^\ast+\ch_7^\ast)+\\
&& (\ch_6+\ch_7) (\ch_2+\ch_3+\ch_4+\ch_5)^\ast ,\\
N_7&=&\ch_0\ch_7^\ast+\ch_7\ch_0^\ast+\ch_1\ch_6^\ast+\ch_6\ch_1^\ast+
(\ch_2+\ch_3+\ch_4+\ch_5) (\ch_6^\ast+\ch_7^\ast)+ \\
&& (\ch_6+\ch_7)(\ch_2+\ch_3+\ch_4+\ch_5)^\ast.
\end{eqnarray*}

In particular all the primary fields are self-conjugate and
moreover the Frobenius-Schur indicator FS$_\la=1$ for all
$\la\in\NXN$.

\subsection{The modular invariants of the quantum $S_3$ double}

We present here the quantum $S_3$ double modular invariants.
There are precisely 48 modular invariants. The identity matrix
and the permutation matrix corresponding to the interchange
$2\leftrightarrow 3$ are the permutation invariants $Z_1$ and
$Z_2$ respectively. Then we get three new modular symmetric
invariants:
$$|\ch_0+\ch_1|^2+2|\ch_2|^2+2|\ch_3|^2+2|\ch_4|^2+2|\ch_5|^2+
k(\ch_2\ch_3^\ast+\ch_3\ch_2^\ast-|\ch_2|^2-|\ch_3|^2)$$
where $k=0,1,2$, producing the matrices
$Z_5, Z_6$ and $Z_7$ respectively.
We have two more symmetric modular invariants ($k=\pm 1$), whose
matrices we denote by $Z_3, Z_4$, respectively:
$$|\ch_0+\ch_1+\ch_2+\ch_3|^2+k(\ch_2\ch_3^\ast+
\ch_3\ch_2^\ast-|\ch_2|^2+|\ch_3|^2).$$ Set
$s_1=\ch_0+\ch_1+\ch_2+\ch_3$, $s_2=\ch_0+\ch_1+2\ch_2$,
$s_3=\ch_0+\ch_1+2\ch_3$, $s_4=\ch_0+\ch_2+\ch_6$ and $s_5=
\ch_0+\ch_3+\ch_6$. Then for every pair $(i,j)$, $i,
j=1,2,3,4,5$, $Z_{ij}:=s_is_j^\ast$ is a modular invariant
partition function. These 32 modular invariants were
found in \cite{Gan} but the list is not complete.
The extra sixteen are listed below:
\begin{eqnarray*}
Z_{(1)}&=&|\ch_0+\ch_1+\ch_3|^2+\ch_2\ch_3^\ast+\ch_3\ch_2^\ast+
|\ch_3|^2+|\ch_4|^2+|\ch_5|^2,\\
Z_{(2)}&=&|\ch_0+\ch_1+\ch_3|^2+|\ch_2|^2+2|\ch_3|^2+|\ch_4|^2+|\ch_5|^2,\\
Z_{(3)}&=&|\ch_0+\ch_1+\ch_2|^2+|\ch_2|^2+\ch_2\ch_3^\ast+\ch_3\ch_2^\ast+
|\ch_4|^2+|\ch_5|^2,\\
Z_{(4)}&=&|\ch_0+\ch_1+\ch_2|^2+2|\ch_2|^2+|\ch_3|^2+|\ch_4|^2+|\ch_5|^2,\\
Z_{(22)}&=&|\ch_0+\ch_2|^2+|\ch_1+\ch_2|^2+|\ch_6|^2+|\ch_7|^2,\\
Z_{(33)}&=&|\ch_0+\ch_3|^2+|\ch_1+\ch_3|^2+|\ch_6|^2+|\ch_7|^2,\\
Z_{(44)}&=&|\ch_0+\ch_2|^2+\ch_1\ch_6^\ast+\ch_2\ch_6^\ast+
\ch_6\ch_1^\ast+\ch_6\ch_2^\ast+|\ch_7|^2,\\
Z_{(55)}&=&|\ch_0+\ch_3|^2+\ch_1\ch_6^\ast+\ch_3\ch_6^\ast+
\ch_6\ch_1^\ast+\ch_6\ch_3^\ast+|\ch_7|^2,\\
Z_{(32)}&=&(\ch_0+\ch_3)(\ch_0+\ch_2)^\ast+(\ch_1+\ch_3)(\ch_1+\ch_2)^\ast+
|\ch_6|^2+|\ch_7|^2\\
Z_{(23)}&=&Z_{(32)}^t,\\
Z_{(45)}&=&(\ch_0+\ch_3)(\ch_0+\ch_2)^\ast+\ch_1\ch_6^\ast+\ch_3\ch_6^\ast+
\ch_6\ch_1^\ast+\ch_6\ch_2^\ast+|\ch_7|^2,\\
Z_{(54)}&=&Z_{(45)}^t,\\
Z_{(21)}&=&(\ch_0+\ch_1+\ch_3)(\ch_0+\ch_1+\ch_2)^\ast+2\ch_3\ch_2^\ast+
\ch_2\ch_3^\ast+|\ch_4|^2+|\ch_5|^2,\\
Z_{(12)}&=&Z_{(21)}^t,\\
Z_{(61)}&=&(\ch_0+\ch_1+\ch_3)(\ch_0+\ch_1+\ch_2)^\ast+
\ch_3\ch_2^\ast+\ch_2\ch_2^\ast+|\ch_3|^2+|\ch_4|^2+
|\ch_5|^2\\
Z_{(16)}&=&Z_{(61)}^t.
\end{eqnarray*}

After a basis change as in \cite[Page 698]{Gan}, the $S$ and $T$ 
matrices become permutation matrices, corresponding to the permutations
$(06)(23)(45)$ and $(345)(67)$ of $S_8$. Then one easily sees that
the dimension of the commutant $\{S,T\}^\prime$ is 11.

\subsection{Nimless modular invariants}

In the sequel we find nimless modular invariants as an application
of the nimrep theory developed in \cite{BEK1} and \cite{Gan2}. We
use the graphs of norm 2 whose number of vertices is less or equal
to six (see \cite{GHJ} or \cite[Page 58]{Gan2}) together with their
Perron-Frobenius eigenvalues as displayed in \fig{graphs}
(which can be derived from Table \ref{graphs-norm2}).
The notation $^0A_n^0, D_n^0$ refers to tadpole graphs as 
illustrated in \fig{graphs}.
In any possible nimrep, the matrices $G_{[\la]}$ are symmetric
because the primary fields are self-conjugate.

\unitlength 1.4mm
\thinlines
\begin{table}[htb]
\begin{center}
{\footnotesize
\begin{tabular}{|c|c|c|}  \hline\hline
Graph & Eigenvalues & Range \\ \hline\hline
$A_n^{(1)}$, $n\ge 1$  & $2\cos(2\pi k/(n+1))$& $0\le k\le n$ \\
$D_n^{(1)}$, $n\ge 4$ & $0,0,2\cos(\pi k/(n-2))$ &
$0\le k\le n-2$\\
$E_6^{(1)}$ & $\pm 2,\pm 1,\pm 1,0$ & {}\\
$E_7^{(1)}$ & $\pm2,\pm\sqrt{2},\pm 1,0,0$ &{} \\
$E_8^{(1)}$ & $\pm 2,\pm 2\cos(\pi/5),\pm 1,\pm 2\cos(2\pi/5),0$&{}\\
${}^0\!A_n^0$, $n\ge 1$ & $2\cos(k\pi/n)$ &  $0\le k<n$  \\
$D_n^0$, $n\ge 3$ & $0,2\cos(2\pi k/(2n-3))$ &
$0\le k\le n-2 $ \\
\hline
\end{tabular}
}
\end{center}
\caption{Connected graphs with norm two and their extended ADE graphs}
\label{graphs-norm2}
\end{table}

\begin{proposition}
The modular invariants $Z_6$, $Z_{11}$, 
$Z_{14}$, $Z_{15}$ and the transposes $Z_{41}$, $Z_{51}$
are nimless.
\lablth{2dnimless}
\end{proposition}
\bproof Let us assume that there is a nimrep for the modular
invariant $Z_6$, then its exponents are  $\{0,1,2,3,4^2,5^2\}$. So
the spectrum of the matrix $G_{[5]}$ is $\{-1^4,2^4\}$. In the
sequel if $a$ is primary field or eigenvalue with multiplicity
$n$, we abbreviate this by $a^n$. Hence $G_{[5]}$ has four
indecomposable components each of them with norm two. By
\fig{graphs} (or Table \ref{graphs-norm2}) we conclude that there
is no $(n+1)$ by $(n+1)$ matrix with spectrum is $\{-1^n,2\}$ for
$n=1,2,3,4,5$. Therefore $Z_6$ is nimless.

If $Z_{11}$ is nimble then as its exponents are the primary
fields $\{1,2,3\}$, we find that the spectrum of $G_{[2]}$ is
$\{-1,2^3\}$. Hence, $G_{[2]}$ has 3 irreducible components, one
of them with spectrum being $\{-1,2\}$. There is no 2 by 2 matrix
whose spectrum is given by this set as seen in \fig{graphs}. Let
us suppose that there is a nimrep for $Z_{15}$, so Exp=$\{0,3\}$.
Then we are looking for a 2 by 2 symmetric non-negative matrix
$G_{[3]}$  with eigenvalues
$S_{20}/S_{00}=2 \hbox{ and } S_{23}/S_{03}=-1$. Again such a
matrix cannot exist. The same can be carried out (it is actually
in \cite[Page 38]{Gan2}) for $Z_{14}$ whose exponents are
$\{0,2\}$. \eproof

\begin{figure}[htb]
\begin{center}
\unitlength 0.7mm
\begin{picture}(120,140)
\thinlines

\put(-5,0){\makebox(0,0){$\bullet$}}
\put(5,0){\makebox(0,0){$\bullet$}} \path(-5,-0.5)(5,-0.5)
\path(-5,0.5)(5,0.5) {\scriptsize
\put(0,-5){\makebox(0,0){$ A_1^{(1)}\quad\{-2,2\}$}} }


\put(45,0){\makebox(0,0){$\bullet$}}
\put(55,0){\makebox(0,0){$\bullet$}} \put(43,0){\arc{4}{0}{6.300}}
\put(57,0){\arc{4}{0}{6.300}} \path(45,0)(55,0) {\scriptsize
\put(50,-5){\makebox(0,0){$ {}^0A_2^0\quad\{0,2\}$}} }

\put(92,0){\makebox(0,0){$\bullet$}}
\put(94,0){\arc{4}{0}{6.300}}
\put(95,0){\arc{6}{0}{6.300}}{\scriptsize
\put(95,-5){\makebox(0,0){${}^0A_1^0\quad\{2\}$}} }


\put(-8,20){\makebox(0,0){$\bullet$}}
\put(2,20){\makebox(0,0){$\bullet$}}
\put(12,20){\makebox(0,0){$\bullet$}} \path(-8,20)(5,20)(12,20)
\put(2,22){\arc{4}{0}{6.300}} {\scriptsize
\put(2,15){\makebox(0,0){$ D_3^0\quad\{-1,0,2\}$}} }

\put(40,20){\makebox(0,0){$\bullet$}}
\put(50,20){\makebox(0,0){$\bullet$}}
\put(60,20){\makebox(0,0){$\bullet$}} \path(40,20)(40,20)(60,20)
\put(38,20){\arc{4}{0}{6.300}} \put(62,20){\arc{4}{0}{6.300}}
{\scriptsize \put(50,15){\makebox(0,0)
{$ {}^0A_3^0\quad\{-1,1,2\}$}} }

\put(95,20){\makebox(0,0){$\bullet$}}
\put(100,30){\makebox(0,0){$\bullet$}}
\put(105,20){\makebox(0,0){$\bullet$}}
\path(95,20)(105,20)(100,30)(95,20) {\scriptsize
\put(100,16){\makebox(0,0){$A_2^{(1)}\quad\{-1,-1,2\}$}} }

\put(-20,50){\makebox(0,0){$\bullet$}}
\put(-10,50){\makebox(0,0){$\bullet$}}
\put(-20,60){\makebox(0,0){$\bullet$}}
\put(-10,60){\makebox(0,0){$\bullet$}}
\path(-20,50)(-10,50)(-10,60)(-20,60)(-20,50)  {\scriptsize
\put(-15,45){\makebox(0,0){$ A_3^{(1)}\quad\{-2,0,0,2\}$}} }

\put(30,50){\makebox(0,0){$\bullet$}}
\put(40,50){\makebox(0,0){$\bullet$}}
\put(50,50){\makebox(0,0){$\bullet$}}

\put(47,58){\makebox(0,0){$\bullet$}}
\path(30,50)(40,50)(50,50)
\path(40,50)(47,58)
\put(28,50){\arc{4}{0}{6.300}} {\scriptsize
\put(40,45){\makebox(0,0){$ D_4^0\quad\{0,(-1\pm\sqrt{5})/2,2\}$}}}

\put(80,50){\makebox(0,0){$\bullet$}}
\put(90,50){\makebox(0,0){$\bullet$}}
\put(100,50){\makebox(0,0){$\bullet$}}
\put(110,50){\makebox(0,0){$\bullet$}} \path(80,50)(110,50)
\put(78,50){\arc{4}{0}{6.300}} \put(112,50){\arc{4}{0}{6.300}}
{\scriptsize \put(95,45){\makebox(0,0)
{${}^0A_4^0\quad\{0,\pm\sqrt{2},2\}$}}}

\put(130,50){\makebox(0,0){$\bullet$}}
\put(130,60){\makebox(0,0){$\bullet$}}
\put(135,55){\makebox(0,0){$\bullet$}}
\put(140,60){\makebox(0,0){$\bullet$}}
\put(140,50){\makebox(0,0){$\bullet$}}
\path(130,50)(140,60)
\path(130,60)(140,50)
{\scriptsize
\put(138,45){\makebox(0,0){$ D_4^{(1)}\quad\{-2,0^3,2\}$}} }


\put(-30,85){\makebox(0,0){$\bullet$}}
\put(-20,85){\makebox(0,0){$\bullet$}}
\put(-10,85){\makebox(0,0){$\bullet$}}
\put(-15,95){\makebox(0,0){$\bullet$}}
\put(-25,95){\makebox(0,0){$\bullet$}}
\path(-30,85)(-20,85)(-10,85)(-15,95)(-25,95)(-30,85)

{\scriptsize \put(10,79){\makebox(0,0)
{$ A_4^{(1)}\quad\{(-1\pm\sqrt{5})/2,
(-1\pm\sqrt{5})/2,2 \}\quad {}^0A_5^0$}} }


\put(30,85){\makebox(0,0){$\bullet$}}
\put(35,90){\makebox(0,0){$\bullet$}}
\put(40,85){\makebox(0,0){$\bullet$}}
\put(45,90){\makebox(0,0){$\bullet$}}
\put(50,85){\makebox(0,0){$\bullet$}}
\path(30,85)(35,90)(40,85)(45,90)(50,85)
\put(28,85){\arc{4}{0}{6.300}}
\put(52,85){\arc{4}{0}{6.300}}

\put(100,85){\makebox(0,0){$\bullet$}}
\put(110,85){\makebox(0,0){$\bullet$}}
\put(120,85){\makebox(0,0){$\bullet$}}
\put(127,91){\makebox(0,0){$\bullet$}}
\put(130,85){\makebox(0,0){$\bullet$}}
\path(100,85)(110,85)(120,85)(127,91) \path(120,85)(130,85)
\put(98,85){\arc{4}{0}{6.300}} {\scriptsize
\put(115,79){\makebox(0,0){$ D_5^0\quad\{0,2,\simeq 1.2,
\simeq -1.8, \simeq -0.44\}$}} }


\put(-30,120){\makebox(0,0){$\bullet$}}
\put(-20,120){\makebox(0,0){$\bullet$}}
\put(-10,120){\makebox(0,0){$\bullet$}}
\put(-30,130){\makebox(0,0){$\bullet$}}
\put(-20,130){\makebox(0,0){$\bullet$}}
\put(-10,130){\makebox(0,0){$\bullet$}}
\path(-30,120)(-20,120)(-10,120)(-10,130)
(-20,130)(-30,130)(-30,120) {\scriptsize
\put(-23,115){\makebox(0,0){$ A_5^{(1)}\quad\{\pm1,\pm1,\pm2\}$}} }

\put(10,120){\makebox(0,0){$\bullet$}}
\put(15,125){\makebox(0,0){$\bullet$}}
\put(20,120){\makebox(0,0){$\bullet$}}
\put(25,125){\makebox(0,0){$\bullet$}}
\put(30,120){\makebox(0,0){$\bullet$}}
\put(35,125){\makebox(0,0){$\bullet$}}
\path(10,120)(15,125)(20,120)(25,125)(30,120)(35,125)
\put(8,120){\arc{4}{0}{6.300}} \put(37,125){\arc{4}{0}{6.300}}
{\scriptsize
\put(20,115){\makebox(0,0){$ {}^0A_6^0\quad\{\pm1,\pm\sqrt{3},0,2\}$}} }

\put(60,120){\makebox(0,0){$\bullet$}}
\put(70,120){\makebox(0,0){$\bullet$}}
\put(80,120){\makebox(0,0){$\bullet$}}
\put(90,120){\makebox(0,0){$\bullet$}}
\put(95,127){\makebox(0,0){$\bullet$}}
\put(98,120){\makebox(0,0){$\bullet$}}
\path(60,120)(70,120)(80,120)(90,120)(98,120)
\path(90,120)(95,127) \put(58,120){\arc{4}{0}{6.300}} {\scriptsize
\put(80,115){\makebox(0,0){$ D_6^0\ \{-1,0,2,\simeq 1.53, \simeq
0.34,\simeq -1.87\}$}} }

\put(129,120){\makebox(0,0){$\bullet$}}
\put(129,130){\makebox(0,0){$\bullet$}}
\put(134,125){\makebox(0,0){$\bullet$}}
\put(143,125){\makebox(0,0){$\bullet$}}
\put(150,130){\makebox(0,0){$\bullet$}}
\put(150,120){\makebox(0,0){$\bullet$}}
\path(129,120)(134,125)(143,125)(150,120) \path(129,130)(134,125)
\path(143,125)(150,130) {\scriptsize
\put(140,115){\makebox(0,0){$ D_5^{(1)}\quad\{0,0,\pm1,\pm2\}$}} }

\end{picture}
\end{center}
\caption{Connected graphs with norm two and $\#$vertices $\leq 6$
together with their eigenvalues (cf.\ Table \ref{graphs-norm2})}
\label{graphs}
\end{figure}
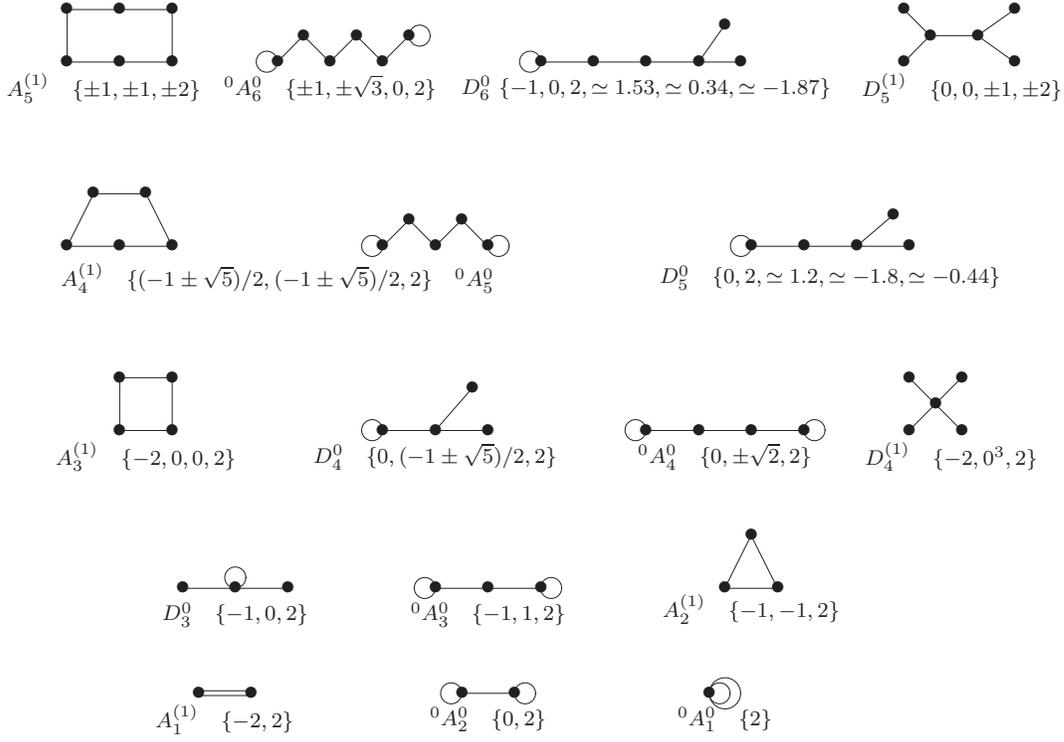

\begin{proposition}
The modular invariants $Z_{(21)}$, 
$Z_{(61)}$, $Z_{(3)}$, $Z_{(4)}$, $Z_{(1)}$, $Z_{(12)}$,
$Z_{(2)}$ and $Z_{(16)}$
are nimless.
\lablth{3dnimless}
\end{proposition}

\bproof The exponents of the trace 6 modular invariants
$Z_{(61)}$, $Z_{(16)}=Z_{(61)}^t$ are $\{0,1,2,3,4,5\}$. If a
nimrep exists then the spectrum of $G_{[2]}$ is $\{-1^3,2^3\}$ so
$G_{[2]}$ has three connected components. However there is no norm
two matrix with spectrum $\{-1,2\}$ or $\{-1^3,2\}$ (see
\fig{graphs}). Therefore $Z_{(61)}$ and $Z_{(16)}$ are nimless.
The exponents of the trace 6 matrix $Z_{(3)}$ are
$\{0,1,2^2,4,5\}$. The spectrum of $G_{[7]}$ in a nimrep (if it
exists) is also $\{-1^3,2^3\}$, therefore $Z_{(3)}$ is nimless by
the above argument. The exponents of the trace 8 matrix $Z_{(4)}$
are $\{0,1,2^3,3,4,5\}$. If a nimrep exists, then the spectrum of
$G_{[2]}$ is $\{-1^3,2^5\}$. Hence $G_{[2]}$ has five connected
components. By the above arguments such a matrix $G_{[2]}$ cannot
exist. So $Z_{(4)}$ is nimless. The interchange of the primary
fields 2 and 3 gives us the exponents of $Z_{(2)}$. Therefore
$Z_{(2)}$ is also nimless. The exponents of $Z_{(21)}$ and
$Z_{(12)}$ are $\{0,1,2^2,6,7\}$. If they are nimble, then the
spectrum of $G_{[5]}$ is $\{-1,2^3\}$. There is however no norm 2
graph whose spectrum is $\{-1,2\}$. Hence these modular invariants
are nimless. Finally, from $Z_{(3)}$ we obtain the exponents of
$Z_{(1)}$ by the interchange of the primary fields 2 and 3. So
$Z_{(1)}$ is nimless. \eproof

\subsection{Nimble modular invariants}

We produce here by hand the list of quantum $S_3$ double modular
invariants that have a compatible nimrep.

\begin{proposition}
The trace one modular invariants $Z_{25}, Z_{52}$, $Z_{34}$ and $Z_{43}$ have
unique nimreps.
\end{proposition}
\bproof The primary field 0 is the only exponent. We have
trivially the (unique) nimrep: $G_{[0]}= G_{[1]}=1,\ $ $G_{[2]}=
G_{[3]}= G_{[4]}= G_{[5]}=2,\ $ $G_{[6]}= G_{[7]}=3.$ \eproof

\begin{proposition} The trace two modular invariants $Z_{3}, Z_{23}, Z_{32}$,
$Z_{45}$, $Z_{54}$, $Z_{(45)}$ and $Z_{(54)}$ have unique nimreps.
\end{proposition}
\bproof The exponents of $Z_3$, $Z_{23}$ and $Z_{32}$ are
$\{0,1\}$, and we obtain (unique) nimrep:
$$G_{[0]}=G_{[1]}=\mathrm{id},\quad
G_{[2]}=G_{[3]}=G_{[4]}=G_{[5]}=2\mathrm{id},\quad
G_{[6]}=G_{[7]}=\left(\matrix {0&3\cr 3&0}\right).$$ We likewise
construct the unique nimrep for both $Z_{45}$ and $Z_{54}$. Their
exponents are $\{0,6\}$. The spectrum of every $G_{[\la]}$ is
provided in Table \ref{Z45Z(2)}. {\footnotesize
 \unitlength 1.4mm
\thinlines
\begin{table}[hhh]
{\footnotesize $$
\begin{array}{|c|c c|} \hline\hline
\emph{$Z_{45}, Z_{54}$}& 0 &6\\ \hline\hline
G_{[0]}&1&1\\
G_{[1]}&1&-1\\
G_{[2]}&2&0\\
G_{[3]}&2&0\\
G_{[4]}&2&0\\
G_{[5]}&2&0\\
G_{[6]}&3&1\\
G_{[7]}&3&-1\\ \hline
\end{array}
\qquad\qquad
\begin{array}{|c|c c|}  \hline\hline
\emph{$Z_{(45)},Z_{(54)}$}& 0& 7\\ \hline\hline
G_{[0]}&1&1\\
G_{[1]}&1&-1\\
G_{[2]}&2&0\\
G_{[3]}&2&0\\
G_{[4]}&2&0\\
G_{[5]}&2&0\\
G_{[6]}&3&-1\\
G_{[7]}&3&1\\ \hline
\end{array}
$$
}
\caption{Eigenvalues of $G_{[\la]}$ for $Z_{45}$, $Z_{54}$ and
$Z_{(45)}$, $Z_{(54)}$} \label{Z45Z(2)}
\end{table}
}
 We therefore find that $G_{[0]}=\mathrm{id},$
$$G_{[1]}={\scriptscriptstyle\left(\matrix{0&1\cr 1&0}\right)},\quad
G_{[2]}=G_{[3]}=G_{[4]}=
G_{[5]}={\scriptscriptstyle\left(\matrix{1&1\cr 1&1}\right)}, 
G_{[6]}={\scriptscriptstyle\left(\matrix{1&2\cr 2&1}\right)},
\hbox{ and }\quad G_{[7]}=G_{[1]}G_{[6]}.$$

Finally, the exponents of the trace 2 matrices $Z_{(45)}$ and
$Z_{(54)}$ are $\{0,7\}$. The list of the eigenvalues of every
matrix $G_{\la}$ is as in Table \ref{Z45Z(2)}. We find the
following (unique) nimrep for these modular invariants:
$$G_{[0]}=\mathrm{id},\quad G_{[1]}={\scriptscriptstyle\pmatrix{0&1\cr 1&0}},
\quad
G_{[2]}=
G_{[3]}= G_{[4]}=
G_{[5]}={\scriptscriptstyle\pmatrix{1&1\cr1&1}}$$
$$G_{[6]}={\scriptscriptstyle\pmatrix{2&1\cr 1&2}},\quad G_{[7]}=
{\scriptstyle\pmatrix{1&2\cr2&1}}.$$

\eproof

\begin{proposition}
The trace three modular invariants $Z_{24}$, $ Z_{42}$, $Z_{35}$,
$ Z_{53}$, $Z_{44}$, $Z_{(55)}$ and $Z_{(44)}$ have unique
nimreps.
\end{proposition}
\bproof The exponents of $Z_{44}$ are $\{0, 2, 6\}$. So the
spectrum of each non-negative matrix $G_{[\la]}$ is provided in
Table \ref{Z44Z35}. Using \fig{graphs}, we can fix (i.e.\ up to
permutations) $G_{[2]}$ as
$$G_{[2]}=\left(\matrix{2& 0&0\cr 0&1&1\cr 0&1&1}\right).$$
Then by the fusion rules,
$G_{[2]}\cdot G_{[3]}=G_{[4]}+ G_{[5]}=2 G_{[3]}.$
This equation holds true only when $G_{[3]}=\left(\matrix{1& 1&1\cr
1&0&0\cr 1&0&0}\right)$ using \fig{graphs}. Moreover, we obtain
$$\left(\matrix{4& 0&0\cr 0&2&2\cr 0&2&2}\right)=
G_{[2]}\cdot G_{[2]}= G_{[0]}+ G_{[1]}+ G_{[2]}=G_{[1]}+
\left(\matrix{3& 0&0\cr 0&2&1\cr 0&1&2}\right)$$
which implies that $G_{[1]}=\left(\matrix{1&0&0\cr 0&0&1\cr 0&1&0}\right)$.

Now we look for a 3 by 3 symmetric matrix $G_{[7]}$ with entries
from $\{0,1,2,3\}$ whose trace equals $3-1=2$. From the fusion
rules
$$G_{[7]}\cdot
G_{[7]}=G_{[0]}+G_{[2]}+G_{[3]}+G_{[4]}+G_{[5]}= G_{[0]}+
G_{[2]}+3G_{[3]}={\scriptstyle\left(\matrix{6&3&3\cr 3&2&1\cr
3&1&2}\right)}$$ so in particular the entries of $G_{[7]}$ are
from $\{0,1,2\}$. We can now easily check that the trace 2 matrix
$G_{[7]}$ has to be $G_{[7]}={\scriptstyle\left(\matrix{2&1&1\cr
1&0&1\cr 1&1&0}\right)}$. Finally, let us look for a 3 by 3 matrix
$$G_{[6]}={\scriptstyle\left(\matrix{a&b&c\cr b&y&f\cr c&f&d}\right)}$$
whose entries belong to $\{0,1,2,3\}$. In this way, by the fusion
rules,
$${\scriptstyle\left(\matrix{2a&b+c&b+c\cr 2b&y+f&y+f\cr
2c&d+f&d+f}\right)}= G_{[6]}\cdot G_{[2]}=G_{[6]}+G_{[7]}=
{\scriptstyle\left(\matrix{2+a&b+1&c+1\cr b+1&y&f+1\cr
c+1&f+1&d}\right)}.$$ But we also know that the fusion rules are
commutative, so
$${\scriptstyle\left(\matrix{2a&b+c&b+c\cr 2b&y+f&y+f\cr
2c&d+f&d+f}\right)=
\left(\matrix{2a&2b&2c\cr b+c&y+f&f+d\cr b+c&y+f&f+d}\right) }$$
leading to the conclusion: $a=2, b=c=d=y=1, f=0$.

\unitlength 1.4mm
\thinlines
\begin{table}[hhh]
{\footnotesize
$$
\begin{array}{|c|c c c|}  \hline\hline
\emph{$Z_{44}$}& 0& 2& 6\\ \hline\hline
G_{[0]}&1&1&1\\
G_{[1]}&1&1&-1\\
G_{[2]}&2&2&0\\
G_{[3]}&2&-1&0\\
G_{[4]}&2&-1&0\\
G_{[5]}&2&-1&0\\
G_{[6]}&3&0&1\\
G_{[7]}&3&0&-1\\
\hline
\end{array}
\qquad\qquad
\begin{array}{|c|c c c|}  \hline\hline
\emph{$Z_{35}, Z_{53}$}&0&3&3 \\
\hline
G_{[0]}&1&1&1\\
G_{[1]}&1&1&1\\
G_{[2]}&2&-1&-1\\
G_{[3]}&2&2&2\\
G_{[4]}&2&-1&-1\\
G_{[5]}&2&-1&-1\\
G_{[6]}&3&0&0\\
G_{[7]}&3&0&0\\
\hline
\end{array}
$$
}
\caption{Eigenvalues of $G_{[\la]}$ for $Z_{44}$, $Z_{35}$ and
$Z_{53}$} \label{Z44Z35}
\end{table}

Next we consider $Z_{35}$ (and $Z_{53}$), with its exponents
$\{0,3^2\}$. The spectra of the matrices $G_{[\la]}$ are in Table
\ref{Z44Z35}. So using \fig{graphs}:
$G_{[0]}=G_{[1]}=\mathrm{id},\quad G_{[3]}= 2\mathrm{id},$
$$G_{[2]}=G_{[4]}=G_{[5]}=\left(\matrix {0&1&1\cr 1&0&1\cr 1&1&0}\right).$$
Now we look for fusion matrices
$G_{[6]}=G_{[7]}={\scriptstyle\left(\matrix {a&b&c\cr b&d&e\cr
c&e&f}\right)}$ with $a,b,c,d,e,f=0,1,2,3$. By the fusion rules,
$$G_{[6]}\cdot G_{[6]}=G_{[0]}+G_{[2]}+G_{[3]}+G_{[4]}+G_{[5]}
=G_{[0]}+3G_{[2]}+ G_{[3]} ={\scriptstyle\left(\matrix {3&3&3\cr
3&3&3\cr 3&3&3}\right)}.$$ Now $G_{[6]}$ is symmetric with
positive spectrum and its square is the matrix with 3 everywhere.
So by spectral theory, $G_{[6]}$ is unique and has to be the
matrix with 1 everywhere. Now it is easy to check that this
$G_{[6]}$ is consistent with the other fusion rules. Therefore
the modular invariants $Z_{35}$ and $Z_{53}$ are nimble. The
exponents of $Z_{24}$ and $Z_{42}$ are $\{0,2^2\}$. The spectra
of the 3 by 3 matrices $G_{[\la]}$ are in Table \ref{Z(10)Z24}.
We then compute as in the $Z_{44}$ case and get the following
nimrep: $G_{[0]}=G_{[1]}=\mathrm{id},\quad G_{[2]}=2\mathrm{id},$
$$G_{[3]}=G_{[4]}= 
G_{[5]}={\scriptstyle\left(\matrix{0&1&1\cr 1&0&1 \cr 1&1&0}\right)},
\quad G_{[5]}=G_{[6]}={\scriptstyle\left(\matrix{1&1&1\cr 1&1&1\cr
1&1&1}\right)}.$$

\unitlength 1.4mm
\thinlines
\begin{table}[hhh]
{\footnotesize
$$
\begin{array}{|c|c c c|}  \hline\hline
\emph{$Z_{(55)}$}& 0 & 3&{7} \\ \hline\hline
G_{[0]}&1&1&1\\
G_{[1]}&1&1&-1\\
G_{[2]}&2&-1&0\\
G_{[3]}&2&2&0\\
G_{[4]}&2&-1&0\\
G_{[5]}&2&-1&0\\
G_{[6]}&3&0&-1\\
G_{[7]}&3&0&1\\ \hline
\end{array}
\qquad
\begin{array}{|c|c c c|}  \hline\hline
\emph{$Z_{24}, Z_{42}$}&0&2&2 \\ \hline\hline
G_{[0]}&1&1&1\\
G_{[1]}&1&1&1\\
G_{[2]}&2&2&2\\
G_{[3]}&2&-1&-1\\
G_{[4]}&2&-1&-1\\
G_{[5]}&2&-1&-1\\
G_{[6]}&3&0&0\\
G_{[7]}&3&0&0\\ \hline
\end{array}
\qquad
\begin{array}{|c|c c c|}  \hline\hline
\emph{$Z_{(44)}$}&0 & 2 &7 \\ \hline\hline
G_{[0]}&1&1&1\cr
G_{[1]}&1&1&-1\cr
G_{[2]}&2&2&0\\
G_{[3]}&2&-1&0\\
G_{[4]}&2&-1&0\\
G_{[5]}&2&-1&0\\
G_{[6]}&3&0&-1\\
G_{[7]}&3&0&1\\ \hline
\end{array}
$$
}
\caption{Eigenvalues of $G_{[\la]}$ for $Z_{(44)}$, $Z_{(55)}$,
$Z_{24}$ and $Z_{42}$} \label{Z(10)Z24}
\end{table}

The exponents of the trace 3 matrix $Z_{(55)}$ are $\{0,3,7\}$.
The list of the eigenvalues of every matrix $G_{[\la]}$ is as in
Table \ref{Z(10)Z24}. Since this table is the same (up to a
permutation of the primary fields) as that for $Z_{44}$, we
conclude that $Z_{(55)}$ is nimble.
We have Exp=$\{0,2,7\}$ for the trace 3 matrix $Z_{(44)}$. The
eigenvalues of $G_{\la}$ are provided in Table \ref{Z(10)Z24}.
Thus, a permutation of matrices of $Z_{(55)}$ give rise to a
nimrep of $Z_{(44)}$.

\eproof
In the following we use $A\oplus B$ for the block matrix $\pmatrix{A&0\cr 0& B}$ and
$A\otimes\pmatrix{0&1\cr 1& 0}$ for block matrix $\pmatrix{0& A\cr A& 0}$.

\begin{proposition}
The trace four modular invariants $Z_{13}$, $Z_{31}$,
$Z_{12}$, $Z_{21}$, $Z_{(32)}$ and $Z_{(23)}$ have unique nimreps.
\end{proposition}

\bproof The primary fields $\{0,1,2^2\}$ are the exponents of
both $Z_{12}$ and $Z_{21}$. The eigenvalues are in Table
\ref{Z12Z(1)}. We have following nimrep for $Z_{21}$:
$$G_{[0]}=G_{[1]}=\mathrm{id},\quad G_{[2]}=2\mathrm{id},$$
$$G_{[3]}=G_{[4]}= G_{[5]}=
{\scriptstyle\left(\matrix{2&0&0&0\cr 0&0&1&1 \cr 0&1&0&1\cr 
0&1&1&0}\right)},\quad
G_{[6]}= G_{[7]}={\scriptstyle\left(\matrix{0&2&2&2\cr 2&0&0&0 \cr
2&0&0&0\cr 2&0&0&0}\right)}.$$ A nimrep for $Z_{13}$ is obtained
from that of $Z_{12}$ by interchanging the primary fields 2 and 3.

The exponents of the trace 4 matrices $Z_{(32)}$ and $Z_{(23)}$
are $\{0,1,6,7\}$. The spectra of the matrices $G_{[\la]}$ are in
Table \ref{Z12Z(1)}. We find the following nimrep for these
modular invariants:
$$G_{[0]}=\mathrm{id},\quad G_{[1]}=
{\scriptstyle \pmatrix{0&1\cr 1&0}\oplus\pmatrix{0&1\cr 1&0}},\quad
G_{[2]}= G_{[3]}= G_{[4]}=G_{[5]}={\scriptstyle \pmatrix{1&1\cr
1&1}\oplus\pmatrix{1&1\cr 1&1}},$$
$$ G_{[6]}={\scriptstyle\pmatrix{2&1\cr 1&2} \otimes\pmatrix{0&1\cr 1&0}},
\quad G_{[7]}={\scriptstyle
\pmatrix{1&2\cr 2&1}\otimes\pmatrix{0&1\cr 1&0}}.$$ \eproof

\unitlength 1.4mm
\thinlines
\begin{table}[hhh]
{\footnotesize
$$
\begin{array}{|c|c c c c|}  \hline\hline
\emph{$Z_{12}, Z_{21}$}&0 &1 & 2&2\\ \hline\hline
G_{[0]}&1&1&1&1\\
G_{[1]}&1&1&1&1\\
G_{[2]}&2&-2&2&2\\
G_{[3]}&2&2&-1&-1\\
G_{[4]}&2&2&-1&-1\\
G_{[5]}&2&2&-1&-1\\
G_{[6]}&3&-3&0&0\\
G_{[7]}&3&-3&0&0\\ \hline
\end{array}
\qquad\qquad
\begin{array}{|c|c c c c|}  \hline\hline
\emph{$Z_{(32)}, Z_{(23)}$}&0&1&6&7 \\ \hline\hline
G_{[0]}&1&1&1&1\cr
G_{[1]}&1&1&-1&-1\cr
G_{[2]}&2&2&0&0\\
G_{[3]}&2&2&0&0\\
G_{[4]}&2&2&0&0\\
G_{[5]}&2&2&0&0\\
G_{[6]}&3&-3&1&-1\\
G_{[7]}&3&-3&-1&1\\ \hline
\end{array}
$$
}
\caption{Eigenvalues of $G_{[\la]}$ for $Z_{12}$, $Z_{21}$,
$Z_{(32)}$ and $Z_{(23)}$} \label{Z12Z(1)}
\end{table}

\begin{proposition}
The trace six modular invariants $Z_{(22)}$, $Z_2$, $Z_4$, $Z_7$,
$Z_{22}$, $Z_{33}$ and $Z_{(33)}$ have unique nimreps.
\lablth{z22z33}
\end{proposition}

\bproof The exponents of the trace 6 matrix $Z_{(22)}$ are
$\{0,1,2^2,6,7\}$. The spectra of the matrices $G_{[\la]}$ are as
in Table \ref{Z(4)}. Then, $G_{[0]}=\mathrm{id}$ and we can fix
$$G_{[3]}=G_{[4]}=G_{[5]}={\scriptstyle\pmatrix{1&1&1\cr 1&0&0\cr
1&0&0}\oplus \pmatrix{1&1&1\cr 1&0&0\cr
1&0&0}}.$$

\unitlength 1.4mm
\thinlines
\begin{table}[hhh]
{\footnotesize
$$
\begin{array}{|c|c c c c c c|}  \hline\hline
\emph{$Z_2$}& 0 &{1}&{4}&{5}&{6}&{7} \\ \hline\hline
G_{[0]}&1&1&1&1&1&1\\
G_{[1]}&1&1&1&1&-1&-1\\
G_{[2]}&2&2&-1&-1&0&0\\
G_{[3]}&2&2&-1&-1&0&0\\
G_{[4]}&2&2&-1&2&0&0\\
G_{[5]}&2&2&2&-1&0&0\\
G_{[6]}&3&-3&0&0&1&-1\\
G_{[7]}&3&-3&0&0&-1&1\\ \hline
\end{array}
\qquad
\begin{array}{|c|c c c c c c|}  \hline\hline
\emph{$Z_{(22)}$}& 0 &{1}&{2}&{2}&{6}&{7} \\ \hline\hline
G_{[0]}&1&1&1&1&1&1\\
G_{[1]}&1&1&1&1&-1&-1\\
G_{[2]}&2&2&2&2&0&0\\
G_{[3]}&2&2&-1&-1&0&0\\
G_{[4]}&2&2&-1&-1&0&0\\
G_{[5]}&2&2&-1&-1&0&0\\
G_{[6]}&3&-3&0&0&1&-1\\
G_{[7]}&3&-3&0&0&-1&1\\ \hline
\end{array}
$$
}
\caption{Eigenvalues of $G_{[\la]}$ for $Z_2$ and $Z_{(22)}$}
\label{Z(4)}
\end{table}
By the fusion rules
$G_{[3]}^2- G_{[3]}-G_{[0]}=G_{[1]}$
we see that
$$G_{[1]}={\scriptstyle\pmatrix{1&0&0\cr 0&0&1\cr
0&1&0}\oplus \pmatrix{1&0&0\cr 0&0&1\cr
0&1&0}}.$$ Next
let us use the fusion rules $G_{[2]}G_{[1]}=G_{[2]}$ and
$G_{[2]}\cdot G_{[3]}=G_{[4]}+ G_{[5]}=2 G_{[3]}$. We then obtain
the (unique) matrix
$$G_{[2]}=\pmatrix{2&0&0\cr 0&1&1\cr 0&1&1}\oplus 
\pmatrix{2&0&0\cr 0&1&1\cr 0&1&1}.$$
Next we use the commutation of $G_{[6]}$ with
$G_{[1]}$, $G_{[2]}$ and $G_{[3]}$, together with the fusion rule
$G_{[6]}^2=G_{[0]}+G_{[2]}+3G_{[3]}$ to get the matrix
$$G_{[6]}={\scriptstyle\pmatrix{2&1&1\cr 1&0&1\cr
1&1&0\cr}\otimes\pmatrix{0&1\cr 1&0}}.$$
Finally, the fusion rule $G_{[6]}G_{[1]}=G_{[7]}$ alone
determines $G_{[7]}$, namely:
$$G_{[7]}={\scriptstyle\pmatrix{2&1&1\cr 1&1&0\cr 1&0&1}
\otimes\pmatrix{0&1\cr 1&0}}.$$
Therefore $Z_{(22)}$ is nimble. Since $Z_{(33)}$ is obtained from
$Z_{(22)}$ by permuting the primary fields 2 and 3, we conclude
that $Z_{(33)}$ is also nimble. In a similar manner we obtain the
following nimrep for the permutation invariant $Z_2$:
$G_{[0]}=\rm{id}$,
$$G_{[1]}={\scriptstyle\pmatrix{1&0&0\cr 0&0&1\cr
0&1&0\cr}\oplus\pmatrix{1&0&0\cr 0&0&1\cr
0&1&0\cr}},
G_{[3]}=G_{[2]}={\scriptstyle\pmatrix{1&1&1\cr
1&0&0\cr 1&0&0\cr}\oplus \pmatrix{1&1&1\cr
1&0&0\cr 1&0&0\cr}},$$
$$G_{[4]}={\scriptstyle\pmatrix{2&0&0\cr 0&1&1\cr 0&1&1}\oplus
\pmatrix{1&1&1\cr 1&0&0\cr 1&0&0}},
G_{[5]}={\scriptstyle\pmatrix{1&1&1\cr 1&0&0\cr
1&0&0\cr}\oplus \pmatrix{2&0&0\cr 0&1&1\cr 0&1&1}},$$
$$G_{[6]}={\scriptstyle\pmatrix{2&1&1\cr 1&1&0\cr 1&0&1}
\otimes\pmatrix{0&1\cr 1&0} },
G_{[7]}={\scriptstyle\pmatrix{2&1&1\cr 1&0&1\cr 1&1&0}
\otimes\pmatrix{0&1\cr 1&0}}.$$

The exponents of $Z_4$ are $\{0,1,2^2,3^2\}$. The spectra of the
fusion matrices $G_{[\la]}$ are in Table \ref{Z4Z33}. We can
compute and find that the following set of matrices form a nimrep
for $Z_4$, $G_{[0]}=G_{[1]}=\mathrm{id}$,

\unitlength 1.4mm
\thinlines
\begin{table}[hhh]
{\footnotesize
$$
\begin{array}{|c|c c c c c c|}  \hline\hline
\emph{$Z_{4}$}& 0&1& 2&2 &{3}&3 \\ \hline\hline
G_{[0]}&1&1&1&1&1&1\\
G_{[1]}&1&1&1&1&1&1\\
G_{[2]}&2&2&2&2&-1&-1\\
G_{[3]}&2&2&-1&-1&2&2\\
G_{[4]}&2&2&-1&-1&-1&-1\\
G_{[5]}&2&2&-1&-1&-1&-1\\
G_{[6]}&3&-3&0&0&0&0\\
G_{[7]}&3&-3&0&0&0&0\\
\hline
\end{array}
\quad
\begin{array}{|c|c c c c c c|}  \hline\hline
\emph{$Z_{33}$}&0 &1&3&3&3&3\\ \hline\hline
G_{[0]}&1&1&1&1&1&1\\
G_{[1]}&1&1&1&1&1&1\\
G_{[2]}&2&2&-1&-1&-1&-1\\
G_{[3]}&2&2&2&2&2&2\\
G_{[4]}&2&2&-1&-1&-1&-1\\
G_{[5]}&2&2&-1&-1&-1& -1\\
G_{[6]}&3&-3&0&0&0&0\\
G_{[7]}&3&-3&0&0&0&0\\ \hline
\end{array}
$$
}
\caption{Eigenvalues of $G_{[\la]}$ for $Z_4$ and $Z_{33}$}
\label{Z4Z33}
\end{table}
$$G_{[2]}={\scriptstyle\left(\matrix{0&1&1\cr 1&0&1\cr 1&1&0}\right)\oplus
\pmatrix{2&0&0\cr 0&2&0\cr
0&0&2}},\quad
G_{[3]}={\scriptstyle\left(\matrix{2&0&0\cr 0&2&0\cr
0&0&2}\right)\oplus \pmatrix{0&1&1\cr 1&0&1\cr
1&1&0}},$$
$$ G_{[4]}= G_{[5]}={\scriptstyle \left(\matrix{0&1&1\cr 1&0&1\cr
1&1&0\cr}\right)\oplus \left(\matrix{0&1&1\cr 1&0&1\cr
1&1&0\cr}\right)},\
G_{[6]}=G_{[7]}={\scriptstyle\left(\matrix{1&1&1\cr 1&1&1\cr
1&1&1}\right)\otimes\pmatrix{0&1\cr 1&0}}.$$ The modular
invariant $Z_7$ has exponents $\{0,1,4^2,5^2\}$. Also, the nimrep
for the modular invariant $Z_7$ is obtained from $Z_4$ by
exchanging the exponents 2 with 5, 3 with 4, 4 with 3 and 5 with
2.

The exponents of $Z_{33}$ are $\{0,1,3^4\}$. The eigenvalues are
displayed in Table \ref{Z4Z33}. We thus find that
$G_{[0]}=G_{[1]}=\mathrm{id},\quad G_{[3]}=2\mathrm{id},$
$$G_{[2]}=G_{[4]}= G_{[5]}=
{\scriptstyle\left(\matrix{0&1&1\cr 1&0&1\cr 1&1&0}\right)},\
G_{[6]}=G_{[7]}={\scriptstyle\left(\matrix{1&1&1\cr 1&1&1\cr
1&1&1}\right)\otimes\pmatrix{0&1\cr 1&0}}.$$ When we interchange
the primaries 2 and 3 we get a nimrep for $Z_{22}$. \eproof

\subsubsection{Nimble and insufferable modular invariants}

We fix the system $\NXN$ arising from the quantum $S_3$ double,
regarding the primary field $i$ for $i=0,\dots,7$ as a morphism
$\lambda_i\in\NXN$ (see \cite{BEK2} and \fig{prgraph}). Since a
modular invariant that can be produced by a subfactor is nimble
we immediately get that those from Proposition \ref{2dnimless}
and \ref{3dnimless} are insufferable because they are nimless.

We now find a $\theta\in \Sigma(\NXN)$ that cannot be a dual
canonical endomorphism of any braided subfactor.
\begin{lemma}
The sector $[\lambda_0]\oplus[\lambda_1]\oplus[\lambda_2]\oplus [\lambda_3]$
cannot be that of the dual canonical endomorphism of a subfactor
\lablth{nondualend}
\end{lemma}

\bproof Let us suppose that
$[\theta]=[\lambda_0]\oplus[\lambda_1]\oplus[\lambda_2]\oplus[\lambda_3]$
is a dual canonical sector of some type III
subfactor with inclusion map $\iota$. Using Verlinde fusion
matrices $\lan\theta\lambda_0,\lambda_0\ran=1$ so
$[\iota\la_0]=[\iota]$ is irreducible,
$\lan\theta\la_0,\la_2\ran=\lan\theta\la_0, \la_3\ran=1,$ thus
$[\iota\lambda_0]$ is in the irreducible decomposition of
$[\iota\lambda_2]$ and $[\iota\lambda_3]$. Moreover
$\lan\theta\lambda_2,\lambda_2\ran=\lan\theta\lambda_2,
\lambda_2\ran=3,$\ $\lan\theta\lambda_2,\lambda_3\ran=0$ which is
incompatible with the fact that $[\iota\lambda_0]$ appears in the
decomposition of both $[\iota\la_2]$ and $[\iota\la_3]$.
Therefore, $\theta$ cannot be a dual canonical endomorphism of
any subfactor. \eproof

\begin{proposition}
The nimble modular invariants $Z_{3}$, $Z_4$, $Z_{12}$, $Z_{13}$,
$Z_{21}$ and $Z_{31}$ cannot be realised by subfactors.
\end{proposition}
\bproof Assume any one of the listed invariants is realised by a
subfactor $N\subset M$. Then
 $s_1$ is visible in either the first row or first column.The type I parent theorem
\cite[Theorem 4.7]{BE4}, yields an intermediate subfactor
$N\subset P\subset M$ such that $N\subset P$ has dual canonical
sector described by  $s_1$, i.e.\
$[\lambda_0]\oplus[\lambda_1]\oplus[\lambda_2]\oplus [\lambda_3]$.
This is however impossible by Lemma \ref{nondualend}. \eproof

\subsection{Subfactors producing quantum $S_3$ double modular invariants}
\lablth{submodinv}

We propose in this Section to study the remaining modular
invariants of $S_3$. As above we fix the system $\NXN$ arising
from the group $S_3$ and look for subfactors $N\subset M$ of type
III with inclusion map $\iota$ whose dual canonical endomorphism
$\theta=\bar{\iota}\iota\in\Sigma(\NXN)$ such that the associated
modular invariant $Z_{N\subset
M}=\lan\alpha_\lambda^+,\alpha_\mu^-\ran$ produces an invariant
from the list of $S_3$. The fundamental inclusion (the quantum
double of the subfactor $M_0\subset M_0\rtimes S_3$) is the subfactor
$N:=M_0\rtimes \Delta(S_3)$ $\subset$ $ M_0\rtimes(S_3\times
S_3)=:M$ where $M_0$ is a type III factor \cite{alain2}. In Table
\ref{psuff} we provide the list of 28 quantum $S_3$ double
modular invariants that are nimble and possibly sufferable that
will be studied further here.
\begin{table}
{\footnotesize
$$\begin{array}{||c||c|c|c|c|c|c|} \hline
\Tr(Z) & \multicolumn{6}{|c|} {\Tr(ZZ^t)} \\ \hline\hline
 {}&8& 9& 12& 18& 20& 36\\ \hline\hline
1&{}&{}& {} &Z_{25},Z_{34}, Z_{43},Z_{52} & {}& {}\\  \hline 2&
{}& Z_{45}, Z_{54}, Z_{(45)}, Z_{(54)} & {}& {} & {} & Z_{23},
Z_{32} \\ \hline 3& {}& Z_{44}, Z_{55}, Z_{(44)}, Z_{(55)}&
{}&Z_{24},Z_{42},Z_{35},Z_{53}&{} &{}\\ \hline 4& {}& {}&
Z_{(32)}, Z_{(23)}& {} & {} &{}\\ \hline 6& Z_{2}&{}& Z_{(22)},
Z_{(33)}& {} & Z_{7} & Z_{22}, Z_{33}\\ \hline 8& Z_{1}&{}& {}&
{} & {} &{}\\ \hline 10& {}& {}& {}& {} & Z_{5} &{}\\ \hline
\end{array}$$
}
\caption{The $(\Tr(Z),\Tr(ZZ^t))$ of nimble $S_3$ modular invariants}
\label{psuff}
\end{table}

\begin{proposition} The above fundamental subfactor 
realises the symmetric modular
invariant $Z_{55}$.
The dual canonical endomorphism
$\theta$ is visible in the vacuum block of $Z_{55}$ (so that chiral
locality holds). \lablth{fundZ55}
\end{proposition}

\bproof The irreducible $N$-$N$ sectors of $N\subset M$ are
labelled by $\widehat{S_3}\sqcup
\widehat{\bbZ_3}\sqcup\widehat{\bbZ_2}$ whereas the $N$-$M$ ones
are labelled by $\widehat{S_3}$ and finally the $M$-$M$
irreducible sectors by $\widehat{S_3}\times\widehat{S_3}$
\cite{KY}. The $N$-$M$ star vertex (i.e.\ the trivial
representation of $S_3$) is connected to the $N$-$N$ star vertex,
i.e.\ $[\lambda_0]$, the trivial representation of $\bbZ_3$, i.e.\
$[\lambda_3]$, and trivial representation of $\bbZ_2$, i.e.\
$[\lambda_6]$. Hence the dual canonical sector of $N\subset M$ is
$[\theta]=[\lambda_0]\oplus[\lambda_3]\oplus[\lambda_6]$. Let us
compute the modular invariant $Z_{N\subset M}$ for this braided
subfactor. From Verlinde fusion matrices, $N_0$, $N_3$ and $N_6$,
we obtain $\lan\theta\lambda_0,\lambda_0\ran=\lan\theta\lambda_1,
\lambda_1\ran=$ $\lan\theta\la_2,\la_2\ran=
\lan\theta\la_4,\la_4\ran=$ $\lan\theta\la_5,\la_5\ran=1,$
$\lan\theta\la_3,\la_3\ran=$ $\lan\theta\la_6,\la_6\ran=$
$\lan\theta\la_7,\la_7\ran=2,$
$\lan\theta\la_0,\la_1\ran=\lan\theta\la_0,\la_2\ran=$
$\lan\theta\la_0,\la_4\ran=\lan\theta\la_0, \lambda_5\ran=0,$
$\lan\theta\la_0,\la_3\ran=\lan\theta\la_0, \lambda_6\ran=1$,
$\lan\theta\la_0,\la_7\ran=\lan\theta\la_1,\la_2\ran=$
$\lan\theta\la_1,\la_4\ran=0$,
$\lan\theta\la_1,\la_5\ran=\lan\theta\la_1,\la_6\ran=0,$
$\lan\theta\la_1,\la_3\ran=\lan\theta\la_1,\la_7\ran=1.$ So
$[\iota\la_0], [\iota\la_1], [\iota\la_2]$ are irreducible $N$-$M$
sectors and $[\iota\la_3]=[\iota\la_0]\oplus[\iota\la_1]$.
Computing further, $\lan\theta\la_2,\la_3\ran=0,$
$\lan\theta\la_2,\la_4\ran=\lan\theta\la_2,\la_5\ran=$
$\lan\theta\la_2,\la_6\ran=\lan\theta\la_2,\la_7\ran=1,$ hence,
$[\iota \la_4]=[\iota \la_5]=$ $[\iota\la_2], [\iota\la_6]=[\iota
\la_0]\oplus[\iota\la_2]$ and
$[\iota\la_7]=[\iota\la_1]\oplus[\iota\la_2]$. Therefore
$\hbox{Tr}(Z_{N\subset M})$=$3$ which implies that $Z_{N\subset
M}=Z_{(44)}, Z_{(55)}, Z_{44}, Z_{55},$ $Z_{24}, Z_{35},$ $Z_{42}$
or  $Z_{53}$. Next we use the inequality
$\lan\alpha_\lambda^\pm,\alpha_\mu^\pm\ran\leq\lan\theta\lambda,\mu\ran$
valid by \cite[\eerf{37}]{BEK2} to eliminate more matrices from
the above list: since $\lan\theta\lambda_0,\lambda_2\ran=0$ but
$[Z_{44}]_{0,2}=[Z_{24}]_{0,2}=1$, $Z_{N\subset M}$ cannot be
either $Z_{44}$ nor $Z_{24}$; also,
$\lan\theta\lambda_0,\lambda_1\ran=0$ and $[Z_{42}]_{0,1}=
[Z_{53}]_{0,1}=1$, so $Z_{N\subset M}$ cannot be $Z_{42}$ or
$Z_{53}$; finally, as $\lan\theta\la_1,\la_0\ran=0$ and
$[Z_{35}]_{1,0}=1$, $Z_{N\subset M}$ cannot be $Z_{35}$. Similarly
we use the 02-entry of $ Z_{(44)}$ to rule it out. Finally, using
the curious identity \cite[Proposition 3.3]{E1} we conclude that
indeed $Z_{N\subset M}=Z_{55}$. \eproof

\begin{proposition}
The modular invariants $Z_1$ and $Z_5$ are realised by subfactors.
\lablth{z1z5}
\end{proposition}
\bproof The identity modular invariant $Z_1$ can always
be obtained by the trivial subfactor $N=M$.
The
simple currents are $\la_0$ and $\la_1$, which form the
commutative group $\bbZ_2$. By Lemma \ref{newfrob},
$[\theta]=[\la_0]\oplus[\la_1]$ is dual sector since
$w_{\la_1}^2=1$. Now we use Verlinde fusion matrices $N_0$ and
$N_1$ and obtain: $\lan\theta\lambda_0,\lambda_0\ran=$
$\lan\theta\lambda_1,\lambda_1\ran=$
$\lan\theta\lambda_6,\lambda_6\ran=$
$\lan\theta\lambda_7,\lambda_7\ran$
$=\lan\theta\lambda_0,\lambda_1\ran=\lan\theta\lambda_6,\lambda_7\ran=1.$
Therefore, $[\iota \lambda_0]=[\iota \lambda_1]$ and $[\iota
\lambda_6]= [\iota \lambda_7]$ provide two irreducible $N$-$M$
sectors where $\iota$ denotes the inclusion map $N\hookrightarrow
M$. On the other hand, $\lan\theta\lambda_2,\lambda_2\ran=
\lan\theta\lambda_3,\lambda_3\ran=
\lan\theta\lambda_4,\lambda_4\ran=\lan\theta\lambda_5,\lambda_5\ran=2,$
so the sectors $[\iota \lambda_2], [\iota \lambda_3], [\iota
\lambda_4], [\iota \lambda_5]$ decompose into two irreducible
$N$-$M$ sectors (which at this stage might be equal). However
since $\lan\theta\la_0,\la_2\ran=$ $\lan\theta\la_0,\la_3\ran=$
$\lan\theta\la_0,\la_4\ran=$ $\lan\theta\la_0,\la_5\ran=$
$\lan\theta\la_6,\la_2\ran=$ $\lan\theta\la_6,\la_3\ran=$
$\lan\theta\la_6,\la_4\ran=$ $\lan\theta\la_6,\la_5\ran=$
$\lan\theta\la_2,\la_3\ran=$ $\lan\theta\la_2,\la_4\ran=$
$\lan\theta\la_2,\la_5\ran=$ $\lan\theta\la_3,\la_2\ran=$
$\lan\theta\la_3,\la_4\ran=$ $\lan\theta\la_3,\la_5\ran=0,$
$\lan\theta\la_4,\la_2\ran=$ $\lan\theta\la_4, \la_3\ran=$
$\lan\theta\la_4,\la_5\ran=0,$ $\lan\theta\la_5,\la_2\ran=$
$\lan\theta\la_5,\la_3\ran=$ $\lan\theta\la_5,\la_4\ran=0,$ we
get a further 2+2+2+2=8 new irreducible $N$-$M$ sectors, making a
total of 10. Hence $Z_{N\subset M}=Z_5$ since it is the only
trace 10 quantum $S_3$ double modular invariant. \eproof

\begin{proposition} The
modular invariant $Z_{(33)}$ is produced by the subfactor
$N=M_0\rtimes S_3\subset M_0\rtimes(S_3\cdot
(\bbZ_3\times\bbZ_3))$. \lablth{nicetheta}
\end{proposition}
\bproof We consider the following intermediate subfactor
$$N=M_0\rtimes\Delta(S_3)\subset_{\iota_1}
M_0\rtimes(S_3\cdot (\bbZ_3\times\bbZ_3))\subset_{\iota_2}
M_0\rtimes (S_3\times S_3)=M.$$ Then $[\bar{\iota}_1\iota_1]$ is
a subsector of $[\theta_{N\subset
M}]=[\overline{\iota_2\iota_1}\iota_2\iota_1]=[\bar{\iota}_1\bar{\iota}_2
\iota_2\iota_1]$ because $[\la_0]$ is a subsector of
$[\bar{\iota}_2\iota_2]$. Since d$_{\bar{\iota}_1\iota_1}=3$ we
have that the dual canonical sector of 
$ M_0\rtimes S_3\subset M_0\rtimes(S_3\cdot (\bbZ_3\times\bbZ_3))$ is
$[\theta]=[\la_0]\oplus[\la_3]$. So now we use again the Verlinde
fusion matrices $N_0$ and $N_3$ and compute the dimensions of the
intertwiner spaces in order to get the $N$-$M$ induced system:
$\lan\theta\la_0,\la_0\ran=\lan\theta\la_1,\la_1\ran=\lan\theta\la_2,\la_2\ran=
\lan\theta\la_4,\la_4\ran=\lan\theta\la_5,\la_5\ran =1,$
$\lan\theta\la_3,\la_3\ran=$ $\lan\theta\la_6,\la_6\ran=
\lan\theta\la_7,\la_7\ran=2,$ $\lan\theta\la_0,\la_1\ran=$
$\lan\theta\la_0,\la_2\ran=$ $\lan\theta\la_0,\la_4\ran=$
$\lan\theta\la_0,\la_5\ran=$ $\lan\theta\la_0,\la_6\ran=$
$\lan\theta\la_0,\la_7\ran=0$, $\lan\theta\la_0,\la_3\ran=$
$\lan\theta\la_1,\la_3\ran=1$, $\lan\theta\la_1,\la_2\ran=$
$\lan\theta\la_1,\la_4\ran=$ $\lan\theta\la_1,\la_5\ran=$
$\lan\theta\la_1,\la_6\ran=$ $\lan\theta\la_1,\la_7\ran=0$,
$\lan\theta\la_2,\la_3\ran=$ $\lan\theta\la_2,\la_6\ran=$
$\lan\theta\la_2,\la_7\ran=0$, $\lan\theta\la_2,\la_4\ran=$
$\lan\theta\la_2,\la_5\ran= \lan\theta\la_6,\la_7\ran=1.$ Then
$[\iota\la_0]$, $[\iota\la_1]$, $[\iota\la_2]$ are distinct
irreducible sectors and
$[\iota\la_3]=[\iota\la_0]\oplus[\iota\la_1]$,
$[\iota\la_5]=[\iota\la_4]=[\iota\la_2]$,
$[\iota\la_6]=[a]\oplus[b]$ and $[\iota\la_7]=[a]\oplus[c]$ with
$[a]$, $[b]$ and $[c]$ three new irreducible sectors. Therefore,
$\Tr(Z_{M_0\rtimes S_3\subset M_0\rtimes(S_3\cdot
(\bbZ_3\times\bbZ_3))})=6$, implying that $Z_{M_0\rtimes S_3\subset
M_0\rtimes(S_3\cdot (\bbZ_3\times\bbZ_3))}=Z_{(22)},
Z_{(33)},Z_2,Z_7,Z_{22} \hbox{ or } Z_{33}$. Recall that in
general $\lan\a_\la^\pm,\a_\mu^\pm\ran\leq\lan\theta\la,\mu\ran$
(see \cite[\eerf{37}]{BEK2}). So since for our concrete $\theta$,
$\lan\theta\la_0,\la_1\ran=0$ then $Z_{M_0\rtimes S_3\subset
M_0\rtimes(S_3\cdot (\bbZ_3\times\bbZ_3))}$ cannot be $Z_7$,
$Z_{22}$, $Z_{33}$ because the 01-entry of those matrices is 1.
The $Z_{(22)}$ matrix is ruled out using its 02-entry. Hence
$Z=Z_2, Z_{(33)}$. We will use the curious identity
\cite[Proposition 3.3]{E1} (see \sect{mixing}) to eliminate one
of them. The RHS of that identity using $Z_2$ is
$$\bigoplus_{\la,\mu}{(Z_2)}_{\la,\mu}[\la\mu]=6[\la_0]\oplus2[\la_1]\oplus
2[\la_2]\oplus 2[\la_3]\oplus 5[\la_4]\oplus 5[\la_5].$$
The LHS is (with $[\theta]=[\la_0]\oplus[\la_3]$):
$$\bigoplus_{x\in\MXN}[x\bar{x}]=[\la_0\theta\la_0]\oplus[\la_1\theta\la_1]
\oplus[\la_2\theta\la_2]\oplus[a\bar{a}]\oplus[b\bar{b}]\oplus[c\bar{c}]=$$
$$3[\la_0]\oplus[\la_1]\oplus[\la_2]\oplus
4[\la_3]\oplus[\la_4]\oplus[\la_5]
\oplus[a\bar{a}]\oplus[b\bar{b}]\oplus [c\bar{c}].$$ Since
$\la_3$ appears at least four times on the LHS, and only twice
on the RHS using the modular invariant $Z_2$, we conclude
using the curious identity \cite[Proposition 3.3]{E1} that indeed
$Z=Z_{(33)}$.
 \eproof

\begin{proposition}
The permutation modular invariant $Z_2$ is sufferable.
\lablth{anothernicetheta}
\end{proposition}
\bproof Let us consider the endomorphism
$\theta=\la_0\oplus\la_4\in\Sigma(\NXN)$. We prove now
that $\theta$ is indeed a dual canonical endomorphism using the
results of \sect{frobsub}. First we equip $\theta$ with a
Q-system structure.
The final Eq.\ in (\ref{mult}) reduces to
\begin{eqnarray}
w_{ij}^p w_{pk}^l= w_{jk}^0 w_{i0}^l{\rm F}_{0p}^{(ijk)l}+
w_{jk}^4 w_{i4}^l{\rm F}_{4p}^{(ijk)l}
\label{mult2}
\end{eqnarray}
with $i,j,k,l,p=0,4$ and the ${\rm F}$'s are the $6j$-symbols
arising from the fusion rules of the the quantum $S_3$ double
model. In the sequel we will use the normalisations from
\cite[\eerf{2.38}]{frs} (cf.\ \cite[\eerf{3.4}]{PZ}: ${\rm
F}_{lj}^{(0jk)l}=\delta_{j,l},$ ${\rm
F}_{ki}^{(i0k)l}=\delta_{k,i},$ ${\rm
F}_{jl}^{(ij0)l}=\delta_{j,l}$. As in \cite[\eerf{3.11}]{PZ},
${\rm F}_{k0}^{(iij)j}=\sqrt{d_k/(d_id_j)}$ whenever
$N_{ij}^k\not=0$ by considering the diagonal case $G_\la=N_\la$.
Then we also have the identity ${\rm F}_{i0}^{(iii)i}{\rm
F}_{0i}^{(iii)i}=$ ${\rm F}_{00}^{(iii)i}=1/d_i$ which holds true
by the pentagon equation
$${\rm F}_{0i}^{(iii)i}{\rm F}_{i0}^{(iii)i}=\sum_s{\rm F}_{is}^{(iii)0}
{\rm F}_{00}^{(isi)i}{\rm F}_{si}^{(iii)0}={\rm
F}_{ii}^{(iii)0}{\rm F}_{00}^{(iii)i}{\rm F}_{ii}^{(iii)0}={\rm
F}_{00}^{(iii)i}.$$ Now we find the unitary matrix ${\rm
F}^{(444)4}$. Note first that since $\la_0,\la_1$ and $\la_4$ are
self conjugate, the entries of ${\rm F}^{(444)4}$ are real
numbers by \cite[\eerf{3.13}]{PZ}. So far we have ${\rm
F}_{40}^{(444)4}=\sqrt{d_4/(d_4d_4)}=1/\sqrt{2}= {\rm
F}_{04}^{(444)4}$, ${\rm F}_{00}^{(444)4}=1/d_4=1/2$, ${\rm
F}_{01}^{(444)4}={\rm F}_{10}^{(444)4}=1/2$. Note that
$$\sum_{s=0,1,4}{\rm F}_{1s}^{(444)4}{\rm F}_{s4}^{(444)4}=0$$
so ${\rm F}_{14}^{(444)4}\not=0$ and we have the following four
equations by the orthogonality relations
\begin{eqnarray}
{\rm F}_{00}^{(444)4}
{\rm F}_{04}^{(444)4}+{\rm F}_{01}^{(444)4} {\rm F}_{14}^{(444)4}+
{\rm F}_{04}^{(4444)4}{\rm F}_{44}^{(444)4}=
\sqrt{2}/4+{\rm F}_{14}^{(444)4}/2+\sqrt{2}{\rm F}_{44}^{(444)4}/2=0
\label{1f}
\end{eqnarray}
\begin{eqnarray}
{\rm F}_{40}^{(444)4}
{\rm F}_{04}^{(444)4}+ {\rm F}_{41}^{(444)4}{\rm F}_{14}^{(444)4}+
{\rm F}_{44}^{(444)4}{\rm F}_{44}^{(444)4}=
1/2+{\rm F}_{14}^{(444)4}\nonumber \\
={\rm F}_{41}^{(444)4}+{\rm F}_{44}^{(444)4}{\rm F}_{44}^{(444)4}=1,
\label{2f}
\end{eqnarray}
\begin{eqnarray}
{\rm F}_{10}^{(444)4}
{\rm F}_{01}^{(444)4}+ {\rm F}_{11}^{(444)4}{\rm F}_{11}^{(444)4}+
{\rm F}_{14}^{(444)4}{\rm F}_{41}^{(444)4}=1,
\label{11f}
\end{eqnarray}
\begin{eqnarray}
{\rm F}_{10}^{(444)4}
{\rm F}_{04}^{(444)4}+ {\rm F}_{11}^{(444)4}{\rm F}_{14}^{(444)4}+
{\rm F}_{14}^{(444)4}{\rm F}_{44}^{(444)4}=0.
\label{14f}
\end{eqnarray}
Now we can insert the value of ${\rm F}_{14}^{(444)4}$ from
\erf{1f} in \erf{2f} to obtain $2{\rm F}_{44}^{(444)4}+3({\rm
F}_{44}^{(444)4})^2=0$. Then ${\rm F}_{44}^{(444)4}=0$ or ${\rm
F}_{44}^{(444)4}=-2/3$. Let us suppose that ${\rm
F}_{44}^{(444)4}=-2/3$. Then replacing this value in \erf{1f} we
find ${\rm F}_{14}^{(444)4}=\sqrt{2}/3$. But then solving
\erf{11f} and \erf{14f} we get different values for ${\rm
F}_{11}^{(444)4}$. Hence ${\rm F}_{44}^{(444)4}=0$ and thus ${\rm
F}_{14}^{(444)4}=-\sqrt{2}/2$ and ${\rm F}_{11}^{(444)4}=1/2$. So
the matrix ${\rm F}^{(444)4}$
$$ {\rm F}^{(444)4}=\pmatrix{{\rm F}_{00}^{(444)4}&{\rm F}_{01}^{(444)4}&
{\rm F}_{04}^{(444)4}\cr
{\rm F}_{10}^{(444)4}&{\rm F}_{11}^{(444)4}&{\rm F}_{14}^{(444)4}\cr
{\rm F}_{40}^{(444)4}&{\rm F}_{41}^{(444)4}&{\rm
F}_{44}^{(444)4}}=\pmatrix{1/2&1/2&
\sqrt{2}/2\cr 1/2&1/2&-\sqrt{2}/2\cr \sqrt{2}/2&-\sqrt{2}/2&0}.$$

Since we get exactly the same (above) equations when we replace
the $\la_3$ by $\la_2$, we conclude that ${\rm F}^{(444)4}={\rm
F}^{(333)3} $, hence since $[\la_0]\oplus[\la_3]$ is a dual
canonical sector by Proposition \ref{nicetheta}, we can conclude
that $[\theta]=[\la_0]\oplus[\la_4]$ is a dual canonical sector
of a braided subfactor $N\subset M$ by \sect{frobsub}. We now
find the modular invariant attached to $\theta$. We compute and
see that $\lan\theta\la_0,\la_0\ran=$ $\lan\theta\la_1,\la_1\ran=$
$\lan\theta\la_2,\la_2\ran=$ $\lan\theta\la_3,\la_3\ran=$
$\lan\theta\la_5,\la_5\ran=1$, $\lan\theta\la_4,\la_4\ran=$
$\lan\theta\la_6,\la_6\ran=$ $\lan\theta\la_7,\la_7\ran=2$,
$\lan\theta\la_0,\la_1\ran=$ $\lan\theta\la_0,\la_2\ran=$
$\lan\theta\la_0,\la_3\ran=$ $\lan\theta\la_0,\la_5\ran=$
$\lan\theta\la_0,\la_6\ran=$ $\lan\theta\la_0,\la_7\ran=0$,
$\lan\theta\la_0,\la_4\ran=1$, $\lan\theta\la_1,\la_2\ran=$
$\lan\theta\la_1,\la_3\ran=$ $\lan\theta\la_1,\la_5\ran=$
$\lan\theta\la_1,\la_6\ran=$ $\lan\theta\la_1,\la_7\ran=0$,
$\lan\theta\la_1,\la_4\ran=1$, $\lan\theta\la_2,\la_3\ran=$
$\lan\theta\la_2,\la_5\ran=$ $\lan\theta\la_3,\la_5\ran=$
$\lan\theta\la_6,\la_7\ran=1$. Then we conclude that $[\iota\la_0]$,
$[\iota\la_1]$ and $[\iota\la_2]$ are distinct irreducible $N$-$M$
sectors.
We also have $[\iota\la_5]=[\iota\la_3]=[\iota\la_2]$,
$[\iota\la_6]=[a]\oplus[b]$
and $[\iota\la_7]=[a]\oplus[c]$ with $[a]$, $[b]$ and $[c]$ new
irreducible $N$-$M$ sectors. Then $\Tr(Z)=6$ so $Z=Z_{(22)},
Z_{(33)},Z_2,Z_7,Z_{22} \hbox{ or } Z_{33}$. Since
$\lan\theta\la_0,\la_1\ran=0$, $Z$ cannot be $Z_7$, $Z_{22}$,
$Z_{33}$ by \cite[\eerf{37}]{BEK2}. Also because
$\lan\theta\la_0,\la_2\ran=0=\lan\theta \la_0,\la_3\ran$, $Z$ cannot be
$Z_{(22)}$ and $Z_{(33)}$ again by \cite[\eerf{37}]{BEK2}. Hence
$Z=Z_2$. \eproof

\begin{remark}
We can likewise prove that
$\theta=\la_0\oplus\la_5\in\Sigma(\NXN)$ is a dual canonical
endomorphism. It also produces the permutation modular invariant
$Z_2$
\end{remark}
\begin{proposition}
The modular invariants $Z_{(22)}$ and $Z_{(44)}$ are sufferable.
\lablth{ext}
\end{proposition}
\bproof The modular invariant $Z_{(22)}$ is sufferable because is
a product of sufferable modular invariants $Z_{(22)}=Z_2
Z_{(33)}Z_2$ (see previous Propositions \ref{nicetheta} and
\ref{anothernicetheta}). We now prove that $Z_{(44)}$ is
sufferable. We already know that $Z_{(22)}$ is sufferable and
moreover it is a type I modular invariant with
$[\theta]=[\la_0]\oplus[\la_2]$. The global indices of $\MXM$,
$\MXMpm$ and $\MXMo$ are $\omega=36, \omega_\pm=36/3=12,
\omega_0=12^2/36=4$ respectively. We compute the
$\mathcal{X}^\pm$-chiral systems and conclude that the
irreducible decompositions are as follows: $[\a_0]$,
$[\a_1^\pm]$, $[\a_2^\pm]= [\a_0]\oplus[\a_1^\pm]$,
$[\a_3^\pm]=[\a_4^\pm]=[\a_5^\pm]$,
$[\a_6^\pm]=[\a_6^{\pm(1)}]\oplus[\a_6^{\pm(2)}]$ and
$[\a_7^\pm]= [\a_6^{\pm(1)}]\oplus[\a_7^{\pm(2)}]$. From the
entries of the modular invariant $Z_{(22)}$, we see that
$[\a_1^+]=[\a_1^-]$ denoted henceforth by $[\a_1]$,
$[\a_6^{+(2)}]=[\a_6^{-(2)}]$ denoted from now by $[\a_6^{(2)}]$
and similarly $[\a_7^{(2)}]$. Thus $\MXMo=\{\a_0,\a_1,
\a_6^{(2)}, \a_7^{(2)}\}$. Also the $\mathcal{X}^\pm$-chiral
systems are $\MXMpm=\{\a_0,\a_1, \a_6^{(2)}, \a_7^{(2)},\a_3^\pm,
\a_6^{\pm(1)}\}$. In particular, the branching coefficient matrix
$b=[b_{\tau,\la}]$ with $b_{\tau,\la}=\lan\tau,\a^\pm_\la\ran$,
$\tau\in\MXMo$:
\begin{eqnarray}
b={\scriptstyle\pmatrix{1&0&1&0&0&0&0&0\cr 0&1&1&0&0&0&0&0\cr 
0&0&0&0&0&0&1&0\cr 0&0&0&0&0&0&0&1}}.
\label{bb}
\end{eqnarray}

By \cite{BE6}, there is a non-degenerate braiding on the
neutral system $\MXMo$ of $Z_{(22)}$. Let $S^{{\rm ext}}$ and
$T^{\rm ext}$ be the modular matrices of the extended system
$\MXMo$. So by \cite[Theorem 6.5]{BE4}, $S^{{\rm ext}}b=bS,
T^{{\rm ext}}b=bT.$ Since we are in the fortunate situation that
$$bb^t={\scriptstyle\pmatrix{2&1&0&0\cr 1&2&0&0\cr 0&0&1&0\cr 0&0&0&1}}$$ is
invertible, we can solve the above equations and get $S^{{\rm
ext}}$ and $ T^{{\rm ext}}$ uniquely. They are the modular data
of the quantum double of $\bbZ_2$ presented in \cite{BE5} or
\cite[Page 285]{BE4} or \erf{stgroup} by permuting the second $v$
and fourth $c$ labels.
Since $$M={\scriptstyle\pmatrix{1&0&0&0\cr 0&0&1&0\cr 0&1&0&0\cr
0&0&0&1}}$$ is sufferable by \cite{BE5}, so is $b^tMb=Z_{(44)}$.
If $N\subset M$ produces $Z_{(22)}$,
then $N\subset M\rtimes_v \bbZ_2$ realises $Z_{(44)}$ whose dual
canonical endomorphism is $\theta=\sigma_{\theta^{\rm ext}}$ with
$[\theta^{\rm ext}]=[\la_0]\oplus [\a_7^{(2)}]$ and $\sigma$ is
the $\sigma$-restriction \cite{BE5} defined by
$\sigma_{\beta}=\bar{\iota}\beta\iota$. We compute and find that
$[\theta]=[\la_0]\oplus[\la_2]\oplus[\la_7]$ is the dual
canonical sector of $N\subset M\rtimes\bbZ_2$. \eproof
\begin{table}
{\scriptsize
$$\begin{array}{|c||c|c|c|c|c|c|c|c|c|c|c|} \hline
\times &Z_1&Z_2&Z_5&Z_7&Z_{22}& Z_{33}& Z_{44}&Z_{55}& Z_{32}&Z_{42}&Z_{52}\\ \hline\hline
Z_2&Z_2&Z_1&Z_7&Z_5&Z_{32}&Z_{23}& Z_{54}& Z_{45}& Z_{33}& Z_{34}& Z_{35}\\ \hline
Z_5&Z_{5}& Z_{7}&2 Z_{5}&2Z_{7}&2 Z_{22}&2 Z_{33}& Z_{24}& Z_{35}&2Z_{23}&2 Z_{24}&2Z_{25}\\  \hline
Z_7& Z_{7}& Z_{5}&2Z_{7}&2Z_{5}&2Z_{32}&2Z_{23}& Z_{34}& Z_{25}& Z_{33}&2 Z_{34}&2 Z_{35}\\  \hline
Z_{22}& Z_{22}& Z_{23}&2Z_{22}&2 Z_{23}&6 Z_{22}&2 Z_{23}&3 Z_{24}& Z_{25}&6 Z_{23}&6 Z_{24}&6 Z_{25} \\  \hline
Z_{33}& Z_{33}& Z_{32}&2Z_{33}&2 Z_{32}&2 Z_{32}&6 Z_{33}& Z_{34}&3 Z_{35}&2 Z_{33}&2 Z_{34}&
2 Z_{35}\\ \hline
Z_{44}& Z_{44}& Z_{42}& Z_{42}& Z_{43}&3Z_{42}& Z_{43}&3 Z_{44}&2 Z_{45}&3 Z_{43}&2 Z_{44}&2Z_{45}\\ \hline
Z_{55}& Z_{55}& Z_{54}& Z_{53}& Z_{52}& Z_{52}&3 Z_{53}&2 Z_{52}&3 Z_{55}& Z_{53}& Z_{54}& Z_{55}\\ \hline
Z_{23}& Z_{23}& Z_{22}& 2Z_{23}&2 Z_{22}&2 Z_{22}&6 Z_{23}& Z_{24}&3 Z_{25}&
2Z_{23}& 2Z_{24}&2 Z_{25}\\ \hline
Z_{24}& Z_{24}& Z_{25}& Z_{22}& Z_{23}&3 Z_{22}& Z_{23}&3 Z_{24}&2 Z_{25}&3 Z_{23}&3 Z_{24}& 3Z_{25}\\ \hline
Z_{25}& Z_{25}& Z_{24}& Z_{23}& Z_{22}& Z_{22}&3 Z_{23}&2 Z_{24}&3 Z_{25}& Z_{23}& Z_{24}& Z_{25}\\  \hline
Z_{32}& Z_{32}& Z_{33}&2 Z_{32}&2 Z_{33}&6 Z_{32}&2 Z_{33}&3 Z_{34}& Z_{35}&6 Z_{33}&6 Z_{34}&6 Z_{35}\\  \hline
Z_{34}& Z_{34}& Z_{35}& Z_{32}& Z_{33}&3 Z_{32}& Z_{33}&3 Z_{34}&2 Z_{35}&3 Z_{33}&3 Z_{34}&3 Z_{35}\\  \hline
Z_{35}& Z_{35}& Z_{34}& Z_{33}& Z_{32}& Z_{32}&3 Z_{33}&2 Z_{34}&3 Z_{35}& Z_{33}& Z_{34}& Z_{35}\\  \hline
Z_{42}& Z_{42}& Z_{43}&2 Z_{42}&2 Z_{43}&6 Z_{42}&2 Z_{43}&3 Z_{44}& Z_{45}&6 Z_{43}&6 Z_{44}&6 Z_{45}\\  \hline
Z_{43}& Z_{43}& Z_{42}&3 Z_{42}&3 Z_{43}&2 Z_{42}&6 Z_{43}& Z_{44}&3 Z_{45}&2 Z_{43}&2 Z_{44}&2 Z_{45}\\  \hline
Z_{45}& Z_{45}& Z_{44}& Z_{43}& Z_{42}& Z_{42}&3 Z_{43}&2Z_{44}&3Z_{45}& Z_{43}& Z_{44}& Z_{45}\\  \hline
Z_{52}& Z_{52}& Z_{53}&2Z_{52}&2 Z_{53}&6Z_{52}&2 Z_{53}&3 Z_{54}& Z_{55}&6 Z_{53}&6 Z_{54}&6 Z_{55}\\  \hline
Z_{53}& Z_{53}& Z_{52}&2Z_{53}&2 Z_{52}&2 Z_{52}&6 Z_{53}& Z_{54}&3 Z_{55}&2 Z_{53}&2 Z_{54}&2 Z_{55}\\  \hline
Z_{54}& Z_{54}& Z_{55}& Z_{52}& Z_{53}&3 Z_{52}& Z_{53}&3 Z_{54}&2 Z_{55}&3 Z_{53}&3 Z_{54}&3 Z_{55}\\  \hline
Z_{(32)}& Z_{(32)}& Z_{(33)}& Z_{32}& Z_{33}&3 Z_{32}& Z_{33}& Z_{34}+ Z_{54}& Z_{55}&3 Z_{33}&3 Z_{34}&3 Z_{35}\\  \hline
Z_{(45)}& Z_{(45)}& Z_{(55)}& Z_{52}& Z_{53}&3 Z_{52}& Z_{53}& Z_{54}+ Z_{34}& Z_{35}&3 Z_{53}&3 Z_{54}&3 Z_{55}\\  \hline
Z_{(44)}& Z_{(44)}& Z_{42}& Z_{42}& Z_{43}&3 Z_{42}& Z_{43}& Z_{44}+ Z_{24}& Z_{25}&3 Z_{43}&3 Z_{44}&3 Z_{45}\\  \hline
Z_{(22)}& Z_{(22)}& Z_{(23)}& Z_{22}& Z_{23}&3 Z_{22}& Z_{23}& Z_{44}+ Z_{45}& Z_{45}&3 Z_{23}&3Z_{24}&3 Z_{25}\\  \hline
Z_{(33)}& Z_{(33)}& Z_{(32)}& Z_{33}& Z_{32}& Z_{32}&3 Z_{33}& Z_{54}& Z_{55}+ Z_{35}& Z_{33}& Z_{34}& Z_{35}\\  \hline
Z_{(55)}& Z_{(55)}& Z_{(45)}& Z_{53}& Z_{52}& Z_{52}&3 Z_{53}& Z_{34}& Z_{55}+ Z_{35}& Z_{53}& Z_{54}& Z_{55}\\  \hline
Z_{(54)}& Z_{(54)}& Z_{(44)}& Z_{43}& Z_{42}& Z_{42}&3 Z_{43}& Z_{44}& Z_{25}+ Z_{45}& Z_{43}& Z_{44}& Z_{45}\\  \hline
Z_{(23)}& Z_{(23)}& Z_{(22)}& Z_{23}& Z_{22}& Z_{22}&3 Z_{23}& Z_{44}&
 Z_{45}+ Z_{25}& Z_{23}& Z_{24}& Z_{25}\\  \hline\hline
\end{array}$$
}
\caption{Fusion $Z_a Z_b^t$ of $S_3$ modular invariants from Table \ref{psuff}}
\label{fusions3}
\end{table}
 In Tables \ref{fusions3} and \ref{fusions33} we present
products $Z_aZ_b^t$ of the 28 matrices of Table \ref{psuff}
together with their (unique) decomposition into normalised modular
invariants.
\begin{table}
{\scriptsize
$$\begin{array}{||c||c|c|c|c|c|c|c|c|c|c|} \hline
\times&Z_{23}&Z_{43}& Z_{53}& Z_{24}& Z_{34}& Z_{54}& Z_{25}& Z_{35}& Z_{45}& Z_{(23)}\\ \hline\hline
Z_2&Z_{22}&Z_{24}& Z_{25}& Z_{52}& Z_{53}& Z_{55}& Z_{42}& Z_{43}& Z_{44}& Z_{(22)}\\ \hline
Z_5&2Z_{32}&2 Z_{34}&2Z_{35}& Z_{22}& Z_{23}& Z_{25}& Z_{32}& Z_{33}& Z_{34}& Z_{32}\\  \hline
Z_7&2 Z_{22}&2 Z_{24}&2 Z_{25}& Z_{32}& Z_{33}& Z_{35}& Z_{22}& Z_{23}& Z_{24}& Z_{22}\\  \hline
Z_{22}&2 Z_{22}&2 Z_{24}&2 Z_{25}&3 Z_{22}&3 Z_{23}&3 Z_{25}& Z_{22}& Z_{23}& Z_{24}& Z_{22}\\  \hline
Z_{33}&6 Z_{32}&6 Z_{34}&6 Z_{35}& Z_{32}& Z_{33}& Z_{35}&3 Z_{32}&3 Z_{33}&
3Z_{34}&3 Z_{32}\\ \hline
Z_{44}& Z_{42}& Z_{44}& Z_{45}&3 Z_{42}&3 Z_{43}&3 Z_{45}&2 Z_{42}&2 Z_{43}&
2Z_{44}& Z_{44}\\ \hline
Z_{55}&2 Z_{52}&3 Z_{54}&3 Z_{55}&2 Z_{52}&2 Z_{53}&2 Z_{53}&3Z_{52}&3 Z_{53}&
3Z_{54}& Z_{52}+ Z_{54}\\ \hline
Z_{23}&6 Z_{22}&6 Z_{24}&6 Z_{25}& Z_{22}& Z_{23}& Z_{25}&3 Z_{22}&3 Z_{23}&
3Z_{24}&3 Z_{22}\\ \hline
Z_{24}& Z_{22}& Z_{24}& Z_{25}&3 Z_{22}&3 Z_{23}&3 Z_{25}&2Z_{22}&2 Z_{23}&
2Z_{24}& Z_{24}\\ \hline
Z_{25}&3 Z_{22}&3 Z_{24}&3 Z_{25}&2 Z_{22}&2 Z_{23}&2 Z_{25}&3 Z_{22}&3 Z_{23}&
3Z_{24}& Z_{22}+ Z_{24}\\  \hline
Z_{32}&2 Z_{32}&2 Z_{34}&2Z_{35}&3 Z_{32}&3 Z_{33}&3 Z_{35}& Z_{32}& Z_{33}& Z_{34}& Z_{32}\\  \hline
Z_{34}& Z_{32}& Z_{34}& Z_{35}&3 Z_{32}&3Z_{33}&3 Z_{35}&2 Z_{32}&2 Z_{33}&
2Z_{34}& Z_{34}\\  \hline
Z_{35}&3 Z_{32}&3 Z_{34}&3 Z_{35}&2 Z_{32}&2 Z_{33}&2 Z_{35}&3 Z_{32}&3 Z_{33}&
3Z_{34}& Z_{32}+ Z_{34}\\  \hline
Z_{42}&2 Z_{42}&2 Z_{44}&2 Z_{45}&3 Z_{42}&3 Z_{43}&3 Z_{45}& Z_{42}& Z_{43}& Z_{44}& Z_{42}\\  \hline
Z_{43}&6 Z_{42}&6 Z_{44}&6 Z_{45}& Z_{42}& Z_{43}& Z_{45}&3 Z_{42}&3 Z_{43}&
3Z_{44}&3 Z_{42}\\  \hline
Z_{45}&3 Z_{42}&3 Z_{44}&3 Z_{45}&2 Z_{42}&2 Z_{43}&2 Z_{45}&3 Z_{42}&3 Z_{43}&
3Z_{44}& Z_{42}+ Z_{44}\\  \hline
Z_{52}&2 Z_{52}&2Z_{54}&2 Z_{55}&3 Z_{52}&3 Z_{53}&3 Z_{55}& Z_{52}& Z_{53}& Z_{54}& Z_{52}\\  \hline
Z_{53}&6 Z_{52}&6 Z_{54}&6 Z_{55}& Z_{52}& Z_{53}& Z_{55}&3 Z_{52}&3 Z_{53}&
3Z_{54}&3 Z_{52}\\  \hline
Z_{54}& Z_{52}& Z_{54}& Z_{55}&3 Z_{52}&3 Z_{53}&3 Z_{55}&2 Z_{52}&2 Z_{53}&
2Z_{54}& Z_{54}\\  \hline
Z_{(32)}& Z_{32}& Z_{34}& Z_{35}& Z_{52}+ Z_{32}& Z_{53}+ Z_{33}& Z_{35}+ Z_{55}& Z_{52}& Z_{53}& Z_{54}& Z_{(32)}\\  \hline
Z_{(45)}& Z_{52}& Z_{54}& Z_{55}& Z_{52}+ Z_{32}& Z_{53}+ Z_{33}& Z_{55}+  Z_{35}& Z_{32}& Z_{33}& Z_{34}& Z_{(45)}\\  \hline
Z_{(44)}& Z_{42}& Z_{44}& Z_{45}& Z_{22}+ Z_{42}& Z_{23}+ Z_{43}& Z_{45}+  Z_{25}& Z_{22}& Z_{23}& Z_{24}& Z_{(44)}\\  \hline
Z_{(22)}& Z_{22}& Z_{24}& Z_{25}& Z_{42}+ Z_{22}& Z_{43}+ Z_{23}& Z_{45}+ Z_{25}& Z_{42}& Z_{43}& Z_{44}& Z_{(22)}\\  \hline
Z_{(33)}&3 Z_{32}&3 Z_{34}&3 Z_{35}& Z_{52}& Z_{53}& Z_{55}& Z_{52}+ \atop Z_{32}& Z_{53}+ \atop Z_{33}& Z_{54}+ Z_{34}& Z_{32}+ Z_{(32)}\\  \hline
Z_{(55)}&3 Z_{52}&3 Z_{54}&3 Z_{55}& Z_{32}& Z_{33}& Z_{35}& Z_{52}+ \atop Z_{32}& Z_{53}+ \atop Z_{33}& Z_{54}+  Z_{34}& Z_{52}+ Z_{(45)}\\  \hline
Z_{(54)}&3 Z_{42}&3 Z_{44}&3 Z_{45}& Z_{22}& Z_{23}& Z_{25}& Z_{42}+ \atop Z_{22}& Z_{43}+ \atop Z_{23}& Z_{44}+  Z_{24}& Z_{42}+ Z_{(44)}\\  \hline
Z_{(23)}&3 Z_{22}&3 Z_{24}&3 Z_{25}& Z_{42}&Z_{43}& Z_{45}& Z_{42}+ \atop Z_{22}& Z_{43}+ \atop Z_{23}& Z_{44}+  Z_{24}& Z_{22}+ Z_{(22)}\\ \hline\hline
\end{array}$$
}
{\scriptsize
$$\begin{array}{||c||c|c|c|c|c|c|c|} \hline
\times&Z_{(54)}& Z_{(44)}& Z_{(22)}& Z_{(33)}& Z_{(55)}& Z_{(45)}& Z_{(32)}\\ \hline\hline
Z_2&Z_{(44)}& Z_{(45)}& Z_{(32)}& Z_{(23)}& Z_{(54)}& Z_{(55)}& Z_{(33)}\\ \hline
Z_5& Z_{34}& Z_{24}& Z_{22}& Z_{33}& Z_{35}& Z_{25}& Z_{23}\\  \hline
Z_7& Z_{24}& Z_{34}& Z_{32}& Z_{23}& Z_{25}& Z_{35}& Z_{33}\\  \hline
Z_{22}& Z_{24}&3 Z_{24}&3 Z_{22}& Z_{23}& Z_{25}&3 Z_{25}&3 Z_{23}\\  \hline
Z_{33}&3 Z_{34}& Z_{34}& Z_{32}&3 Z_{33}&3 Z_{35}& Z_{35}& Z_{33}\\ \hline
Z_{44}& Z_{42}& Z_{44}+ Z_{42}& Z_{44}+ Z_{42}& Z_{44}& Z_{43}& Z_{43}+ Z_{45}& Z_{43}+ Z_{45}\\ \hline
Z_{55}& Z_{52}+ Z_{54}& Z_{52}& Z_{54}& Z_{53}+ Z_{55}& Z_{53}+ Z_{55}& Z_{53}& Z_{55}\\ \hline
Z_{23}& Z_{24}& Z_{24}& Z_{22}&3 Z_{23}&3 Z_{25}& Z_{25}& Z_{23}\\ \hline
Z_{24}& Z_{22}& Z_{22}+ Z_{24}& Z_{22}+ Z_{24}& Z_{25}& Z_{23}& Z_{23}+ Z_{25}& Z_{23}+ Z_{25}\\ \hline
Z_{25}& Z_{22}+ Z_{24}& Z_{22}& Z_{24}& Z_{23}+ Z_{25}& Z_{23}+ Z_{25}& Z_{23}& Z_{25}\\  \hline
Z_{32}& Z_{34}&3 Z_{34}&3 Z_{32}& Z_{33}& Z_{35}&3 Z_{35}&3 Z_{33}\\  \hline
Z_{34}& Z_{32}& Z_{32}+ Z_{34}& Z_{32}+ Z_{34}& Z_{35}& Z_{33}& Z_{33}+ Z_{35}& Z_{33}+ Z_{35}\\  \hline
Z_{35}& Z_{32}+ Z_{34}& Z_{32}& Z_{34}& Z_{33}+ Z_{35}& Z_{33}+ Z_{25}& Z_{33}& Z_{35}\\  \hline
Z_{42}& Z_{44}&3 Z_{44}&3 Z_{42}& Z_{43}& Z_{45}&3 Z_{45}&3 Z_{45}\\  \hline
Z_{43}&3 Z_{44}& Z_{44}& Z_{42}&3 Z_{43}&3 Z_{45}& Z_{45}& Z_{43}\\  \hline
Z_{45}& Z_{42}+ Z_{44}& Z_{42}& Z_{44}& Z_{43}+ Z_{45}& Z_{43}+ Z_{45}& Z_{43}& Z_{45}\\  \hline
Z_{52}& Z_{54}&3 Z_{54}&3 Z_{52}& Z_{53}& Z_{55}&3 Z_{55}&3 Z_{53}\\  \hline
Z_{53}&3 Z_{54}& Z_{54}& Z_{52}&3 Z_{53}&3 Z_{55}& Z_{55}& Z_{53}\\  \hline
Z_{54}& Z_{52}& Z_{54}+ Z_{52}& Z_{54}+ Z_{52}& Z_{55}& Z_{53}& Z_{55}+ Z_{53}& Z_{55}+ Z_{53}\\  \hline
Z_{(32)}& Z_{(45)}& Z_{34}+ Z_{(45)}& Z_{(32)}+ Z_{32}& Z_{(33)}& Z_{(55)}& Z_{35}+ Z_{(55)}& Z_{33}+ Z_{(33)}\\  \hline
Z_{(45)}& Z_{(32)}& Z_{54}+ Z_{(32)}& Z_{52}+ Z_{(45)}& Z_{(55)}& Z_{(33)}& Z_{55}+ Z_{(33)}& Z_{53}+ Z_{(55)}\\  \hline
Z_{(44)}& Z_{(22)}& Z_{44}+ Z_{(22)}& Z_{42}+ Z_{(45)}&Z_{(54)}& Z_{(23)}& Z_{45}+ Z_{(23)}& Z_{43}+ Z_{(54)}\\  \hline
Z_{(22)}& Z_{(44)}& Z_{24}+ Z_{(44)}& Z_{22}+ Z_{(22)}&Z_{(23)}& Z_{(54)}& Z_{25}+ Z_{(54)}& Z_{23}+ Z_{(23)}\\  \hline
Z_{(33)}& Z_{34}+ Z_{(45)}& Z_{(45)}& Z_{(32)}& Z_{33}+Z_{(33)}&Z_{35}+Z_{(55)}& Z_{(55)}& Z_{(33)}\\  \hline
Z_{(55)}& Z_{54}+ Z_{(32)}& Z_{(32)}& Z_{(45)}& Z_{53}+ Z_{(55)}& Z_{55}+Z_{(33)}& Z_{(33)}& Z_{(55)}\\  \hline
Z_{(54)}& Z_{54}+ Z_{(22)}& Z_{(22)}& Z_{(44)}& Z_{43}+ Z_{(54)}& Z_{45}+Z_{(22)}& Z_{(23)}& Z_{(54)}\\  \hline
Z_{(23)}& Z_{24}+Z_{(44)}& Z_{(44)}& Z_{(22)}& Z_{23}+ Z_{(23)}& Z_{25}+ Z_{(54)}& Z_{(54)}& Z_{(23)}\\  \hline\hline
\end{array}$$
}
\caption{Fusion $Z_a Z_b^t$ of $S_3$ modular invariants from Table \ref{psuff} (cont.)}
\label{fusions33}
\end{table}
Since by \sect{frobsub} products of sufferable modular invariants
are sufferable, we conclude that all the 28 nimble modular
invariants of the Table \ref{psuff} are sufferable as follows. We
already know that $Z_1$, $Z_2$, $Z_5$, $Z_{55}$, $Z_{(44)}$,
$Z_{(22)}$, $Z_{(33)}$ are sufferable. The modular invariant $Z_7$
is sufferable because $Z_7=Z_2Z_5$. Then
$Z_{(23)}=Z_{(22)}Z_{(33)}$ is sufferable and thus so is its
conjugate $Z_{(32)}=Z_{(23)}^t$. As $Z_7Z_{55}=Z_{25}$ we also
conclude that $Z_{25}$ and its conjugate $Z_{52}$ are sufferable;
so is $Z_{22}$ by the type I parent theorem \cite{BE4}. As
$Z_5Z_{52}= Z_{32}$, $Z_{32}$ and its conjugate $Z_{23}$ are
sufferable; so is $Z_{33}$ by the type I theorem. As
$Z_{55}Z_5=Z_{35}$, $Z_{35}$ and its conjugate $Z_{53}$ are also
sufferable. Then $Z_5Z_{(44)}=Z_{24}$, $ Z_{(44)} Z_{33}= Z_{43}$
and $Z_{55} Z_{24}=Z_{54}$, so $Z_{24}, Z_{42}$, $Z_{34}, Z_{43},
Z_{54},Z_{45}$ are sufferable. Again by the type I parent theorem
\cite{BE4}, $Z_{44}$ is sufferable. As $Z_{(44)}Z_{(33)}=Z_{(54)}$
we conclude that $Z_{(54)}$ and its conjugate $Z_{(45)}$ are
sufferable. The same with $Z_{(55)}$ since
$Z_{(55)}=Z_{(33)}Z_{(54)}$.

\begin{remark}{\rm  In the framework of \cite{ost}, every (conjugacy class)
subgroup of $S_3\times S_3$ gives rise to a sufferable modular
invariant (although different subgroups may well give the same
modular invariant). The trace of $Z$ attached to $H$ is easily
computed: we choose representatives $\{g\}$ of the double cosets
$\Delta(S_3)\setminus S_3\times S_3/H$ then
$$\Tr(Z)=\sum_g\#\hbox{irred}\Big(\Delta(S_3)\cap
gHg^{-1}\Big).$$ With our approach we managed to show that the
trace 4 modular invariants $Z_{(32)}$ and $Z_{(23)}$ are
sufferable which do not appear in Ostrik's Table 4.4 in
\cite{ost} where there are 20 subgroups of $S_3\times S_3$.
However, there are two more (conjugacy classes) subgroups
$H_{21}=(\bbZ_3\times 0)\cdot (\Delta(\bbZ_2))$ and
$H_{22}=(\Delta(\bbZ_2))\cdot(0\times\bbZ_3)$, both
isomorphic to $S_3$ which are missing in \cite[Table 4.1]{ost}, which produce
$\Tr(Z)=4$ modular invariants. Therefore $Z=Z_{(23)}, Z_{(32)}$.}
\end{remark}

For the sake of completeness, we summarize the 48 modular invariants
in each category: nimless, spurious nimble and sufferable.
The span of the sufferable modular invariants is $\{S,T\}^\prime$.
\begin{corollary}
\begin{enumerate}
\item The nimless $S_3$ modular invariants are the 14 matrices\\
{\footnotesize
$\begin{array}{|c|c|c|c|c|c|c|}\hline\hline
Z_6&Z_{11}& Z_{14}&Z_{15}&Z_{41}& Z_{51}& Z_{(21)}\\ \hline
Z_{(3)}&Z_{(4)}& Z_{(1)}& Z_{(12)}& Z_{(2)}&Z_{(16)}& Z_{(61)}\\
\hline\hline\end{array}$
}
\item The nimble but insufferable modular invariants are the
6 matrices\\
{\footnotesize
$\begin{array}{|c|c|c|c|c|c|}\hline\hline
Z_3& Z_4& Z_{12}&
Z_{13}& Z_{21}& Z_{31}\\ \hline\hline
\end{array}$
}
\item The sufferable modular invariants are the 28 matrices\\
{\footnotesize $\begin{array}{|c|c|c|c|c|c|c|c|c|c|c|c|c|c|}
\hline\hline Z_1&Z_2& Z_7& Z_5& Z_{22}& Z_{23}& Z_{32}& Z_{33}&
Z_{(22)}& Z_{(33)}& Z_{55}& Z_{24}& Z_{25}& Z_{34}\\ \hline
Z_{35}& Z_{42}&Z_{43}&Z_{44}& Z_{45}&
Z_{52}&Z_{53}&Z_{54}&Z_{(32)}& Z_{(45)}& Z_{(44)}& Z_{(55)}&
Z_{(54)}& Z_{(23)}\\ \hline\hline
\end{array}$
}
\end{enumerate}
\lablth{list}
\end{corollary}

\subsection{The full systems for the sufferable modular invariants}
\label{pictures}

The global index $\omega=\sum_\eta S_{0\eta}^2/S_{00}^2$ of our system
$\NXN$ is 36.

\subsubsection*{Case $Z_1$}

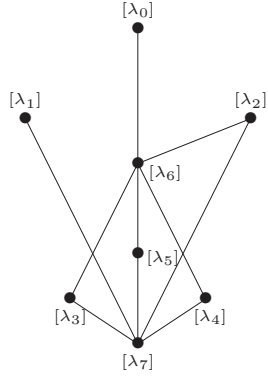
\begin{figure}[htb]
\begin{center}
\unitlength 0.6mm
\begin{picture}(70,90)
\thinlines

\put(20,0){\makebox(0,0){$\bullet$}}
\put(20,20){\makebox(0,0){$\bullet$}}
\put(20,40){\makebox(0,0){$\bullet$}}
\put(20,70){\makebox(0,0){$\bullet$}}

\put(35,10){\makebox(0,0){$\bullet$}}
\put(45,50){\makebox(0,0){$\bullet$}}
\put(5,10){\makebox(0,0){$\bullet$}}
\put(-5,50){\makebox(0,0){$\bullet$}}

\path(20,70)(20,40)(20,20)(20,0)(35,10)(20,40)(5,10)(20,0)
\path(-5,50)(20,0)(45,50)(20,40)

\put(25.5,19){\makebox(0,0){{\tiny $[\la_5]$}}}
\put(20,74){\makebox(0,0){{\tiny $[\la_0]$}}}
\put(26,38){\makebox(0,0){{\tiny $[\la_6]$}}}
\put(20,-4){\makebox(0,0){{\tiny $[\la_7]$}}}
\put(-5,54){\makebox(0,0){{\tiny $[\la_1]$}}}
\put(45,54){\makebox(0,0){{\tiny $[\la_2]$}}}
\put(5,6){\makebox(0,0){{\tiny $[\la_3]$}}}
\put(36,6){\makebox(0,0){{\tiny $[\la_4]$}}}

\end{picture}
\end{center}
\caption{$Z_1$, fusion graph of $[\a_{\la_6}]$}
\label{qsZ1}
\end{figure}

The structure of the induced full $M$-$M$ system for the trivial
modular invariant $Z_1$ is given by the original Verlinde algebra
of modular data. Since
$\lan\a_\la^+,\a_\mu^-\ran=\delta_{\la,\mu}$,
$\a_\la^+=\a_\la^-$, and all disjoint. Then $\la\mapsto\a_\la$ is
an isomorphism. Hence $\MXM\simeq\MXMpm\simeq\MXMo\simeq\NXN$
for any subfactor $N \subset M$ realising the trivial modular invariant.
We display the fusion graph of $[\la_6]$ in \fig{qsZ1}. The full
systems for the permutation invariant $Z_2$ are similar, obtained by permuting
$\la_2$ and $\la_3$ in those for $Z_1$.

\subsubsection*{Case $Z_5$}
\label{chiralz5}

The simple current invariant $Z_5$ yields 20 irreducible $M$-$M$
sectors, whose fusion is non-commutative by \cite[Corollary
6.9]{BEK1}. The global indices of $\omega_\pm=18, \omega_0=9.$ We
compute now the chiral systems $\MXMpm$ and then the branching
coefficients. By the fusion rules, Frobenius reciprocity and
homomorphism property of the $\a$-induction, we have
$$\lan \a_\la^+,\a_\mu^+\ran=\lan \a_\la^+\a_\mu^+,\a_0\ran=
\sum_\eta N_{\la\mu}^\eta Z_{\eta 0}= N_{\la\mu}^0+ N_{\la\mu}^1=
\lan \a_\la^-,\a_\mu^-\ran .$$ We therefore find, using the
original Verlinde fusion matrices, that $[\a_0^\pm]=
[\a_1^\pm],\quad [\a_2^\pm]=$
$[\a_2^{\pm(1)}]\oplus[\a_2^{\pm(2)}],\quad
[\a_2^\pm]=[\a_2^{\pm(1)}]\oplus[\a_2^{\pm(2)}],$
$[\a_4^\pm]=[\a_4^{\pm(1)}]\oplus[\a_4^{\pm(2)}],\quad
[\a_5^\pm]=[\a_5^{\pm(1)}]\oplus[\a_5^{\pm(2)}],\quad
[\a_6^\pm]=[\a_7^\pm],$ with $[\a_0^\pm], [\a_j^{\pm(i)}]$ and
$[\a_6^\pm]$ irreducible sectors ($i=1,2; j=2,3,4,5$). Since the
dimension function is additive, $d_{\a_j^{\pm(i)}}=1$ for all
$i=1,2; j=2,3,4,5$. Of course $d_{\a_6^\pm}=3$. We conclude that
the commutative neutral system as sectors is formed with 9
automorphisms $[\a_0], [\a_j^{(i)}],\quad \hbox{with}\ i=1,2;
j=2,3,4,5$ which is either $\bbZ_3\times\bbZ_3$ or $\bbZ_9$.
However, $\hbox{irred}([\a_j^{(i)}]^n)\in\{[\a_0], [\a_j^{(i)}]\},
i=1,2; j=2,3,4,5.$ This implies that it cannot be $\bbZ_9$
therefore is has to be $\bbZ_3\times\bbZ_3$. Hence the system
$\MXMpm=\MXMo\cup\{\a_6^\pm\},$ with the other fusion rules given
by $[\a_j^{(i)}] [\a_6^\pm]=[\a_6^\pm] [\a_j^{(i)}]= [\a_6^\pm],
i=1,2;j=2,3,4,5.$

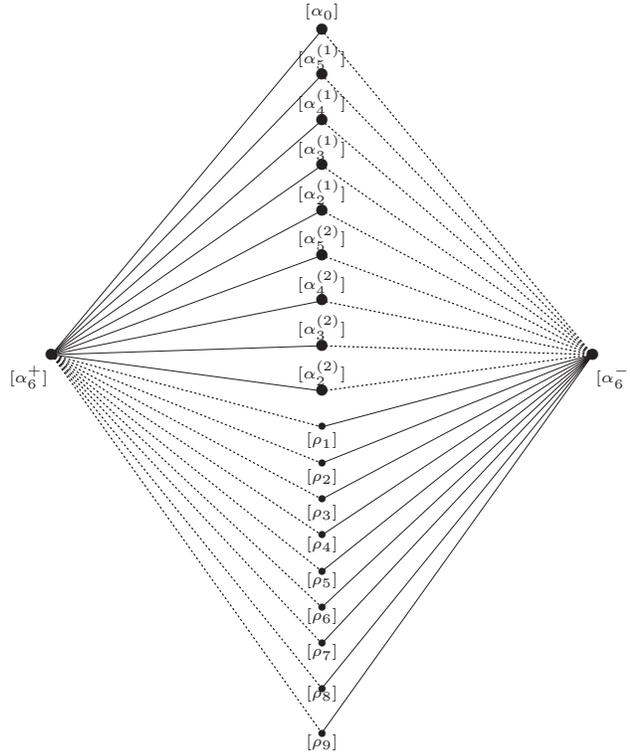
\begin{figure}[htb]
\begin{center}
\unitlength 0.6mm
\begin{picture}(70,160)
\thinlines

\put(20,-4){\makebox(0,0){{\tiny $\bullet$}}}
\put(20,6){\makebox(0,0){{\tiny $\bullet$}}}
\put(20,16){\makebox(0,0){{\tiny $\bullet$}}}
\put(20,24){\makebox(0,0){{\tiny $\bullet$}}}
\put(20,32){\makebox(0,0){{\tiny $\bullet$}}}
\put(20,40){\makebox(0,0){{\tiny $\bullet$}}}
\put(20,48){\makebox(0,0){{\tiny $\bullet$}}}
\put(20,56){\makebox(0,0){{\tiny $\bullet$}}}
\put(20,64){\makebox(0,0){{\tiny $\bullet$}}}

\put(20,72){\makebox(0,0){{$\bullet$}}}
\put(20,82){\makebox(0,0){{$\bullet$}}}
\put(20,92){\makebox(0,0){{$\bullet$}}}
\put(20,102){\makebox(0,0){{$\bullet$}}}

\put(20,112){\makebox(0,0){{$\bullet$}}}
\put(20,122){\makebox(0,0){{$\bullet$}}}
\put(20,132){\makebox(0,0){{$\bullet$}}}
\put(20,142){\makebox(0,0){{$\bullet$}}}

\put(20,152){\makebox(0,0){{$\bullet$}}}

\put(20,156){\makebox(0,0){{\tiny $[\a_0]$}}}

\put(20,146){\makebox(0,0){{\tiny $[\a_5^{(1)}]$}}}
\put(20,136){\makebox(0,0){{\tiny $[\a_4^{(1)}]$}}}
\put(20,126){\makebox(0,0){{\tiny $[\a_3^{(1)}]$}}}
\put(20,116){\makebox(0,0){{\tiny $[\a_2^{(1)}]$}}}

\put(20,106){\makebox(0,0){{\tiny $[\a_5^{(2)}]$}}}
\put(20,96){\makebox(0,0){{\tiny $[\a_4^{(2)}]$}}}
\put(20,86){\makebox(0,0){{\tiny $[\a_3^{(2)}]$}}}
\put(20,76){\makebox(0,0){{\tiny $[\a_2^{(2)}]$}}}

\put(80,80){\makebox(0,0){{$\bullet$}}}
\put(-40,80){\makebox(0,0){{$\bullet$}}}

\path(-40,80)(20,72)
\path(-40,80)(20,82)
\path(-40,80)(20,92)
\path(-40,80)(20,102)
\path(-40,80)(20,112)
\path(-40,80)(20,122)
\path(-40,80)(20,132)
\path(-40,80)(20,142)
\path(-40,80)(20,152)

\dottedline(20,72)(80,80)
\dottedline(20,82)(80,80)
\dottedline(20,92)(80,80)
\dottedline(20,102)(80,80)
\dottedline(20,112)(80,80)
\dottedline(20,122)(80,80)
\dottedline(20,132)(80,80)
\dottedline(20,142)(80,80)
\dottedline(20,152)(80,80)

\dottedline(20,64)(-40,80)(20,56)
\dottedline(20,48)(-40,80)(20,40)
\dottedline(20,32)(-40,80)(20,24)
\dottedline(20,16)(-40,80)(20,6)
\dottedline(-40,80)(20,-4)

\path(20,64)(80,80)(20,56) \path(20,48)(80,80)(20,40)
\path(20,32)(80,80)(20,24) \path(20,16)(80,80)(20,6)
\path(80,80)(20,-4)

\put(85,75){\makebox(0,0){{\tiny $[\a_6^-]$}}}
\put(-45,75){\makebox(0,0){{\tiny $[\a_6^+]$}}}

\put(20,-6){\makebox(0,0){{\tiny $[\rho_9]$}}}
\put(20,4.8){\makebox(0,0){{\tiny $[\rho_8]$}}}
\put(20,14){\makebox(0,0){{\tiny $[\rho_7]$}}}
\put(20,22){\makebox(0,0){{\tiny $[\rho_6]$}}}
\put(20,30){\makebox(0,0){{\tiny $[\rho_5]$}}}
\put(20,38){\makebox(0,0){{\tiny $[\rho_4]$}}}
\put(20,45){\makebox(0,0){{\tiny $[\rho_3]$}}}
\put(20,53){\makebox(0,0){{\tiny $[\rho_2]$}}}
\put(20,61){\makebox(0,0){{\tiny $[\rho_1]$}}}

\end{picture}
\end{center}
\caption{$Z_{5}$, fusion graphs of $[\a_6^\pm]$, where $Z_5^2=2Z_5$}
\label{qsZ5}
\end{figure}

Next we compute the other 9 irreducible induced $M$-$M$ sectors. Using
Frobenius reciprocity and the homomorphism property of the $\a$-induction
$$\lan\a_6^+\a_6^-,\a_6^+\a_6^-\ran=\lan\a_6^+\a_6^+,\a_6^-\a_6^-\ran=
\sum_{\eta,\mu}N_{66}^\eta N_{66}^\mu \lan\a_\eta^+,\a_\mu^-\ran=
\sum_{\eta,\mu}N_{66}^\eta N_{66}^\mu Z_{\eta\mu}=9$$
we get another 9 automorphisms $[\rho_k]$, with $k=1,\dots,9$.
Furthermore, as above
$\lan\a_6^\pm\a_6^+\a_6^-,\a_6^\mp\ran=\sum_{\eta,\mu}N_{66}^\eta
N_{66}^\mu=9$ so the fusion graphs of $[\a_6^\pm]$ are displayed
in \fig{qsZ5}. In this figure, we use straight lines for the
fusion graph with $[\a_6^+]$ whereas dashed lines for that with
$[\a_6^-]$. Vertices from the neutral system $\MXMo$ are
marked with larger dots. The decomposition $Z_5^2 = 2Z_{5}$
is reflected by the layering of the full system as two orbits
$\MXMp$ and $\MXMp\a_{6}^{-}$. The full system has a
C$^\ast$-algebra structure (which is a weak C$^\ast$-Hopf algebra
of dimension 20, see \cite{PZ}) as follows
$$\hbox{Fusion}(\MXM)\simeq\bbC^4\oplus\hbox{Mat}_2
\oplus\hbox{Mat}_2\oplus\hbox{Mat}_2\oplus\hbox{Mat}_2.$$


The branching coefficient matrix is
$$b={\scriptstyle\pmatrix{1&1&0&0&0&0&0&0\cr 0&0&1&0&0&0&0&0\cr 
0&0&1&0&0&0&0&0\cr
0&0&0&1&0&0&0&0\cr 0&0&0&0&0&1&0&0\cr 0&0&0&0&1&0&0&0\cr
0&0&0&1&0&0&0&0\cr 0&0&0&0&1&0&0&0\cr 0&0&0&0&0&1&0&0}.}$$

For the local realistion when the dual canonical endomorphism
$\theta = \la_{0}\oplus \la_{1}$
we can compute the canonical endomorphism $\gamma$
 using the machinery developed in
\cite[Corollary 3.19]{BE3}. It is, as a sector,
$[\gamma]=[\mathrm{id}]\oplus[\alpha_2]$ since
$\lan\a_i^\pm,\gamma\ran=$ $\lan\theta\la_i,\la_0\ran=$
$\delta_{i,0}+\delta_{i,3}.$

The full system of $Z_7$ is the same as that for $Z_5$ (up to
permutation of sectors).

\subsubsection*{Case $Z_{(22)}$}

As in the proof and notation of \ref{ext}, the neutral system
$\MXMo=\{\a_0,\a_1, \a_6^{(2)}, \a_7^{(2)}\}$ is isomorphic to the
quantum double of $\bbZ_2$. The $\mathcal{X}^\pm$-chiral systems
are $\MXMpm=\{\a_0,\a_1, \a_6^{(2)},
\a_7^{(2)},\a_3^\pm,\a_6^{\pm(1)}\}$, which are actually
isomorphic to $\widehat{S_3}\times \bbZ_2$ as sectors.
>From both
chiral systems we have so far eight $M$-$M$ irreducible sectors.

Now $\lan\a_3^+\a_3^-,\a_3^+\a_3^-\ran=\sum N_{3,3}^\xi
N_{3,3}^\eta Z_{\xi,\eta}=2$ and moreover the dimension (of the
intertwiner spaces) of $\a_3^+\a_3^-$ with any sector of both
chiral systems vanishes. Thence the irreducible decomposition of
$\a_3^+\a_3^-$ provides us with two new $M$-$M$ irreducible
sectors. Now $\lan\a_3^+\a_6^-,\a_3^+\a_6^-\ran=3$ and
$\lan\a_3^+\a_6^-,\a_6^+\ran=\lan\a_3^+\a_6^-,\a_7^+\ran=1$ and
the intertwiner spaces of $\a_3^+\a_6^-$ with any other
irreducible $M$-$M$ sector vanishes. Therefore we have two more
$M$-$M$ irreducible sectors $[(\a_3^+\a_6^-)^{(i)}]$. So
$[\a_3^+\a_6^-]= [\a_6^{+(1)}]\oplus [(\a_3^+\a_6^-)^{(1)}]\oplus
[(\a_3^+\a_6^-)^{(2)}]$. The full system $\MXM$ decomposes into
two sheets according to $Z_{(22)}^{2} = Z_{(22)} + Z_{22}$, with six
elements each.

The full systems of $Z=Z_{(32)}, Z_{(33)}, Z_{(23)}$ are obtained
by permutations of those for $Z_{(22)}$.

\subsubsection*{Case $Z_{55}$}

For the sufferable modular invariant $Z_{55}$, the cardinality of
the full $M$-$M$ system is $\Tr(Z_{55} Z_{55}^t)=9$ whose fusion
rules are in turn commutative \cite[Corollary 6.9]{BEK1}, and with
global indices $\omega_\pm=6,\quad \omega_0=1.$ For the
$\mathcal{X}^\pm$-chiral systems, we use
$$\lan\a_\la^+,\a_\mu^+\ran=
\lan\a_\la^+\a_\mu^+, \a_0\ran= \sum_\eta N_{\la\mu}^\eta Z_{\eta
0}= N_{\la\mu}^0+N_{\la\mu}^3+N_{\la\mu}^6
=\lan\a_\la^-,\a_\mu^-\ran.$$ So computing we get that
$\MXMpm=\{\a_0,\a_1^\pm,\a_2^\pm\}$ with
$[\a_3^\pm]=[\a_0]\oplus[\a_1^\pm],\quad [\a_5^\pm]=$
$[\a_4^\pm]=[\a_2^\pm],$ $[\a_6^\pm]=[\a_0]\oplus[\a_2^\pm],
[\a_7^\pm]=[\a_1^\pm]\oplus[\a_2^\pm].$ Since $\omega_\pm=6$, and
$d_{\a_2^\pm}=2$ we easily conclude that the sectors of $\MXMpm$
are $\widehat{S_3}$.
Hence the sectors from
$\MXM$ are $\widehat{S_3}\times \widehat{S_3}$ identifying the
sectors of $\MXMp$ with $\widehat{S_3}\times\{{\bf 1}\}$ and
those from $\MXMm$ with $\{{\bf 1}\}\times\widehat{S_3}$.

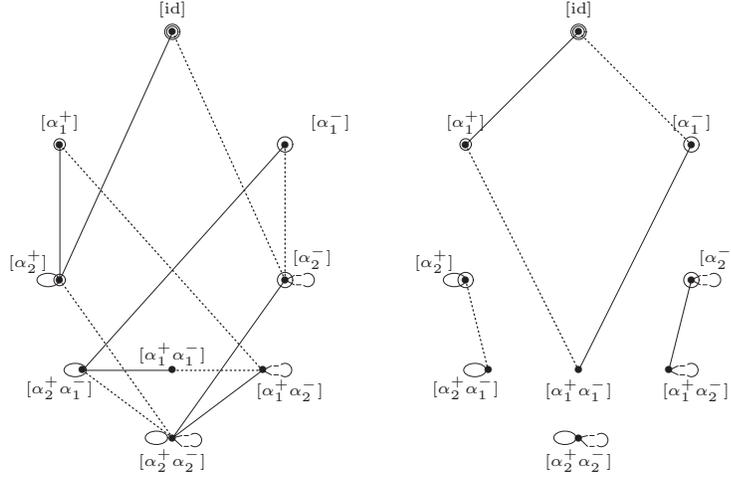
\begin{figure}[htb]
\begin{center}
\unitlength 0.6mm
\begin{picture}(75,110)
\thinlines

\put(-10,0){\makebox(0,0){{\tiny $\bullet$}}}
\put(-10,15){\makebox(0,0){{\tiny $\bullet$}}}
\put(-10,90){\makebox(0,0){{\tiny $\bullet$}}}

\put(10,15){\makebox(0,0){{\tiny $\bullet$}}}
\put(-30,15){\makebox(0,0){{\tiny $\bullet$}}}
\put(15,35){\makebox(0,0){{\tiny $\bullet$}}}
\put(-35,35){\makebox(0,0){{\tiny $\bullet$}}}

\put(-35,65){\makebox(0,0){{\tiny $\bullet$}}}
\put(15,65){\makebox(0,0){{\tiny $\bullet$}}}

\path(-35,65)(-35,35)(-10,90) \path(-10,15)(-30,15)(15,65)
\path(10,15)(-10,0)(15,35) \dottedline(-10,90)(15,35)
\dottedline(15,35)(15,65) \dottedline(-35,65)(10,15)
\dottedline(10,15)(-10,15) \dottedline(-35,35)(-10,0)
\dottedline(-10,0)(-30,15)

\put(-35,65){\arc{2.5}{0}{6.300}}
\put(-35,35){\arc{2.5}{0}{6.300}}
\put(-10,90){\arc{2.5}{0}{6.300}}
\put(-10,90){\arc{3.5}{0}{6.300}} \put(15,65){\arc{3.5}{0}{6.300}}
\put(15,35){\arc{3.5}{0}{6.300}}

\put(-37.5,35){\ellipse{5}{3}} \put(-31.5,15){\ellipse{5}{3}}%
\put(-13.5,0){\ellipse{5}{3}}

\path(-8,0.5)(-10,0)(-8,-2) \path(-7.5,0.5)(-6.5,0.5)
\path(-7.5,-2)(-6.5,-2) \put(-5,-0.5){\arc{3}{4.1}{8.5}}

\path(12,16)(10,15)(12,14) \path(12.5,16)(13.5,16)
\path(12.5,14)(13.5,14) \put(15,15){\arc{3}{4.1}{8.5}}

\path(17,36)(15,35)(17,34) \path(17.5,36)(18.5,36)
\path(17.5,34)(18.5,34) \put(20,35){\arc{3}{4.1}{8.5}}

\put(-10,-5){\makebox(0,0){{\tiny $[\a_2^+\a_2^-]$}}}
\put(-10,18.5){\makebox(0,0){{\tiny $[\a_1^+\a_1^-]$}}}
\put(-10,95){\makebox(0,0){{\tiny $[\mathrm{id}]$}}}

\put(16,10.5){\makebox(0,0){{\tiny $[\a_1^+\a_2^-]$}}}
\put(-35,11){\makebox(0,0){{\tiny $[\a_2^+\a_1^-]$}}}
\put(21,40){\makebox(0,0){{\tiny $[\a_2^-]$}}}
\put(-42,39){\makebox(0,0){{\tiny $[\a_2^+]$}}}

\put(-35,70){\makebox(0,0){{\tiny $[\a_1^+]$}}}
\put(25,70){\makebox(0,0){{\tiny $[\a_1^-]$}}}


\put(80,0){\makebox(0,0){{\tiny $\bullet$}}}
\put(80,15){\makebox(0,0){{\tiny $\bullet$}}}
\put(80,90){\makebox(0,0){{\tiny $\bullet$}}}

\put(100,15){\makebox(0,0){{\tiny $\bullet$}}}
\put(60,15){\makebox(0,0){{\tiny $\bullet$}}}
\put(105,35){\makebox(0,0){{\tiny $\bullet$}}}
\put(55,35){\makebox(0,0){{\tiny $\bullet$}}}

\put(55,65){\makebox(0,0){{\tiny $\bullet$}}}
\put(105,65){\makebox(0,0){{\tiny $\bullet$}}}

\path(80,90)(55,65) \path(105,65)(80,15) \path(100,15)(105,35)

\dottedline(80,90)(105,65) \dottedline(55,65)(80,15)
\dottedline(55,35)(60,15)

\put(55,65){\arc{2.5}{0}{6.300}} 
\put(80,90){\arc{2.5}{0}{6.300}} \put(80,90){\arc{3.5}{0}{6.300}}
\put(105,65){\arc{3.5}{0}{6.300}} \put(55,35){\arc{3.5}{0}{6.300}}
\put(105,35){\arc{3.5}{0}{6.300}}

\put(52.5,35){\ellipse{5}{3}} \put(57,15){\ellipse{5}{3}}
\put(77,0){\ellipse{5}{3}}

\path(82,0.5)(80,0)(82,-2) \path(82.5,0.5)(83.5,0.5)
\path(82.5,-2)(83.5,-2) \put(85,-0.5){\arc{3}{4.1}{8.5}}

\path(102,16)(100,15)(102,14) \path(102.5,16)(103.5,16)
\path(102.5,14)(103.5,14) \put(105,15){\arc{3}{4.1}{8.5}}

\path(107,36)(105,35)(107,34) \path(107.5,36)(108.5,36)
\path(107.5,34)(108.5,34) \put(110,35){\arc{3}{4.1}{8.5}}

\put(80,-5){\makebox(0,0){{\tiny $[\a_2^+\a_2^-]$}}}
\put(80,10.5){\makebox(0,0){{\tiny $[\a_1^+\a_1^-]$}}}
\put(80,95){\makebox(0,0){{\tiny $[\mathrm{id}]$}}}

\put(106,10.5){\makebox(0,0){{\tiny $[\a_1^+\a_2^-]$}}}
\put(55,11){\makebox(0,0){{\tiny $[\a_2^+\a_1^-]$}}}
\put(111,40){\makebox(0,0){{\tiny $[\a_2^-]$}}}
\put(48,39){\makebox(0,0){{\tiny $[\a_2^+]$}}}

\put(55,70){\makebox(0,0){{\tiny $[\a_1^+]$}}}
\put(105,70){\makebox(0,0){{\tiny $[\a_1^-]$}}}

\end{picture}
\end{center}
\caption{$Z_{55}$, fusion graphs of $[\a_2^\pm]$ and $[\a_1^\pm]$
where $Z_{55}^2 = 3Z_{55}$}
\label{qsZ55}
\end{figure}

We display the fusion the all $M$-$M$ system together with the
fusion graphs of both $[\a_2^\pm]$ on the LHS and $[\a_1^\pm]$ on
the RHS in \fig{qsZ55}. In this figure, we use straight lines for
the fusion graphs of $[\a_2^+]$ and $[\a_1^+]$ whereas dashed
lines for those of $[\a_2^-]$ and $[\a_1^-]$. We also encircled
the $\mathcal{X}^+$-chiral sectors with small circles and with
larger circles those for the $\mathcal{X}^-$-chiral system.
The full system decomposes into three sheets $\MXMp, \MXMp\a_1^-,
\MXMp\a_2^-$ according to $Z_{55}^2 = 3Z_{55}$.
The canonical sector of the fundamental inclusion $N\subset M$
is $[\gamma]=[\mathrm{id}]\oplus[\a_1^+\a_1^-]\oplus[\a_2^+\a_2^-]$
by \cite[Corollary 3.19]{BE3} because
$\lan\a_1^+\a_1^-,\gamma\ran=\lan\a_2^+\a_2^-,\gamma\ran=1.$

We also conclude that the full system of $Z_{44}, Z_{45}, Z_{54}$
$\MXM$ as sectors are $\widehat{S}_3\times\widehat{S}_3$.

\subsubsection*{Case $Z_{22}$}

The global indices of $\MXMpm$ and $\MXMo$ are, respectively,
$\omega_\pm=6$  and $\omega_0=1$. The dimension of the
intertwiner spaces in the $\mathcal{X}^+$-chiral system is
obtained by $\lan \a_\la^+,\a_\mu^+\ran=\lan
\a_\la^+\a_\mu^+,\a_0\ran=$
$\sum_\eta N_{\la\mu}^\eta Z_{\eta
0}=N_{\la\mu}^0+N_{\la\mu}^1+2N_{\la\mu}^2.$ We then find that
$\MXMp=\{\a_0, \a_3^{+(1)}, \a_3^{+(2)},$ $\a_6^{+(1)},
\a_6^{+(2)}, \a_6^{+(3)}\}$ with $[\a_0]= [\a_1^+]$,
$[\a_2^+]=2[\a_0]$, $[\a_3^+]=[\a_4^+]=$
$[\a_5^+]=[\a_3^{+(1)}]\oplus[\a_3^{+(2)}],$ $[\a_7^+]=[\a_6^+]=$
$[\a_6^{+(1)}]\oplus[\a_6^{+(2)}]\oplus[\a_6^{+(3)}].$ Since
$Z_{22}$ is symmetric we have $\MXMm$=
$\{\a_0,\a_3^{-(1)},\a_3^{-(2)},\a_6^{-(1)},
\a_6^{-(2)},\a_6^{-(3)}\}$. We can conclude that $\MXMpm$ is as
sectors  ${S_3}$. Hence by the global indices we have $\MXM$ as
sectors $S_3\times S_3$ which decomposes into six sheets
of $S_3$ according to $Z_{22}^2 = 6Z_{22}$.
We also have for $Z=Z_{33}, Z_{23},
Z_{32}$ that ${(\MXM)}_Z\simeq{(\MXM)}_{Z_{22}}$.

\subsubsection*{Case $Z_{25}$}

We easily find that $\MXM$ as sectors $S_3\times\widehat{S}_3$,
$\MXMp$ as sectors $S_3$, $\MXMm$  as sectors $\widehat{S}_3$ and
$\MXMo=\{\a_0\}$. For $Z=Z_{24}, Z_{35}, Z_{34}$
we get isomorphic full systems as for $Z_{25}$. Note that
$\MXMp$ and $\MXMm$ are not isomorphic -- indeed their cardinalities
are different and moreover one system is commutative
whilst the other is not. Thus the decompositions of the full
system according to $Z_{25}Z_{25}^{t} = 3Z_{22}$
and $Z_{25}^{t}Z_{{25}} = 6Z_{55}$ are  different.

\subsubsection*{Case $Z_{(44)}$}

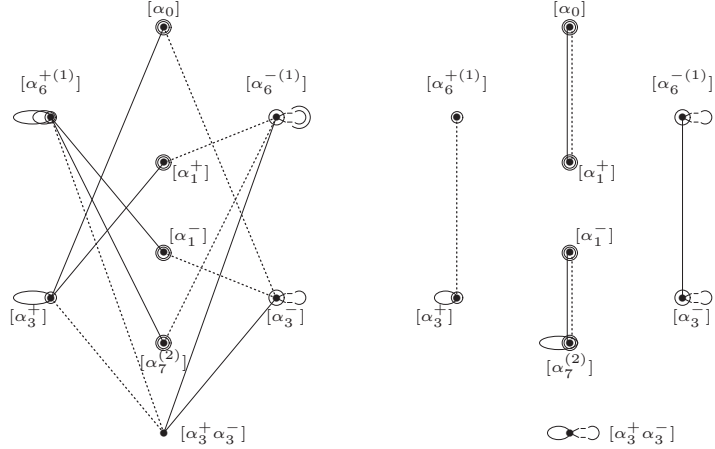
\begin{figure}[htb]
\begin{center}
\unitlength 0.6mm
\begin{picture}(70,110)
\thinlines

\put(-10,10){\makebox(0,0){\tiny $\bullet$}}
\put(-10,30){\makebox(0,0){\tiny $\bullet$}}
\put(-10,50){\makebox(0,0){\tiny $\bullet$}}
\put(-10,80){\makebox(0,0){\tiny $\bullet$}}

\put(15,20){\makebox(0,0){\tiny $\bullet$}}
\put(15,60){\makebox(0,0){\tiny $\bullet$}}
\put(-35,20){\makebox(0,0){\tiny $\bullet$}}
\put(-35,60){\makebox(0,0){\tiny $\bullet$}}
\put(-10,-10){\makebox(0,0){\tiny $\bullet$}}

\put(-10,10){\arc{2.5}{0}{6.300}}
\put(-10,30){\arc{2.5}{0}{6.300}}
\put(-10,50){\arc{2.5}{0}{6.300}}
\put(-10,80){\arc{2.5}{0}{6.300}}
\put(-35,20){\arc{2.5}{0}{6.300}}
\put(-35,60){\arc{2.5}{0}{6.300}}
\put(-10,10){\arc{3.5}{0}{6.300}}
\put(-10,30){\arc{3.5}{0}{6.300}}
\put(-10,50){\arc{3.5}{0}{6.300}}
\put(-10,80){\arc{3.5}{0}{6.300}} \put(15,20){\arc{3.5}{0}{6.300}}
\put(15,60){\arc{3.5}{0}{6.300}}

\put(-4.5,34){\makebox(0,0){{\tiny $[\a_1^-]$}}}
\put(-10,84){\makebox(0,0){{\tiny $[\a_0]$}}}
\put(-4,48){\makebox(0,0){{\tiny $[\a_1^+]$}}}
\put(-10,6){\makebox(0,0){{\tiny $[\a_7^{(2)}]$}}}
\put(-35,68){\makebox(0,0){{\tiny $[\a_6^{+(1)}]$}}}
\put(15,68){\makebox(0,0){{\tiny $[\a_6^{-(1)}]$}}}
\put(-40,16){\makebox(0,0){{\tiny $[\a_3^+]$}}}
\put(17,16){\makebox(0,0){{\tiny $[\a_3^-]$}}}
\put(1,-10){\makebox(0,0){{\tiny $[\a_3^+\a_3^-]$}}}

\put(-39,60){\ellipse{8.5}{3}} \put(-39,20){\ellipse{8.5}{3}}

\put(-36.5,60){\ellipse{5}{3}}%

\path(-10,80)(-35,20)(-10,50) \path(-10,30)(-35,60)(-10,10)
\path(15,60)(-10,-10)(15,20) \dottedline(15,20)(-10,80)
\dottedline(15,20)(-10,30) \dottedline(15,60)(-10,50)
\dottedline(15,60)(-10,10) \dottedline(-10,-10)(-35,60)
\dottedline(-10,-10)(-35,20)

\path(17,61)(15,60)(17,59) \path(17.5,61)(18.5,61)
\path(17.5,59)(18.5,59) \put(20,60){\arc{3}{4.1}{8.5}}
\put(20,60){\arc{5}{4}{8.5}}

\path(17,21)(15,20)(17,19) \path(17.5,21)(18.5,21)
\path(17.5,19)(18.5,19) \put(20,20){\arc{3}{4.1}{8.5}}


\put(80,10){\makebox(0,0){\tiny $\bullet$}}
\put(80,30){\makebox(0,0){\tiny $\bullet$}}
\put(80,50){\makebox(0,0){\tiny $\bullet$}}
\put(80,80){\makebox(0,0){\tiny $\bullet$}}

\put(105,20){\makebox(0,0){\tiny $\bullet$}}
\put(105,60){\makebox(0,0){\tiny $\bullet$}}
\put(55,20){\makebox(0,0){\tiny $\bullet$}}
\put(55,60){\makebox(0,0){\tiny $\bullet$}}
\put(80,-10){\makebox(0,0){\tiny $\bullet$}}

\path(79.5,80)(79.5,50) \path(79.5,10)(79.5,30)
\dottedline(80.5,80)(80.5,50) \dottedline(80.5,10)(80.5,30)
\dottedline(55,60)(55,20) \path(105,60)(105,20)

\put(80,10){\arc{2.5}{0}{6.300}} \put(80,30){\arc{2.5}{0}{6.300}}
\put(80,50){\arc{2.5}{0}{6.300}} \put(80,80){\arc{2.5}{0}{6.300}}
\put(55,20){\arc{2.5}{0}{6.300}} \put(55,60){\arc{2.5}{0}{6.300}}
\put(80,10){\arc{3.5}{0}{6.300}} \put(80,30){\arc{3.5}{0}{6.300}}
\put(80,50){\arc{3.5}{0}{6.300}} \put(80,80){\arc{3.5}{0}{6.300}}
\put(105,20){\arc{3.5}{0}{6.300}}
\put(105,60){\arc{3.5}{0}{6.300}}

\put(77.5,-10){\ellipse{5}{3}} \put(77.5,10){\ellipse{8.5}{3}}
\put(52.5,20){\ellipse{5}{3}}

\path(107,61)(105,60)(107,59) \path(107.5,61)(108.5,61)
\path(107.5,59)(108.5,59) \put(110,60){\arc{3}{4.1}{8.5}}

\path(107,21)(105,20)(107,19) \path(107.5,21)(108.5,21)
\path(107.5,19)(108.5,19) \put(110,20){\arc{3}{4.1}{8.5}}

\path(82,-9)(80,-10)(82,-11) \path(82.5,-9)(83.5,-9)
\path(82.5,-11)(83.5,-11) \put(85,-10){\arc{3}{4.1}{8.5}}

\put(85.5,34){\makebox(0,0){{\tiny $[\a_1^-]$}}}
\put(80,84){\makebox(0,0){{\tiny $[\a_0]$}}}
\put(86,48){\makebox(0,0){{\tiny $[\a_1^+]$}}}
\put(80,5){\makebox(0,0){{\tiny $[\a_7^{(2)}]$}}}
\put(55,68){\makebox(0,0){{\tiny $[\a_6^{+(1)}]$}}}
\put(105,68){\makebox(0,0){{\tiny $[\a_6^{-(1)}]$}}}
\put(50,16){\makebox(0,0){{\tiny $[\a_3^+]$}}}
\put(107,16){\makebox(0,0){{\tiny $[\a_3^-]$}}}
\put(96,-10){\makebox(0,0){{\tiny $[\a_3^+\a_3^-]$}}}

\end{picture}
\end{center}
\caption{$Z_{(44)}$: fusion graph of $[\a_{3}^\pm]$ and
$[\a_1^\pm]$ where $Z_{(44)}^2 = Z_{(22)} + Z_{44}$}\label{qsz(3)}
\end{figure}

The global indices of $\MXMpm$ and $\MXMo$ are $\omega_\pm=12$
and $\omega_0=4$ respectively. We compute the
$\mathcal{X}^\pm$-chiral systems and conclude that the
irreducible decompositions are as follows: $[\a_0]$,
$[\a_1^\pm]$, $[\a_2^\pm]= [\a_0]\oplus[\a_1^\pm]$,
$[\a_3^\pm]=[\a_4^\pm]=[\a_5^\pm]$,
$[\a_6^\pm]=[\a_6^{\pm(1)}]\oplus[\a_6^{\pm(2)}]$ and
$[\a_7^\pm]= [\a_6^{\pm(1)}]\oplus[\a_7^{\pm(2)}]$. By the
entries of the modular invariant $Z_{(44)}$, we see that
$\a_1^+=\a_6^{-(2)}$, $\a_6^{+(2)}=\a_1^-$ and
$\a_7^{+(2)}=\a_7^{-(2)}$ which it will be denoted simply by
$\a_7^{(2)}$. Thus $\MXMo=\{\a_0,\a_1^\pm, \a_7^{(2)}\}$ and
$\MXMpm=\MXMo\cup\{\a_3^\pm, \a_6^{\pm(1)}\}$. 
Since $\Tr(ZZ^t)=9$ we still have
to find one more irreducible $M$-$M$ sector in the full system,
which in turn is $\a_3^+\a_3^-$. In \fig{qsz(3)} we display the
fusion graphs of $[\a_3^\pm]$ and $[\a_1^\pm]$ whose details we
omit but are similar to those for e.g. $Z_5$. We use straight
lines for the fusion graph with $[\a_3^+]$ whereas dashed lines
for that with $[\a_3^-]$. Vertices from the neutral system
$\MXMo$ are marked with larger dots (similarly with the fusion
graph of $[\a_1^\pm]$ drawn on the RHS of \fig{qsz(3)}).
The decomposition $Z_{(44)}^2 = Z_{(22)} + Z_{44}$
means that the full system can be decomposed into a six element sheet
$\MXMp = \MXMo\cup\{\a_3^+, \a_6^{+(1)} \}$ from the chiral part of
$Z_{(22)}$  and a
further three element
sheet
$\{\a_3^{-}, \a_6^{-(1)},  \a_3^+\a_3^-\}$ from the chiral part of
$Z_{44}$.

The branching coefficient matrix is as in \erf{bb}.
Therefore by the proof of Proposition \ref{ext}, the extended
modular matrices $S^{\rm{ext}}$ and $T^{\rm{ext}}$ are the quantum
$\bbZ_2$ double model.

 The full system for $Z=Z_{(45)},
Z_{(55)}, Z_{(54)}$ are similar to that for $Z_{(44)}$.

\begin{remark}{\rm
It was the case with ${\mathit SU}(2)$ modular invariants (which
are all sufferable) that:
$$[\theta]_Z=\bigoplus_{\{\la: \hbox{ {\scriptsize FS}}_\la=1\}}
Z_{\la,\la} [\la]$$
is always a dual canonical sector producing the same $Z$, see \cite{E1}.
So since in the quantum $S_3$ double case, where the Frobenius-Schur
indicator FS$_\la=1$ for all $\la\in\NXN$, we would naturally ask whether
$\theta_Z$ is a dual endomorphism reproducing the sufferable modular
invariant $Z$. However in general this no longer holds true in our
current quantum $S_3$ double model as follows. Consider for example
the following sectors:\par
$[\theta_{Z_1}]=[\la_0]\oplus[\la_1]\oplus[\la_2]\oplus[\la_3]\oplus[\la_4]
\oplus[\la_5]\oplus[\la_6]\oplus[\la_7]$\par
$[\theta_{Z_2}]=[\la_0]\oplus[\la_1]\oplus[\la_4]\oplus[\la_5]
\oplus[\la_6]\oplus[\la_7]$\par
$[\theta_{Z_{55}}]=[\la_0]\oplus[\la_3]\oplus[\la_6]$\\
Amongst these, only
$\theta_{Z_{55}}$ is a dual canonical endomorphism producing
$Z_{55}$ by Proposition \ref{fundZ55}. The endomorphisms
$\theta_{Z_2}$  gives rise to inconsistent
dimensions of intertwiner spaces (since
$\lan\theta\la_4,\la_5\ran=0$ and $\lan\theta\la_0,\la_4\ran=
\lan\theta\la_0,\la_5\ran=1$),
whereas $\theta_{Z_1}$ would give a trace 6 or
trace 9 matrix instead of 8 if it reproduced $Z_1$.}
\label{canformZ}
\end{remark}


\section{The three quantum doubles of $D_5^{(1)}$}

There are precisely three subfactors whose principal graph is the
extended Dynkin diagram $D_5^{(1)}$. Izumi
\cite[Page 622]{iz2} has written down the modular data from the
Longo-Rehren inclusion of these subfactors. One coincides with
that from the quantum $S_3$ double and the other two lead to the
following  modular data :
{\footnotesize
\begin{eqnarray}
S&=&{1\over 6}\pmatrix{ 1&1&2&2&2&2&3&3\cr 1&1&2&2&2&2&-3&-3\cr
2&2&4\cos[4\pi/9]& 4\cos[8\pi/9]& 4\cos[2\pi/9]&-2&0&0\cr
2&2&4\cos[8\pi/9]& 4\cos[2\pi/9]& 4\cos[4\pi/9]&-2&0&0\cr
2&2&4\cos[2\pi/9]& 4\cos[4\pi/9]& 4\cos[8\pi/9]&-2&0&0\cr
2&2&-2&-2&-2&4&0&0\cr 3&-3&0&0&0&0&3&-3\cr
3&-3&0&0&0&0&-3&3}, \nonumber \\
T&=&\hbox{diag}(1,1,\varphi,\varphi^4,\varphi^7,1,1,-1) \nonumber
\end{eqnarray}
}
where $\varphi=\exp[\pm2\pi i/9]$.
All the $N$-$N$ sectors $[\la_i]$, $i=0,\dots,7$, are
self-conjugate and moreover the Frobenius-Schur indicator
FS$_\la=1$ for all $\la\in\NXN$.

\subsection{The modular invariants}

All the $N$-$N$ sectors $[\la_i]$, $i=0,\dots,7$, are
self-conjugate and moreover the Frobenius-Schur indicator
FS$_\la=1$ for all $\la\in\NXN$.

The dimension of the commutant $\{S,T\}^\prime$ is 6. With a
numerical search and employing the estimate of \cite[\eerf{1.3}]{BE5}
we find the list of modular invariants for the above model:
 \begin{eqnarray*}
Z_1&=& |\ch_0|^2+|\ch_1|^2+|\ch_2|^2+
|\ch_3|^2+|\ch_4|^2+|\ch_5|^2+|\ch_6|^2+|\ch_7|^2,\\
Z_2&=&|\ch_0+\ch_1|^2+2|\ch_2|^2+2|\ch_3|^2+2|\ch_4|^2+2|\ch_5|^2,\\
Z_3&=&|\ch_0+\ch_5|^2+(\ch_1+\ch_5)\ch_6^\ast+\ch_6(\ch_1+\ch_5)^\ast+
|\ch_7|^2,\\
Z_4&=&|\ch_0+\ch_5+\ch_6|^2,\\
Z_5&=&|\ch_0+\ch_5|^2+|\ch_1+\ch_5|^2+|\ch_6|^2+|\ch_7|^2,\\
Z_6&=&(\ch_0+\ch_5+\ch_6)(\ch_0+\ch_1+2\ch_5)^\ast,\\
Z_7&=&|\ch_0+\ch_1+\ch_5|^2+|\ch_2|^2+|\ch_3|^2+|\ch_4|^2+2|\ch_5|^2,\\
Z_8&=&(\ch_0+\ch_1+2\ch_5)(\ch_0+\ch_5+\ch_6)^\ast,\\
Z_9&=&|\ch_0+\ch_1+2\ch_5|^2.
\end{eqnarray*}

The fundamental dual canonical sector is
$[\theta]=[\la_0]\oplus[\la_5]\oplus[\la_6]$ producing the trace 3
matrix $Z_4$. Considering the two (even) endpoints of the graph
$D_5^{(1)}$ which both have quantum dimension equals 1 and
applying Izumi's Galois correspondence \cite{iz1} there exists an
intermediate subfactor $N\subset P\subset M$ with $N\subset M$
being the above subfactor that produces $Z_4$. Since the Jones
index of $N\subset P$ is 3 and its dual canonical sector, a
subsector of $[\theta]$ must be $[\theta_{N\subset
P}]=[\la_0]\oplus[\la_5]$. This endomorphism produces a type I
modular invariant, of trace 6, $Z_5$. The simple currents sector
$[\theta]=[\la_0]\oplus[\la_1]$ is a dual canonical sector. It
produces the trace 10 matrix $Z_2.$ Also $Z_8=Z_2Z_4$ is
sufferable by Theorem \ref{nice}, and so is its conjugate $Z_6$.
Then by the type I parent theorem \cite{BE4}, $Z_9$ is sufferable
as well. By the same type of computation used in Proposition
\ref{ext}, we can conclude that $Z_3$ is also sufferable. The
spectrum of $G_{[5]}$ is $\{-1^3,2^5 \}$ if a nimrep exists for
$Z_7$, but there is no graph with this spectrum (as was observed
by Gannon \cite{Gan2}). Hence we have the following classification
for the (normalized) modular invariants of
LR($D_5^{(1)},w=\exp(-2\pi i/3))$. The modular invariant $Z_7$ is
nimless, but the modular invariants $Z_1, Z_2$, $Z_3$, $Z_4$,
$Z_5$, $Z_6$, $Z_8$ and $Z_9$ are all sufferable (hence nimble).
We also have the fusion of sufferable modular invariants in Table
\ref{fusionsuffd51}. The decompositions in this table are unique.

\unitlength 1.4mm
\thinlines
\begin{table}[htb]
\begin{center}
{\scriptsize
\begin{tabular}{|c|c c c c c c c c|}  \hline\hline
&$Z_1$&$Z_2$&$Z_3$& $Z_4$&$Z_5$&$Z_6$&$Z_8$&$Z_9$\\ \hline\hline
$Z_1$& $Z_1$&$Z_2$&$Z_3$&$Z_4$&$Z_5$&$Z_8$&$Z_6$&$Z_9$ \cr
$Z_2$&$Z_2$&$2Z_2$&$Z_8$&$Z_8$&$Z_9$&$2Z_8$&$Z_9$&2$Z_9$ \cr

$Z_3$& $Z_3$&$Z_6$&$Z_4$+$Z_5$&$Z_4+Z_8$&$Z_3+Z_6$&3$Z_4$&$Z_6+Z_9$&$3Z_6$\\
$Z_4$&$Z_4$&$Z_6$&$Z_4+Z_6$& $3Z_9$&$Z_4+Z_6$&$3Z_9$&$3Z_6$&3$Z_6$\\
$Z_5$&$Z_5$&$Z_9$&$Z_3+Z_8$&$Z_4+Z_8$&$Z_5+Z_9$&3$Z_8$&$Z_6+Z_9$&3$Z_9$\\
$Z_6$&$Z_6$&2$Z_6$&3$Z_4$&$3Z_4$&3$Z_6$&6$Z_4$&3$Z_6$&6$Z_6$\\
$Z_8$&$Z_8$&$Z_9$&$Z_8+Z_9$&3$Z_8$&$Z_8+Z_9$&3$Z_8$&3$Z_9$&3$Z_9$\\
$Z_9$&$Z_9$&$2Z_9$&3$Z_8$&3$Z_8$&3$Z_9$&6$Z_8$&3$Z_9$&6$Z_9$
\\\hline
\end{tabular}
}
\end{center}
\caption{Fusion $Z_aZ_b^t$ of sufferable $D_5^{(1)}$ modular invariants}
\label{fusionsuffd51}
\end{table}

\vspace{0.2cm}\addtolength{\baselineskip}{-2pt}
\begin{footnotesize}
\noindent{\it Acknowledgement.}
The first named author is supported by the
EU QSNG network in Quantum Spaces - Noncommutative Geometry,
and the EPSRC network ABC--KLM, whilst the second named
author is supported by FCT (Portugal) under
grant BD/9704/96 and CAMGDS-IST. We are extremely grateful for their
financial assistance.
\end{footnotesize}


\newcommand\biba[7]   {\bibitem{#1} {#2:} {\sl #3.} {\rm #4} {\bf #5,}
                     {#6 } {#7}}
\newcommand\bibx[4]   {\bibitem{#1} {#2:} {\sl #3} {\rm #4}}

\def\ASENS            {Ann. Sci. \'Ec. Norm. Sup.}
\def\AM               {Acta Math.}
\def\AnM              {Ann. Math.}
\def\CMP              {Commun.\ Math.\ Phys.}
\def\IJM              {Internat.\ J. Math.}
\def\JAMS             {J. Amer. Math. Soc.}
\def\JFA              {J.\ Funct.\ Anal.}
\def\JMP              {J.\ Math.\ Phys.}
\def\JRA              {J. Reine Angew. Math.}
\def\LMP              {Lett.\ Math.\ Phys.}
\def\RMP              {Rev.\ Math.\ Phys.}
\def\RNM              {Res.\ Notes\ Math.}
\def\RIMS             {Publ.\ RIMS.\ Kyoto Univ.}
\def\Inv              {Invent.\ Math.}
\def\npbp             {Nucl.\ Phys.\ {\bf B} (Proc.\ Suppl.)}
\def\nupb             {Nucl.\ Phys.\ {\bf B}}
\def\adma             {Adv.\ Math.}
\def\coma             {Con\-temp.\ Math.}
\def\PAMS             {Proc. Amer. Math. Soc.}
\def\PJM              {Pacific J. Math.}
\def\ijmp             {Int.\ J.\ Mod.\ Phys.\ {\bf A}}
\def\jpa              {J.\ Phys.\ {\bf A}}
\def\PLB              {Phys.\ Lett.\ {\bf B}}
\def\RIMS             {Publ.\ RIMS, Kyoto Univ.}
\def\Top               {Topology}
\def\TAMS             {Trans.\ Amer.\ Math.\ Soc.}
\def\Duke              {Duke Math.\ J.}
\def\K                 {K-theory}
\def\JOP               {J.\ Oper.\ Theory}


\end{document}